\documentclass[11pt]{article}
\usepackage{amsthm}
\usepackage[toc,page]{appendix}

\usepackage{hslTR}

\usepackage{graphics}
\usepackage{graphicx}
\usepackage{amssymb}
\usepackage{verbatim}     
\usepackage{amsxtra}
\usepackage{epsfig}
\usepackage{hyperref} 
\usepackage{amsmath}
\usepackage{algpseudocode}
\usepackage{latexsym}
\usepackage{adjustbox}
\usepackage{ifthen}
\usepackage{psfrag}
\usepackage[caption=false,font=footnotesize]{subfig}
\usepackage{color}
\usepackage[normalem]{ulem}
\usepackage{float}
\usepackage{tikz}
\usetikzlibrary{positioning}
\usetikzlibrary{shapes,arrows}
\usepackage{cuted}
\usepackage{multicol}
\usepackage{bm}
\usepackage{mathtools}
\usepackage{etoolbox}
\usepackage[outdir=./]{epstopdf}
\usepackage[modulo,switch,columnwise]{lineno} 
\usepackage{changepage}
\usetikzlibrary{shapes,arrows}

\usepackage{cite}
\usepackage{amsmath,amssymb,amsfonts}
\usepackage{graphicx}
\usepackage{textcomp}

\usepackage{rgsEnvironments}
\usepackage{rgsMacros}  

\usetikzlibrary{shapes,arrows}

\newcommand{\Tmat}{\mathcal{T}}
\newcommand{\DG}{\mathcal{G}}
\newcommand{\nodes}{\mathcal{V}}
\newcommand{\edges}{\mathcal{E}}
\newcommand{\adj}{A}
\newcommand{\nbors}{\mathcal{N}}
\newcommand{\minus}{\scalebox{0.75}[1.0]{$-$}}

\newcommand{\EndNew}{\color{black}}

\newtheorem{problem}{Problem}[section]

\tikzstyle{block} = [draw, fill=blue!20, rectangle, 
    minimum height=3em, minimum width=6em]
\tikzstyle{sum} = [draw, fill=blue!20, circle, node distance=1cm]
\tikzstyle{input} = [coordinate]
\tikzstyle{output} = [coordinate]
\tikzstyle{pinstyle} = [pin edge={to-,thin,black}]

\newbool{rep}
\setbool{rep}{true}

\newbool{conf}
\setbool{conf}{false}

\newbool{two_col}
\setbool{two_col}{false}

\ifbool{rep}{
 \addtolength{\oddsidemargin}{-.875in}
	\addtolength{\evensidemargin}{-.875in}
	\addtolength{\textwidth}{1.75in}

	\addtolength{\topmargin}{-.875in}
	\addtolength{\textheight}{1.75in} 
}{}
\begin{document}

\ifbool{rep}{
\ititle{HyNTP: A Distributed Hybrid Algorithm for Time Synchronization}
\title{\bf HyNTP: A Distributed Hybrid Algorithm for Time Synchronization}
\iauthor{
  Marcello Guarro \\
  {\normalsize mguarro@ucsc.edu} \\
    Ricardo Sanfelice \\
  {\normalsize ricardo@ucsc.edu}}
\idate{\today{}} 
\iyear{2021}
\irefnr{13}
\makeititle}{
\title{\bf HyNTP: A Distributed Hybrid Algorithm for \\ Time Synchronization}

\author{Marcello Guarro and Ricardo G. Sanfelice
\thanks{This research has been partially supported by the National Science Foundation under Grant no. ECS-1710621, and Grant no. CNS-1544396, by the Air Force Office of Scientific Research under Grant no. FA9550-16-1-0015, Grant no. FA9550-19-1-0053, and Grant no. FA9550-19-1-0169, and by CITRIS and the Banatao Institute at the University of California.}
\thanks{Marcello Guarro and Ricardo G. Sanfelice are with Department of Electrical
and Computer Engineering, University of California, Santa Cruz. Email:
        {\tt\small mguarro@ucsc.edu,ricardo@ucsc.edu}}%
}
}

\maketitle

\begin{abstract}
This paper presents a distributed hybrid algorithm that synchronizes the time and rate of a set of clocks connected over a network. Clock measurements of the nodes are given at aperiodic time instants and the controller at each node uses these measurements to achieve synchronization. Due to the continuous and impulsive nature of the clocks and the network, we introduce a hybrid system model to effectively capture the dynamics of the system and the proposed hybrid algorithm. Moreover, the hybrid algorithm allows each agent to estimate the skew of its internal clock in order to allow for synchronization to a common timer rate. We provide sufficient conditions guaranteeing synchronization of the timers, exponentially fast, with robustness. Numerical results illustrate the synchronization property induced by the proposed algorithm as well as its performance against comparable algorithms from the literature.
\end{abstract}

\section{Introduction}

\subsection{Motivation}

Since the advent of asynchronous packet-based networks in communication and information technology, the topic of clock synchronization has received significant attention due to the temporal requirements of packet-based networks for the exchange of information. In recent years, as distributed packet-based networks have evolved in terms of size, complexity, and, above all, application scope, there has been a growing need for new clock synchronization schemes with tractable design conditions to meet the demands of these evolving networks. 

Distributed applications such as robotic swarms, automated manufacturing, and distributed optimization rely on precise time synchronization among distributed agents for their operation; see \cite{7}. For example, in the case of distributed control and estimation over networks, the uncertainties of packet-based network communication requires timestamping of sensor and actuator messages in order to synchronize the information to the evolution of the dynamical system being controlled or estimated. Such a scenario is impossible without the existence of a common timescale among the non-collocated agents in the system. In fact, the lack of a shared timescale among the networked agents can result in performance degradation that can destabilize the system; see \cite{graham2004clock}. Moreover, one cannot  always assume that consensus on time is a given, especially when the network associated to the distributed system is subject to perturbations such as noise, delay, or jitter. Hence, it is essential that these networked systems utilize clock synchronization schemes that establish and maintain a common timescale for their algorithms.

\subsection{Background and Related Work}

For many networked control system settings, each agent in the system is fitted with its own internal hardware clock and one or more software clocks that inherits the dynamics of the hardware clock. In an ideal scenario, the $i$-th agent in the system would have a clock $\tau_i \in \mathbb{R}_{\geq 0}$ such that $\tau_i(t) = t$, where $t$ is the global or real time. However, many hardware clocks utilize quartz-crystal or MEMS oscillators, susceptible to manufacturing imperfections and environmental factors that affect oscillator frequency; see \cite{wu2010clock} and \cite{vig1992introduction}. Due to the variability in oscillator frequency, one generally considers the continuous-time dynamics of the $i$-th hardware clock node given by 
\begin{equation} \label{eqn:clock_dyn}
\dot{\tau}_i = a_i
\end{equation} 
\noindent
where $a_i \in \reals$ defines the clock's rate of change. Solving the differential equation gives the following relationship to the ideal clock or real time reference $t$:
\begin{equation}
\tau_i(t) = a_i t + \tau_i(0)
\end{equation} 
\noindent
where the initial condition $\tau_i(0)$ gives the offset from $t$. For a network of $n$ agents, the notion of clock synchronization can be defined as the state of the networked system such that $\tau_i = \tau_j$ for all $i,j \in \{1,2,\ldots, n \}$, $i \neq j$. 

In an ideal setting with no delay and identical clock  rates (or skews), synchronization between two nodes, Node $1$ and Node $2$, can be achieved by the following simple reference-based algorithm. Node $1$ sends its time to Node $2$. Node $2$ calculates its offset relative to $1$. Node $2$ applies the offset correction to its clock. For the case of non-identical clock skews, a pair of measurements from Node $1$ would allow Node $2$ to calculate its relative skew $\frac{a_1}{a_2}$ and apply a correction accordingly. In a realistic setting, however, network communication between nodes is subject to a variety of delays to which such simple reference-based algorithms are nonrobust; see \cite{freris2010fundamental}. Moreover, these algorithms become cumbersome in terms of network utilization and computation as the number of nodes on the network increases.

The seminal Networking Time Protocol (NTP) presented in \cite{mills1991internet} mitigates these challenges through the implementation of a centralized algorithm.  In particular,  the networked agents  in the system  synchronize to a known reference that is either injected or provided by an  elected  leader agent. The effects of communication delay are mitigated via assumptions on the round-trip delay that occurs in the communication of any two nodes on the network. Conversely, other centralized approaches, such as those in \cite{8} and \cite{9}, assume the communication delay to be negligible and instead utilize least-squares minimization to estimate the errors in the offset and rates of change between the synchronizing nodes and the elected  reference  agent. Unfortunately, these approaches suffer robustness issues when communication with the reference node is lost or if the random delays in the transmission do not follow a normal distribution, see \cite{wu2010clock}. Moreover, algorithms like NTP were not designed for dynamic network topologies as they rely on predefined network hierarchies that define the relationships between the reference nodes and their children. Any change to the topology requires a reconstruction of the hierarchy adding considerable delay to the synchronization of the clocks.

Recently, the observed robustness issues in the centralized protocols have motivated leader-less, consensus-based approaches by leveraging the seminal results on networked consensus in \cite{5}, \cite{fax2004TAC}, and \cite{cao2008reaching}. In particular, the works of \cite{7}, \cite{10}, \cite{garone2015clock}, and, more recently, \cite{1}  employ  average consensus to give asymptotic results on clock synchronization  under asynchronous and asymmetric communication topology. However, due to the constraints on the communication modeling of the system the convergence results do not hold globally. Moreover, the lack of global convergence precludes any guarantees of robust asymptotic stability; see the forthcoming Remark 3.2. 

The work in \cite{carli2008pi} also considers a consensus-based approach by using a controller that uses a proportional gain to compensate for the clock rates and an integrator gain acting on an auxiliary control state that compensates for the clock offsets. Though the solution in \cite{carli2008pi} provides faster convergence than the other approaches using average consensus, the algorithm assumes periodic synchronous communication of the nodes. This assumption is relaxed in \cite{carli2014TAC} by considering asynchronous communication events. The authors in \cite{bolognani2015randomized} consider a similar relaxation but also relax assumptions on the graph structure.
However, in both  \cite{carli2014TAC} and \cite{bolognani2015randomized}  the clocks are slower to converge compared to the synchronous communication setting. Still, both synchronous and asynchronous scenarios require a large number of iterations before synchronization is achieved. Moreover, the algorithm subjects the clocks to significant nonsmooth adjustments in clock rate and offset that may prove undesirable in certain application settings or even prevent the rigorous establishment of robustness properties.

\ifbool{two_col}{}{
Other recent works include \cite{liao2013distributed} where the clocks are communicating over a network that is modeled as a Markov chain, the parameters of the clocks are then estimated from noisy measurements using a distributed algorithm. By making assumptions on the Markov chain, the estimates are shown to be mean square convergent if the union of the communication graphs is connected. The authors in \cite{chaudhari2008maximum} present an algorithm that estimates both skew and offset using a joint Maximum-Likelihood Estimation by assuming an exponential delay model.
}

\subsection{Contributions}

The lack of performance guarantees in the aforementioned works have motivated the design of a hybrid clock synchronization algorithm with tractable design conditions. In particular, this paper introduces a distributed hybrid algorithm that exponentially synchronizes a set of clocks connected over a network via measurements given at aperiodic time instants. 

Inspired by the contributions in \cite{3}, we present a distributed hybrid algorithm to synchronize the network clocks in the presence of non-ideal clock skews while capturing the continuous and impulsive dynamics of the network into a hybrid model. To achieve synchronization with a common rate of change, the algorithm also allows for local estimation of the skew of the internal clock at each agent. The use of a hybrid systems model to solve the problem under consideration allows for the application of a Lyapunov-based analysis to show stability of a desired set of interest.
Using results from \cite{4}, we show that, via a suitable change of coordinates, our distributed hybrid clock synchronization algorithm guarantees synchronization of the timers, exponentially fast, with robustness.
The main contributions of this paper are given as follows: 

\begin{itemize}
\item In Section \ref{sec:controller}, we introduce \textit{HyNTP}, a distributed hybrid algorithm that synchronizes the clock rates and offsets to solve the problem outlined in Section \ref{sec:prob}. Moreover, we present a hybrid systems model to capture the network dynamics for the case of synchronous and aperiodic communication events. In Section \ref{sec:cl_properties}, we present a reduced model of the system and a subsequent auxiliary model that is generated from  an appropriately defined  change of coordinates.
With the auxiliary model, we present necessary and sufficient conditions for which stability of a compact set, representing synchronization, holds. Moreover, we show that the system is robust to perturbations on the communication noise, clock drift, the desired clock rate reference, and to communication delays in Section \ref{sec:robustness}.
\ifbool{two_col}{}{\item In Section \ref{sec:ext}, we extend the algorithm to compensate for the case of asynchronous and aperiodic communication and introduce a new hybrid model to capture this setting. Similarly, to the synchronous setting, we consider a reduced model and then using a separate change of coordinates, generate an auxiliary model for which we present necessary and sufficient conditions for asymptotic stability of the system.}
\item In Section \ref{sec:num}, we compare the merits of our algorithm to competing algorithms in the literature.
\end{itemize}
We inform the reader that some details have been omitted due to space constraints and can be found in the technical report \cite{12}. This work is an extension of our conference paper \cite{225},
which pertains to the nominal case and has no proofs.

\ifbool{two_col}{
\subsection{Notation}

The set of natural numbers including zero, i.e., $\{0,1,2,\ldots\}$ is denoted by $\mathbb{N}$. The set of natural numbers is denoted as $\mathbb{N}_{> 0}$, i.e., $\mathbb{N}_{> 0} = \{1,2,\ldots\}$. The set of real numbers is denoted as $\reals$. The set of nonnegative real numbers is denoted by $\mathbb{R}_{\geq 0}$, i.e., $\mathbb{R}_{\geq 0} = [0, \infty )$. The $n$-dimensional Euclidean space is denoted $\mathbb{R}^n$. Given sets $A$ and $B$, $F: A \rightrightarrows B$ denotes a set-valued map from $A$ to $B$. For a matrix $A \in \mathbb{R}^{n \times m}$, $A^\mathsf{T}$ denotes the transpose of $A$. Given a vector $x \in \mathbb{R}^n$, $|x|$ denotes the Euclidean norm. Given a vector $x \in \mathbb{R}^n$ and  a nonempty  set $\Sigma \subset \reals^n$, $|x|_{\Sigma}$ denotes the Euclidean point-to-set distance, i.e., $|x|_{\Sigma} \doteq {\rm inf}_{y \in \Sigma} |x-y|$. Given two vectors $x \in \mathbb{R}^n$ and $y \in \mathbb{R}^m$, we use the equivalent notation $(x,y) = [ x^\mathsf{T} \hspace{2mm} y^\mathsf{T} \hspace{1mm} ]^\mathsf{T}$. Given a matrix $A \in \mathbb{R}^{n \times m}$, $|A| := \mbox{max} \{\sqrt{|\lambda|} : \lambda \in \mbox{eig}(A^\mathsf{T}A)\}$. For two symmetric matrices $A \in \mathbb{R}^{n \times m}$ and $B \in \mathbb{R}^{n \times m}$, $A \succ B$ means that $A - B$ is positive definite; conversely, $A \prec B$ means that $A - B$ is negative definite. A vector of $N$ ones is denoted $\textbf{1}_{N}$. The matrix $I_n$ is used to denote the identity matrix of size $n \times n$.
}{
\subsection{Organization and Notation}

This paper is organized as follows. Section \ref{sec:prelim} presents preliminary material on graph theory and hybrid systems. Section \ref{sec:prob} introduces the clock synchronization problem and the system being studied.  Section \ref{sec:controller}  outlines the algorithm under consideration, the associated hybrid model of the closed-loop system, including the main result. Section \ref{sec:cl_properties} gives properties to the nominal closed-loop system. Section \ref{sec:num} provides comparisons to other comparable algorithm through numerical examples.

\ifbool{conf}{\textit{Notation}: The set of natural numbers including zero, i.e., $\{0,1,2,\ldots\}$ is denoted by $\mathbb{N}$. The set of natural numbers is denoted as $\mathbb{N}_{> 0}$, i.e., $\mathbb{N}_{> 0} = \{1,2,\ldots\}$. The set of real numbers is denoted as $\reals$. The set of non-negative real numbers is denoted by $\mathbb{R}_{\geq 0}$, i.e., $\mathbb{R}_{\geq 0} = [0, \infty )$. The $n$-dimensional Euclidean space is denoted $\mathbb{R}^n$. Given sets $A$ and $B$, $F: A \rightrightarrows B$ denotes a set-valued map from $A$ to $B$. For a matrix $A \in \mathbb{R}^{n \times m}$, $A^\mathsf{T}$ denotes the transpose of $A$. Given a vector $x \in \mathbb{R}^n$, $|x|$ denotes the Euclidean norm. Given two vectors $x \in \mathbb{R}^n$ and $y \in \mathbb{R}^m$, we use the equivalent notation $(x,y) = [ x^\mathsf{T} \hspace{2mm} y^\mathsf{T} \hspace{1mm} ]^\mathsf{T}$. Given a matrix $A \in \mathbb{R}^n$, $|A| := \mbox{max} \{\sqrt{|\lambda|} : \lambda \in \mbox{eig}(A^\mathsf{T}A)\}$. For two symmetric matrices $A \in \mathbb{R}^n$ and $B \in \mathbb{R}^n$, $A \succ B$ means that $A - B$ is positive definite; conversely, $A \prec B$ means that $A - B$ is negative definite. A vector of $N$ ones is denoted $\textbf{1}_{N}$. The matrix $I_n$ is used to denote the identity matrix of size $n \times n$.}{
\textit{Notation}: The set of natural numbers including zero, i.e., $\{0,1,2,\ldots\}$ is denoted by $\mathbb{N}$. The set of natural numbers is denoted as $\mathbb{N}_{> 0}$, i.e., $\mathbb{N}_{> 0} = \{1,2,\ldots\}$. The set of real numbers is denoted as $\reals$. The set of non-negative real numbers is denoted by $\mathbb{R}_{\geq 0}$, i.e., $\mathbb{R}_{\geq 0} = [0, \infty )$. The $n$-dimensional Euclidean space is denoted $\mathbb{R}^n$. Given sets $A$ and $B$, $F: A \rightrightarrows B$ denotes a set-valued map from $A$ to $B$. For a matrix $A \in \mathbb{R}^{n \times m}$, $A^\mathsf{T}$ denotes the transpose of $A$. Given a vector $x \in \mathbb{R}^n$, $|x|$ denotes the Euclidean norm. Given a vector $x \in \mathbb{R}^n$ and set $\Sigma \subset \reals^n$, $|x|_{\Sigma}$ denotes the Euclidean point-to-set distance, i.e., $|x|_{\Sigma} \doteq {\rm inf}_{y \in \Sigma} |x-y|$. We denote $\mathbb{B}$ is a closed unit ball in $\reals^n$, and by $\A + \delta \mathbb{B}$ the set of all $x \in \reals^n$ such that $|x-a| \leq \delta$ for some $a \in \A$. Given two vectors $x \in \mathbb{R}^n$ and $y \in \mathbb{R}^m$, we use the equivalent notation $(x,y) = [ x^\mathsf{T} \hspace{2mm} y^\mathsf{T} \hspace{1mm} ]^\mathsf{T}$. Given a matrix $A \in \mathbb{R}^{n \times m}$, $|A| := \mbox{max} \{\sqrt{|\lambda|} : \lambda \in \mbox{eig}(A^\mathsf{T}A)\}$. For two symmetric matrices $A \in \mathbb{R}^{n \times m}$ and $B \in \mathbb{R}^{n \times m}$, $A \succ B$ means that $A - B$ is positive definite; conversely, $A \prec B$ means that $A - B$ is negative definite. A vector of $N$ ones is denoted $\textbf{1}_{N}$. The matrix $I_n$ is used to denote the identity matrix of size $n \times n$.}
}

\section{Preliminaries} \label{sec:prelim}

\subsection{Preliminaries on Graph Theory}

Let $\DG = (\nodes, \edges, \adj)$  be a weighted directed graph (digraph) where $\nodes = \{1,2,\ldots,n\}$ represents the set of $n$ nodes, $\edges \subset \nodes \times \nodes$ the set of edges, and $\adj \in \{0,1\}^{n \times n}$ represents the adjacency matrix. An edge of $\DG$ is denoted by $e_{ij} = (i,j)$. The elements of $\adj$ are denoted by $a_{ij}$, where $a_{ij} = 1$ if $e_{ij} \in \edges$ and $a_{ij} = 0$ otherwise. The in-degree and out-degree of a node $i$ are defined by $d^{in}(i) = \sum^n_{k=1} a_{ki}$ and $d^{out}(i) = \sum^n_{k=1} a_{ik}$, respectively. The largest and smallest in-degree of a digraph are given by $\bar{d} = \textmd{max}_{i \in \nodes} d^{in}(i)$ and $\underbar{$d$} = \textmd{min}_{i \in \nodes} d^{in}(i)$, respectively. The in-degree matrix is an $n \times n$ diagonal matrix, denoted $\mathcal{D}$, with elements given by
\ifbool{two_col}{ $d_{ij} = d^{in}(i) \mbox{ if } i = j$, $d_{ij} = 0 \mbox{ if } i \neq j$ for each $i \in \nodes$. }{
\begin{equation*}
d_{ij} = \begin{cases} d^{in}(i) & \mbox{ if } i = j \\ 0 & \mbox{ if } i \neq j \end{cases} \hspace{5mm} \forall i \in \nodes
\end{equation*}}
\noindent
The Laplacian matrix of a digraph $\DG$, denoted by $\mathcal{L}$, is defined as $\mathcal{L} = \mathcal{D} - \adj$ and has the property that $\mathcal{L} \textbf{1}_n = 0$. The set of nodes corresponding to the neighbors that share an edge with node $i$ is denoted by $\nbors(i) := \defset{k \in \nodes}{e_{ki} \in \edges}$. In the context of networks, $\nbors(i)$ represents the set of nodes for which an agent $i$ can communicate with.

\ifbool{two_col}{}{\begin{lemma} \label{lem:olfati}
((Olfati-Saber and Murray, 2004, Theorem 6),(Fax and Murray, 2004, Propositions 1, 3, and 4))
For an undirected graph, $\mathcal{L}$ is symmetric and positive semidefinite, and each eigenvalue of  $\mathcal{L}$ is real. For a directed graph, zero is a simple eigenvalue of  $\mathcal{L}$ if the directed graph is strongly connected.
\end{lemma} 

\begin{lemma} \label{lem:godsil}(Godsil and Royle (2001)) Consider an $n \times n$ symmetric matrix $A = \{a_{ik} \}$ satisfying $\sum_{i=1}^n a_{ik} = 0$ for each $k \in \nodes$. The following statements hold:
\begin{itemize}
\item There exists an orthogonal matrix $U$ such that $U^\top A U = \begin{bmatrix} 0 & 0 \\ 0 & \star \end{bmatrix}$, where $\star$ represents any nonsingular matrix with
appropriate dimensions and 0 represents any zero
matrix with appropriate dimensions.
\item The matrix $A$ has a zero eigenvalue with eigenvector
$\textbf{1}_n \in \reals^n$. 
\end{itemize}
\end{lemma}

\begin{definition} \label{def:consensus}
A weighted digraph is said to be
\begin{itemize}
\item \textit{balanced} {\rm if the in-degree matrix and the out-degree matrix for every node are equal, i.e., $d^{in}(i) = d^{out}(i)$ for each $i \in \nodes$.}
\item \textit{complete} {\rm if every pair of distinct nodes are connected by a unique edge, i.e., $a_{ik} = 1$ for each $i,k \in \nodes, i \neq k$.}
\item \textit{strongly connected} {\rm if and only if for any two distinct nodes there exists a path of directed edges that connects them.}
\end{itemize}
\end{definition}}

\subsection{Preliminaries on Hybrid Systems}

A hybrid system $\cal H$ in $\mathbb{R}^n$ is composed by the following \textit{data}: a set $C \subset \mathbb{R}^n$, called the flow set; a differential equation defined by the function $f: \mathbb{R}^n \to \mathbb{R}^n$ with $C \subset \mbox{dom } f$, called the flow map; a set $D \subset \mathbb{R}^n$, called the jump set; and a set-valued mapping $G: \mathbb{R}^n \rightrightarrows \mathbb{R}^n$ with $D \subset \mbox{dom } G$, called the jump map. Then, a hybrid system $\mathcal{H} := (C,f,D,G)$ is written in the compact form 
\begin{equation} \label{eqn:Hy}
\cal H : \begin{cases} x \in C \hspace{0.5cm} & \dot{x} \hspace{2mm} = f(x) \\
x \in D \hspace{0.5cm} &  x^+ \in G(x) 
\end{cases}
\end{equation}
\noindent
where $x$ is the system state. Solutions to hybrid systems are denoted by $\phi$ and are parameterized by $(t,j)$, where $t \in \mathbb{R}_{\geq 0}$ defines ordinary time and $j \in \mathbb{N}$ is a counter that defines the number of jumps. A solution $\phi$ is defined by a \textit{hybrid arc} on its domain $\mbox{dom } \phi$ with \textit{hybrid time domain} structure \cite{4}. The domain $\mbox{dom } \phi$ is a hybrid time domain if $\mbox{dom } \phi \subset \mathbb{R}_{\geq 0} \times \mathbb{N}$ and for each $(T,J) \in \mbox{dom } \phi,$ $\mbox{dom } \phi \cap ([0,T] \times \{0,1,...,J\})$ is of the form $\bigcup_{j = 0}^J ([t_j, t_{j+1}] \times \{j\})$, with $0 = t_0 \leq t_1 \leq t_2 \leq t_{J + 1}$. A function $\phi : \mbox{dom } \phi \to \reals^n$ is a \textit{hybrid arc} if $\mbox{dom } \phi$ is a hybrid time domain and if for each $j \in \mathbb{N}$, the function $t \mapsto \phi(t,j)$ is locally absolutely continuous on the interval $I^j = \{t : (t,j) \in \mbox{dom } \phi \}$. A solution $\phi$ satisfies the system dynamics; see \cite[Definition 2.6]{4} for more details. A solution $\phi$ is said to be \textit{maximal} if it cannot be extended by flow or a jump, and \textit{complete} if its domain is unbounded. The set of all maximal solutions to a hybrid system $\HS$ is denoted by $\mathcal{S}_{\HS}$ and the set of all maximal solutions to $\HS$ with initial condition belonging to a set $A$ is denoted by $\mathcal{S}_{\HS}(A)$. A hybrid system is \textit{well-posed} if it satisfies the hybrid basic conditions in \cite[Assumption 6.5]{4}. 

\ifbool{two_col}{}{
The following definition introduces stability notions for hybrid systems $\HS$ with state $x \in \reals^n$ to a closed set of interest $\A \subset \reals^n$.}

\begin{definition} \label{def:stability} Given a hybrid system $\HS$ defined on $\reals^n$, the closed set $\A \subset \reals^n$ is said to be \textit{globally exponentially stable} (GES) for $\HS$ if there exist $\kappa, \alpha > 0$ such that every maximal solution $\phi$ to $\HS$ is complete and satisfies $| \phi(t,j) |_{\A} \leq \kappa e^{\minus \alpha (t+ j)} |\phi(0,0)|_{\A}$ for each $(t,j) \in \mbox{dom } \phi$.
\end{definition}

\section{Problem Statement} \label{sec:prob}

Consider a group of $n$ sensor nodes connected over a network represented by a digraph $\DG = (\nodes,\edges, \adj)$. Two clocks are attached to each node $i$ of $\DG$: an (uncontrollable) internal clock $\tau^*_i \in \reals_{\geq 0}$ whose dynamics are given by 
\begin{equation} \label{eqn:clk_dynamics1}
\dot{\tau}^*_i = a_i
\end{equation} 
\noindent
and an adjustable clock $\tilde{\tau}_i \in \reals_{\geq 0}$ with dynamics
\begin{equation} \label{eqn:clk_dynamics2}
\begin{aligned}
\dot{\tilde{\tau}}_i & = a_i + u_i
\end{aligned}
\end{equation} 
\noindent
where $u_i \in \reals$ is a control input. In both of these models, the (unknown) constant $a_i$ represents the unknown drift of the internal clock. 

At times $t_j$ for $j \in \mathbb{N}_{> 0}$ (we assume $t_0 = 0$), node $i$ receives measurements $\tilde{\tau}_k$ from its neighbors, namely, for each $k \in \mathcal{N}(i)$. The resulting sequence of time instants $\{t_j\}^{\infty}_{j=1}$ is assumed to be strictly increasing and unbounded. Moreover, for such a sequence, the time elapsed between each time instant when the clock measurements are exchanged satisfies 
\ifbool{two_col}{
\begin{equation} \label{eqn:t_bounds}
T_1 \leq t_{j+1} - t_j \leq T_2, \hspace{3mm}
0 \leq t_1 \leq T_2 \hspace{5mm} \forall j \in \mathbb{N}_{> 0}
\end{equation}}{
\begin{equation} \label{eqn:t_bounds}
\begin{aligned}
T_1 & \leq t_{j+1} - t_j \leq T_2 \hspace{5mm} \forall j \in \mathbb{N}_{> 0} \\
0 & \leq t_1 \leq T_2
\end{aligned}
\end{equation}}
\noindent
where $0 < T_1 \leq T_2$, with $T_1$ defining a minimum time between consecutive measurements and $T_2$ defines the maximum allowable transfer interval (MATI).

\begin{remark}
The models for the clocks are based on the hardware and software relationship of the real-time system that implements them. That is, the internal clock $\tau_i^*$ is treated as a type of hardware oscillator while the adjustable clock $\tilde{\tau}_i$ is treated as a virtual clock, implemented in software (as part of the proposed algorithm), that evolves according to the dynamics of the hardware oscillator. Any virtual clock implemented in node $i$ inherits the drift parameter $a_i$ of the internal clock, which cannot be controlled. More importantly, this drift parameter is not known due to the fact that universal time information is not available to any node. The input $u_i$ is unconstrained as allowed by hardware platforms.
\end{remark}

\ifbool{rep}{
\begin{remark}
Note that the proposed communication strategy does not capture the situation when the agents in the network communicate at different times (asynchronous). The decision for modeling a synchronous communication strategy is due to the challenges in guaranteeing global, robust clock synchronization in an asynchronous setting.

For instance, in the problem settings of the work by \cite{7}, \cite{garone2015clock}, and \cite{stankovic2018distributed} the algorithms therein allow for the nodes to communicate at different time between nodes:
\begin{enumerate}
\item These articles consider $n$ networked nodes whose interconnections are represented by a undirected graph $\mathcal{G}$;
\item Each node $i$ is equipped with a local hardware clock $\tau_i$ and software clock $\hat{\tau}_i$;
\item Information is exchanged between an agent pair $(i,k)$ at times instants $t_j^{i,k}$.
\end{enumerate}
\noindent
However, the convergence properties guaranteed are neither global nor robust; see  \cite[Remark 3.2]{12}. In this paper, we demonstrate through the use of a simpler network model commonly used in the literature (see \cite{carli2008pi} and \cite{bolognani2015randomized}) and the use of hybrid systems tools that such properties are possible.

Let us point out that there are several algorithms in the literature that use the same network structure that we use.  These include, \cite{bolognani2015randomized}, \cite{carli2008pi}, and \cite{carli2014TAC}. In those works, they consider synchronous communication of the nodes in their algorithms which similarly exposes them to the issue of synchronization to a common sequence of time instants $t_j$.

On the other hand, in the problem settings of the work by \cite{7}, \cite{garone2015clock}, and \cite{stankovic2018distributed} the algorithms therein allow for the nodes to communicate at different time between nodes:
\begin{enumerate}
\item These articles consider $n$ networked nodes whose interconnections are represented by a undirected graph $\mathcal{G}$;
\item Each node $i$ is equipped with a local hardware clock $\tau_i$ and software clock $\hat{\tau}_i$;
\item Information is exchanged between an agent pair $(i,k)$ at times instants $t_j^{i,k}$.
\end{enumerate}
However, the convergence properties guaranteed are neither global not robust, unlike 
the results in our submission, where attaining such properties is possible through the use of a simpler network model and the use of hybrid systems tools.
In fact, consider the algorithm in \cite{7}, where  the communication event condition is strictly between nodes $i$ and $k$ at time instants\footnote{Interestingly, the authors in \cite{stankovic2018distributed}, consider these time instants $t_j^{i,k}$ to be driven by a ``global virtual communication clock" (as defined therein) that ticks according to each $t_j^{i,k}$ for each nodes $i$ and $k$ with $i \neq k$.} $t_j^{i,k}$.
The algorithm therein is given as follows, at times $t_j^i$, node $i$ broadcasts its time $\hat{\tau}_i(t_j^i)$ to its neighbors $\mathcal{N}(i)$. Upon receipt of the timestamp by nodes $k$ in the set of neighboring nodes $\mathcal{N}(i)$, each node updates its clocks and clock rates as follows:

\begin{equation}
\begin{cases} 
\hat{\tau}_k(t_{j+1}^i) = \nu_k(t_{j}^i) \tau_k(t_{j}^i) + \hat{o}_k(t_{j}^i) \\
\nu_k(t_{j+1}^i) = \rho_v \nu_k(t_{j}^i) + (1 - \rho_v) \eta_{ik} \nu_k \\
\eta_{ki} (t_{j+1}^i) = \rho_{\eta} \eta_{ki} (t_{j}^i) + (1 - \rho_{\eta}) \frac{\tau_{k}^{new}  - \tau_{k}^{old}}{\tau_{ik}^{new}  - \tau_{ik}^{old} } \\
\hat{o}_k(t_{j+1}^i) = \hat{o}_k(t_{j}^i)  + (1 - \rho_o) \big ( \hat{\tau}_i (t_{j+1}^i) - \hat{\tau}_k(t_{j+1}^i) \big )
\end{cases} \hspace{5mm} \forall k \in \mathcal{N}(i)
\end{equation}
\noindent
where $\tau_{k}^{new}=\tau_i(t_{j}^i)$, $\tau_{k}^{old}=\tau_i(t_{j-1}^i)$, $\tau_{ik}^{new}=\tau_k(t_{j}^i)$ $\tau_{ik}^{old}=\tau_k(t_{j-1}^i)$ and  $\rho_v \in (0,1)$, $\rho_{\eta} \in (0, 1)$, and $\rho_o \in (0,1)$. Note that since the quantities $\tau_{ik}^{new}$ and $\tau_{ik}^{old}$ involved in the update law for $\eta_{ij}$ in \cite[(11)]{7}  are part of the state of the controller, as the first line of (11) indicates. Hence, if at any communication time,  it is the case that $\tau_{ik}^{new}$ and $\tau_{ik}^{old}$ are very close to each other, then the update law for $\eta_{ik}$ would be rather large in norm.  In fact, the norm of the update law for $\eta_{ik}$ becomes arbitrarily large with $\tau_{ik}^{new} - \tau_{ik}^{old}$ decreasing (not necessarily equal). As a consequence,  the following issues arise:
\begin{itemize}
    \item Lack of globality: The initial condition for $\tau_{ij}^{new}$ and $\tau_{ij}^{old}$ needs to be properly chosen, which precludes a global result -- note that the limiting property established by \cite[Theorem 4]{7} is said to hold for initial condition for $\eta_{ij}$, but the initial conditions for the other state components is not specified.
    \item Lack of general robustness: The update law is not necessarily bounded.  From \cite[Theorem 7.21]{4}, we know that in order to have a general robust asymptotic stability property, the update law needs to be bounded.  Due to this, general robustness cannot be achieved.
\end{itemize}
\end{remark}
}{
\begin{remark}
Note that the proposed communication strategy does not capture the situation when the agents in the network communicate at different times (asynchronous). The decision for modeling a synchronous communication strategy is due to the challenges in guaranteeing global, robust clock synchronization in an asynchronous setting.

For instance, in the problem settings of the work by \cite{7}, \cite{garone2015clock}, and \cite{stankovic2018distributed} the algorithms therein allow for the nodes to communicate at different time between nodes:
\begin{enumerate}
\item These articles consider $n$ networked nodes whose interconnections are represented by a undirected graph $\mathcal{G}$;
\item Each node $i$ is equipped with a local hardware clock $\tau_i$ and software clock $\hat{\tau}_i$;
\item Information is exchanged between an agent pair $(i,k)$ at times instants $t_j^{i,k}$.
\end{enumerate}
\noindent
However, the convergence properties guaranteed are neither global nor robust; see  \cite[Remark 3.2]{12}. In this paper, we demonstrate through the use of a simpler network model commonly used in the literature (see \cite{carli2008pi} and \cite{bolognani2015randomized}) and the use of hybrid systems tools that such properties are possible.
\end{remark}
}

Under such a setup, our goal is to design a distributed hybrid controller that, without knowledge of the drift parameter and of the communication times in advance, assigns the input $u_i$ to drive each clock $\tilde{\tau}_i$ to synchronization with every other clock $\tilde{\tau}_k$, with $\tilde{\tau}_k$ evolving at a common prespecified constant rate of change $\sigma^* > 0$ for each $k \in \nodes$. This problem is formally stated as follows:


\begin{problem} \label{prob:1}
Given a network of $n$ agents with dynamics as in (\ref{eqn:clk_dynamics1}) and (\ref{eqn:clk_dynamics2}) represented by a directed graph $\mathcal{G}$ and $\sigma^* > 0$, design a distributed hybrid controller that achieves the following properties when information $\tilde{\tau}_k$ for each $k \in \mathcal{N}(i)$ is received by node $i$ at times $t_j$ satisfying (\ref{eqn:t_bounds}):
\begin{itemize}
\item[i)] Global clock synchronization: for each initial condition, the components $\tilde{\tau}_1, \tilde{\tau}_2, \ldots, \tilde{\tau}_n$ of each complete solution to the system satisfy $$\lim_{t \to \infty} | \tilde{\tau}_i(t) - \tilde{\tau}_k(t)| = 0 \hspace{5mm} \forall i,k \in \mathcal{V}, i \neq k$$
\item[ii)] Common clock rate: for each initial condition, the components $\tilde{\tau}_1, \tilde{\tau}_2, \ldots, \tilde{\tau}_n$ of each complete solution to the system satisfy $$\lim_{t \to \infty}  |\dot{\tilde{\tau}}_i(t) - \sigma^*| = 0 \hspace{5mm} \forall i \in \mathcal{V}$$
\end{itemize}
\end{problem}

\section{Distributed Hybrid Controller for \\ Time Synchronization} \label{sec:controller}

We define the hybrid model that provides the framework and a solution to Problem \ref{prob:1}. First, since we are interested in the ability of the rate of each clock to synchronize to a constant rate $\sigma^*$, we propose the following change of coordinates: for each $i \in \nodes$, define $e_i := \tilde{\tau}_i - r$, where $r \in \reals_{\geq 0}$ is an auxiliary variable such that $\dot{r} = \sigma^*$. The state $r$ is only used for analysis. Then, the dynamics for $e_i$ are given by
\begin{equation}
\begin{aligned}
\dot{e}_i & = \dot{\tilde{\tau}}_i - \sigma^* \qquad \forall i \in \nodes \\ 
\end{aligned}
\end{equation} By making the appropriate substitutions, one has
\begin{equation} \label{eqn:e_dynamics}
\begin{aligned}
\dot{e}_i & = a_i + u_i - \sigma^* \qquad \forall i \in \nodes \\ 
\end{aligned}
\end{equation}
\noindent
To model the network dynamics for aperiodic communication events at $t_j$'s satisfying (\ref{eqn:t_bounds}), we consider a timer variable $\tau$ with hybrid dynamics
\begin{equation} \label{eqn:timer}
\begin{aligned}
& \dot{\tau} \hspace{2mm} = \minus 1 \hspace{5mm} \tau \in (0,T_2],
\hspace{5mm} \tau^+ \in [T_1,T_2] \hspace{5mm} \tau = 0 
\end{aligned}
\end{equation}
\noindent
This model is such that when $\tau = 0$, a communication event is triggered, and $\tau$ is reset to a point in the interval $[T_1,T_2]$ in order to preserve the bounds given in (\ref{eqn:t_bounds}); see \cite{97}. Note that $\tau$ is a global variable that models the network dynamics of the system triggering the communication events between the nodes. Moreover, information on $\tau$ is not available to the nodes. One can think of this mechanism as a type of network manager that governs the communication events of the system; see \cite{carli2008pi} and \cite{bolognani2015randomized} for similar network models.


The proposed hybrid algorithm assigns a value to $u_i$ so as to solve Problem \ref{prob:1}, which in the $e_i$ coordinates requires $e_i$ to converge to zero for each $i \in \nodes$. In fact, the algorithm implements two feedback laws: a distributed feedback law and a local feedback law. The distributed feedback law utilizes a control variable $\eta_i \in \reals$ that is impulsively updated at communication event times using both local and exchanged measurement information $\tilde{\tau}_k$. Specifically, it takes the form $$\eta_i^+ = \sum_{k \in \mathcal{N}(i)} K_i^{k}(\tilde{\tau}_i, \tilde{\tau}_k)$$ where $K_i^{k}(\tilde{\tau}_i, \tilde{\tau}_k) := -\gamma_i (e_i - e_k)$ with $\gamma_i > 0$. Between communication event times, $\eta_i$ evolves continuously.\ifbool{rep}{The local feedback strategy utilizes a continuous-time linear adaptive estimator with states $\hat{\tau}_i \in \reals$ and $\hat{a}_i \in \reals$ to estimate the drift $a_i$ of the internal clock; see forthcoming Remark \ref{rem:estim_ai}.}{The local feedback strategy utilizes a continuous-time linear adaptive estimator with states $\hat{\tau}_i \in \reals$ and $\hat{a}_i \in \reals$ to estimate the drift $a_i$ of the internal clock.\footnote{ In \cite{12} we demonstrate the need for such a strategy to estimate $a_i$ since first difference methods would not be viable in our problem setting.}} The estimate of the drift is then injected as feedback to compensate for the effect of $a_i$ on the evolution of $\tilde{\tau}_i$. Furthermore, the local feedback strategy injects $\sigma^*$ to attain the desired clock rate for $\tilde{\tau}_i$.

Inspired by the protocol in \cite[Protocol 4.1]{3}, the dynamics of the $i$-th hybrid controller are given by
\begin{equation} \label{protocol}
\begin{aligned}
& \begin{drcases}
\dot{u}_i & =  h_i \eta_i - \mu_i ( \hat{\tau}_i - \tau^*_i ), \hspace{1mm}
\dot{\eta}_i =  h_i \eta_i \\
\dot{\hat{a}}_i & = - \mu_i ( \hat{\tau}_i - \tau_i^* ), \hspace{1mm}
\dot{\hat{\tau}}_i = \hat{a}_i - (\hat{\tau}_i - \tau_i^*) \hspace{4.5mm}
\end{drcases} \hspace{2mm} \tau \in [0, T_2] \\
& \begin{drcases} u_i^+ & =  {\minus}\gamma_i \sum_{k \in \mathcal{N}(i)} ( \tilde{\tau}_i {-} \tilde{\tau}_k) {-} \hat{a}_i {+} \sigma^*, \hspace{1mm} \hat{a}^+_i = \hat{a}_i \\
\eta_i^+ & =  -\gamma_i \sum_{k \in \mathcal{N}(i)} (\tilde{\tau}_i - \tilde{\tau}_k), \hspace{1mm}
\hat{\tau}^+_i = \hat{\tau}_i
\end{drcases} \hspace{2mm} \tau = 0
\end{aligned}
\end{equation}
\noindent
where $h_i \in \reals$, $\gamma_i > 0$ are controller parameters for the distributed hybrid consensus controller and $\mu_i > 0$ is a parameter for the local parameter estimator. The state $\eta_i$ is an auxiliary controller state that is injected into the control input $u_i$. At communication events $t_j$, i.e. when $\tau = 0$, both $u_i$ and $\eta_i$ reset to new values using measurement information $\tilde{\tau}_k$ from the neighbors $k \in \mathcal{N}(i)$ of node $i$. Moreover, the values for $\hat{a}_i$ and $\sigma^*$ that are injected as feedback into $u_i$ are kept constant across jumps. \EndNew
 
\ifbool{rep}{
\begin{remark} \label{rem:estim_ai}
The need for a linear adaptive estimation strategy to estimate $a_i$ arises from the fact that the reference time or universal time $t$ is inaccessible to the nodes. The lack of a known reference prevents the use of first difference strategies such as $$a_i= \frac{(\tau_i^*(t_j)- \tau_i^*(t_{j-1}))}{(t_j-t_{j-1})}$$ 
Moreover, if such a strategy were to be viable, it would then result in the possibility of not only estimating $a_i$ incorrectly but also resulting in unbounded solutions to the system. To demonstrate the need for an alternate strategy to estimate $a_i$, suppose that information on the global timer $\tau$ is available to all of the nodes. Then, each node could integrate $\tau$ to give access to the instants $t_j$. However, since calculating $a_i$ only requires the time delta between successive $t_j$'s, the dynamics of the global timer could be adjusted as follows
\begin{equation*} \label{eqn:timer}
\begin{aligned}
& \dot{\tau} \hspace{2mm} = 1 \hspace{5mm} \tau \in [0,T_2], \\
& \tau^+ = 0 \hspace{5mm} \tau \in [T_1,T_2]
\end{aligned}
\end{equation*}
\noindent
such that $a_i$ can be calculated via
\begin{equation} \label{protocol}
\begin{aligned}
& \begin{drcases}
\dot{\hat{a}}_i & = 0 \\ 
\dot{\ell}_i & = 0 \hspace{11mm}
\end{drcases} \hspace{2mm} \tau \in [0, T_2] \\ \\
& \begin{drcases} 
\hat{a}^+_i & = \frac{\tau_i^* - \ell_{i}}{\tau} \\
\ell_i^+ & = \tau_i^*
\end{drcases} \hspace{2mm} \tau = 0
\end{aligned}
\end{equation}
\noindent
where $\tau_i^*$ is the current local hardware clock state and $\ell_i$ is a local memory state to store the previous hardware clock state. 

Now, in simulating such a system, if the memory buffer is not properly initialized, that is $\ell_i(0) \neq \tau_i^*(0)$, then $\hat{a}_i^+$ will not yield $a_i$. Moreover, if $\tau = 0$, then the solution to $a_i$ would observe finite escape time due to a possible singularity. This observed singularity precludes the results from applying globally. Moreover, any robustness properties would not be possible since the strategy would violate the condition of boundedness on the jump map.

Note that a similar strategy is employed in \cite{7} and \cite{garone2015clock} rendering the relevant results on the convergence of the clock drift, to not apply globally.
\end{remark}
\EndNew
}{} 
 
With the timer variable and hybrid controller defined in (\ref{protocol}), we construct the hybrid closed-loop system $\HS$ obtained from the interconnection between the distributed hybrid controller and the local adaptive estimator given in error coordinates. The state of the closed-loop system is
\ifbool{two_col}{ 
\begin{equation} \label{eqn:state}
\begin{aligned}
x & = (e, u, \eta, \tau^*, \hat{a}, \hat{\tau}, \tau) \in \mathcal{X}
\end{aligned}
\end{equation}
where $\mathcal{X} := \reals^n \times \reals^n \times \reals^n \times \reals_{\geq 0}^n \times \reals^n \times \reals_{\geq 0}^n \times [0,T_2]$ with }{
\begin{equation} \label{eqn:state}
x = (e, u, \eta, \tau^*, \hat{a}, \hat{\tau}, \tau) \in \reals^n \times \reals^n \times \reals^n \times \reals_{\geq 0}^n \times \reals^n \times \reals_{\geq 0}^n \times [0,T_2] =: \mathcal{X}
\end{equation}
\noindent
where}
$e = (e_1, e_2, \ldots, e_n)$, $u = (u_1, u_2, \ldots u_n)$, $\eta = (\eta_1, \eta_2, \ldots, \eta_n)$, $\tau^* = (\tau^*_1, \tau^*_2, \ldots, \tau^*_N)$, $\hat{\tau} = (\hat{\tau}_1, \hat{\tau}_2, \ldots, \hat{\tau}_N)$, $a = (a_1, a_2, \ldots, a_N)$, and $\hat{a} = (\hat{a}_1, \hat{a}_2, \ldots, \hat{a}_n)$. 
\ifbool{two_col}{Then, let  
\begin{equation} \label{eqn:cl_hysys} 
\HS := (C,f,D,G) 
\end{equation} where the dynamics and data $(C,f,D,G)$ are given by $(\dot{e}, \dot{u}, \dot{\eta}, \dot{\tau}^*, \dot{\hat{a}}, \dot{\hat{\tau}}, \dot{\tau} )$ $= (a + u - \sigma^* \textbf{1}_n, h \eta -\mu ( \hat{\tau} - \tau^* ), h \eta, a, -\mu ( \hat{\tau} - \tau^* ), \hat{a} - (\hat{\tau} - \tau^*), -1) =: f(x)$ for each $x \in C$ and
$(e^+, u^+, \eta^+, {\tau^*}^+, \hat{a}^+, \hat{\tau}^+, \tau^+)$ $= (e, -\gamma \mathcal{L} e - \hat{a} + \sigma^* \textbf{1}_n, -\gamma \mathcal{L} e, \tau^*, \hat{a}, \hat{\tau}, [T_1,T_2]) =: G(x)$ for each $x \in D$
}{The dynamics and data $(C,f,D,G)$ of $\HS$ are given by
\begin{equation} \label{eqn:cl_hysys}
\hspace{-1mm}
\begin{bmatrix}
\dot{e} \\
\dot{u} \\
\dot{\eta} \\
\dot{\tau}^* \\
\dot{\hat{a}} \\
\dot{\hat{\tau}} \\
\dot{\tau}
\end{bmatrix} = \begin{bmatrix}
a + u - \sigma^* \textbf{1}_n  \\ h \eta -\mu ( \hat{\tau} - \tau^* ) \\ h \eta \\ a \\ -\mu ( \hat{\tau} - \tau^* ) \\ \hat{a} - (\hat{\tau} - \tau^*) \\ -1 \\
\end{bmatrix} \hspace{1mm} =: f(x) \hspace{1mm}  x \in C,
\hspace{1mm} 
\begin{bmatrix}
e^+ \\
u^+ \\
\eta^+ \\
{\tau^*}^+ \\
\hat{a}^+ \\
\hat{\tau}^+ \\
\tau^+
\end{bmatrix} =
 \begin{bmatrix}
e \\ -\gamma \mathcal{L} e - \hat{a} + \sigma^* \textbf{1}_n  \\ -\gamma \mathcal{L} e \\ \tau^* \\ \hat{a} \\ \hat{\tau} \\ [T_1,T_2]
\end{bmatrix} \hspace{1mm} =: G(x) \hspace{1mm} \hspace{1mm} x \in D
\end{equation}
}
\noindent
where $C := \mathcal{X}$ and $D := \{ x \in \mathcal{X} : \tau = 0 \}$. Note that $\mathcal{X} \subset \reals^m$ where $m = 7n$.

With the hybrid system $\HS$ defined, the next two results establish existence of solutions to $\HS$ and that every maximal solution to $\HS$ is complete. In particular, we show that, through the satisfaction of some basic conditions on the hybrid system data, which is shown first, the system $\HS$ is well-posed and that each maximal solution to the system is defined for arbitrarily large $t + j$. \ifbool{two_col}{}{The next two lemmas hold for any choice of parameters $T_1$, $T_2$, $\sigma^*$, $h$, $\gamma$, $\mu$, and strongly connected digraph $\DG$.}


\ifbool{two_col}{
\begin{lemma} \label{lem:hbc} The hybrid system $\HS$ satisfies the hybrid basic conditions defined in \cite[Assumption 6.5]{4}.
\end{lemma}}{
\begin{lemma} \label{lem:hbc} The hybrid system $\HS$ satisfies the following conditions, defined in \cite[Assumption 6.5]{4} as the hybrid basic conditions. 
\begin{itemize}
\item[(A1)] $C$ and $D$ are closed subsets of $\reals^m$.
\item[(A2)] $f : \mathcal{X} \to \mathcal{X}$ is continuous and locally bounded relative to $C$ and $C \subset \mbox{dom }  f$;
\item[(A3)] $G : \reals^m \rightrightarrows \reals^m$ is outer semicontinuous and locally bounded relative to $D$, and $D \subset \mbox{dom } G$.
\end{itemize}
\end{lemma}}

\ifbool{two_col}{
}{
\noindent See the appendix for proof.
}


\begin{lemma} \label{lem:complete1}
For every $\xi \in C \cup D (= \mathcal{X})$, every maximal solution $\phi$ to $\HS$ with $\phi(0,0) = \xi$ is complete.
\end{lemma}

\ifbool{two_col}{
The properties given in these two lemmas are easily established from the information given in the data of $\HS$; see \cite{12} for full details on the proofs of these results.
}{\noindent See the appendix for proof.
}

With the hybrid closed-loop system $\HS$ in (\ref{eqn:cl_hysys}), the set to asymptotically stabilize so as to solve Problem \ref{prob:1} is 
\ifbool{two_col}{
\begin{equation} \label{set:A}
\begin{aligned} 
\A & := \{x \in \mathcal{X} : e_i = e_k, \eta_i = 0, \hat{a}_i = a_i, \hat{\tau}_i = \tau^*_i, \\
& \hspace{30mm} u_i = \eta_i - \hat{a}_i + \sigma^* \hspace{1mm} \forall i ,k \in \mathcal{V} \}
\end{aligned}
\end{equation}
}{
\begin{equation} \label{set:A}
\begin{aligned} 
\A {:=} \{x \in \mathcal{X} : e_i = e_k, \eta_i = 0, \hat{a}_i = a_i, \hat{\tau}_i = \tau^*_i, u_i = \eta_i - \hat{a}_i + \sigma^* \hspace{1mm} \forall i ,k \in \mathcal{V} \}
\end{aligned}
\end{equation}
}
\noindent
Note that $e_i = e_k$ and $\eta_i = 0$ for all $i,k \in \mathcal{V}$ imply synchronization of the clocks, meanwhile $\hat{a}_i = a_i$ and  $\tau^*_i = \hat{\tau}_i$ for all $i,k \in \mathcal{V}$ ensure no error in the estimation of the clock skew and that the internal and estimated clocks are synchronized, respectively. The inclusion of $u_i = - \hat{a}_i + \sigma^*$ in $\A$ ensures that, for each $i \in \nodes$, $e_i$ remains constant (at zero) so that $e_i$ does not leave the set $\A$. This property is captured in the following result using the notion of forward invariance of a set. \ifbool{conf}{It can be shown that the set $\A$ is forward invariant for the hybrid system $\HS$.}{}

\begin{remark} \label{rem:htds}
Given that each maximal solution $\phi$ to $\HS$ is complete, with the state variable $\tau$ acting as a timer for $\HS$, for every initial condition $\phi(0,0) \in C \cup D$ we can characterize the domain of each solution $\phi$ to $\HS$ as follows:
\begin{equation}
\mbox{\rm dom } \phi = \bigcup_{j \in \mathbb{N}} [t_j, t_{j+1}] \times \{j\}
\end{equation}
\noindent
with $t_0 = 0$ and $t_{j+1} - t_j$ as in (\ref{eqn:t_bounds}). Furthermore, the structure of the above hybrid time domain implies that for each $(t,j) \in \mbox{dom } \phi$ we have
\begin{equation} \label{eqn:t_upper_bound}
t \leq T_2(j+1)
\end{equation}
\end{remark}


\ifbool{conf}{}{
\begin{lemma} \label{lem:complete}
Given a strongly connected digraph $\DG$, the set $\A$ in (\ref{set:A}) is forward invariant for the hybrid system $\HS$, i.e., each maximal solution $\phi$ to $\HS$ with $\phi(0,0) \in \A$ is complete and satisfies $\phi(t,j) \in \A$ for each $(t,j) \in \mbox{dom } \phi$ (see \cite[Chapter 10]{220}).
\end{lemma}

\ifbool{two_col}{}{
\noindent See the appendix for proof.
}
}

With the definitions of the closed-loop system $\HS$ in (\ref{eqn:cl_hysys}) and the set of interest $\A$ in (\ref{set:A}) to asymptotically stabilize in order to solve Problem \ref{prob:1}, we introduce our main result showing global exponential stability of $\A$ to $\HS$. This result is established through an analysis of an auxiliary system $\widetilde{\HS}_{\varepsilon}$ presented in (\ref{eqn:H_vareps}) and its global exponential stability for the auxiliary set $\tilde{\A}_{\varepsilon}$ in (\ref{set:A_tilde}), the details of which, along with a proofs, can be found in Section \ref{sec:thm_proof}.


\ifbool{two_col}{\begin{theorem} \label{thrm1}
Given a strongly connected digraph $\DG$, if the parameters $T_2 \geq T_1 > 0$, $\mu > 0$, $h \in \reals$, and $\gamma > 0$, the positive definite matrices $P_1$, $P_2$, and $P_3$ are such that 
\small
\begin{equation} \label{eqn:lyap1}
P_2 A_{f_3} + A_{f_3}^{\top} P_2 \prec 0, \hspace{2mm} P_3 A_{f_4} + A_{f_4}^{\top} P_3 \prec 0
\end{equation}
\begin{equation} \label{cond:phil_1}
\hspace{-3mm} A^{\top}_{g_2} {\exp(A^{\top}_{f_2} \nu)} P_1 {\exp(A_{f_2} \nu)} A_{g_2} {-} P_1 \prec 0 \hspace{3mm} \forall \nu \in [T_1, T_2]
\end{equation}
\begin{equation} \label{thrm_cond1}
\Big | \exp \Big (\frac{\bar{\kappa}_1}{\alpha_2} \hspace{0.5mm} T_2 \Big ) \Big ( 1 - \frac{\bar{\kappa}_2}{\alpha_2} \Big ) \Big | < 1
\end{equation}
\normalsize
hold,  where $A_{f_2}$, $A_{g_2}$ are given in (\ref{eqn:matrices}) and 
\small
\begin{equation} \label{eqn:thm_consts}
\begin{aligned}
\hspace{-5mm} \bar{\kappa}_1 & = \max \Big \{ \frac{\kappa_1}{2 \epsilon}, \frac{\kappa_1 \epsilon}{2} \minus \beta_2  \Big \}, \hspace{2mm} \bar{\kappa}_2 = \min \{1 , \kappa_2 \} \\
\kappa_1 & = 2 \underset{\nu \in [0,T_2]}{\mathrm{max}} \big | \exp{(A^{\top}_{f_2} \nu )} P_1 \exp{(A_{f_2} \nu )} \big | \\
\kappa_2 & \in \hspace{-1mm} \Big ( 0, {\minus} \hspace{-3mm} \underset{\nu {\in} [T_1,T_2]}{\min} \big \{ \lambda_{min} (A^{\top}_{g_2} {\exp{(A^{\top}_{f_2} \nu )}} P_1 {\exp{(A_{f_2} \nu )}} A_{g_2} {\minus} P_1 ) \big \} \Big ) \\
\alpha_2 & = \underset{\nu \in [0,T_2]}{\mathrm{max}} \Big \{ \exp{(2h \nu)}, \lambda_{max} \big (\exp{(A^{\top}_{f_2} \nu )} P_1 \exp{(A_{f_2} \nu )} \big ), \\
& \hspace{50mm} \lambda_{max}(P_2), \lambda_{max}(P_3) \Big \}
\end{aligned}
\end{equation}
\normalsize
with $\epsilon > 0$, and $\beta_1 > 0$ and $\beta_2 >0$ such that, in light of (\ref{eqn:lyap1}), 
$P_2 A_{f_3} + A_{f_3}^{\top} P_2 \preceq - \beta_1 I_2$ and $P_3 A_{f_4} + A_{f_4}^{\top} P_3 \preceq - \beta_2 I_{2(n-1)}$
\noindent
then, the set $\A$ in (\ref{set:A}) is globally exponentially stable for the hybrid system $\HS$ in (\ref{eqn:cl_hysys}).
\end{theorem}
}{
\begin{theorem} \label{thrm1}
Given a strongly connected digraph $\DG$, if the parameters $T_2 \geq T_1 > 0$, $\mu > 0$, $h \in \reals$, and $\gamma > 0$, the positive definite matrices $P_1$, $P_2$, and $P_3$ are such that \begin{equation} \label{eqn:lyap1}
P_2 A_{f_3} + A_{f_3}^{\top} P_2 \prec 0
\end{equation}
\begin{equation} \label{eqn:lyap2}
P_3 A_{f_4} + A_{f_4}^{\top} P_3 \prec 0
\end{equation}
\begin{equation} \label{cond:phil_1}
\hspace{-3mm} A^{\top}_{g_2} {\exp(A^{\top}_{f_2} \nu)} P_1 {\exp(A_{f_2} \nu)} A_{g_2} {-} P_1 \prec 0 \hspace{3mm} \forall \nu \in [T_1, T_2]
\end{equation}
\begin{equation} \label{thrm_cond1}
\Big | \exp \Big (\frac{\bar{\kappa}_1}{\alpha_2} \hspace{0.5mm} T_2 \Big ) \Big ( 1 - \frac{\bar{\kappa}_2}{\alpha_2} \Big ) \Big | < 1
\end{equation}
hold,  where $A_{f_2}$, $A_{g_2}$ are given in (\ref{eqn:matrices}) and 
\small
\begin{equation} \label{eqn:thm_consts}
\begin{aligned}
\hspace{-5mm} \bar{\kappa}_1 & {=} \max \Big \{ \frac{\kappa_1}{2 \epsilon}, \frac{\kappa_1 \epsilon}{2} \minus \beta_2  \Big \}, \hspace{2mm} \bar{\kappa}_2 = \min \{1 , \kappa_2 \} \\
\kappa_1 & {=} 2 \underset{\nu \in [0,T_2]}{\mathrm{max}} \big | \exp{(A^{\top}_{f_2} \nu )} P_1 \exp{(A_{f_2} \nu )} \big | \\
\kappa_2 & {\in} \big ( 0, {\minus} \hspace{-2mm} \underset{\nu {\in} [T_1,T_2]}{\min} \big \{ \lambda_{min} (A^{\top}_{g_2} {\exp{(A^{\top}_{f_2} \nu )}} P_1 {\exp{(A_{f_2} \nu )}} A_{g_2} {\minus} P_1 ) \big \} \big ) \\
\alpha_2 & {=} \underset{\nu \in [0,T_2]}{\mathrm{max}} \Big \{ \exp{(2h \nu)}, \lambda_{max} \big (\exp{(A^{\top}_{f_2} \nu )} P_1 \exp{(A_{f_2} \nu )} \big ), \\
& \hspace{50mm} \lambda_{max}(P_2), \lambda_{max}(P_3) \Big \}
\end{aligned}
\end{equation}
\normalsize
with $\epsilon > 0$, and $\beta_1 > 0$ and $\beta_2 >0$ such that, in light of (\ref{eqn:lyap1}), $P_2 A_{f_3} + A_{f_3}^{\top} P_2 \leq - \beta_1 I_2$, and $P_3 A_{f_4} + A_{f_4}^{\top} P_3 \leq - \beta_2 I_{2(n-1)}$ then, the set $\A$ in (\ref{set:A}) is globally exponentially stable for the hybrid system $\HS$ in (\ref{eqn:cl_hysys}).
\end{theorem}
}

To validate our theoretical stability result in Theorem \ref{thrm1}, consider five agents with dynamics as in (\ref{eqn:clk_dynamics1}) and (\ref{eqn:clk_dynamics2}) over a strongly connected digraph with the following adjacency matrix 
\ifbool{two_col}{\sloppy 
$\DG_A = ( [0, 1, 1, 0, 1],$ $[1, 0 , 1 , 0 , 0],$ $[1 , 0 , 0 , 1 , 0],$ $[0 , 0 , 1 , 0 , 1],$ $[1 , 0 , 1 , 1 , 0])$. }{
\small
\begin{equation*}
\DG_A = \begin{pmatrix}
0 & 1 & 1 & 0 & 1 \\
1 & 0 & 1 & 0 & 0 \\
1 & 0 & 0 & 1 & 0 \\
0 & 0 & 1 & 0 & 1 \\
1 & 0 & 1 & 1 & 0 \\
\end{pmatrix}
\end{equation*}}
\normalsize
Given $T_1 = 0.01$, $T_2 = 0.1$, and $\sigma^* = 1$, then it can be found that the parameters $h= -1.3$, $\mu = 3$, $\gamma = 0.125$, suitable matrices $P_1$, $P_2$, $P_3$ (see \cite{12} for details),
\ifbool{two_col}{}{ 
\begin{figure*} 
\begin{equation} \label{eqn:P_matrices} \scriptsize
\begin{aligned} 
P_1 & {=} \begin{bmatrix} 33.61  &  0 &  0  &  0  &  4.20  & 0  &  0  & 0 \\
       0 &  28.61  &  0  &  0 &  0 &   5.73  & 0  &  0 \\
      0  &  0  & 25.35 &  0 &   0 &   0  &  4.75 &  0 \\
       0  &  0 &  0  & 28.61 &   0 &  0  &  0 &   5.73 \\
       4.20  & 0 &   0  &  0 &   7.02 & 0   & 0  &  0 \\
      0  &  5.73 &   0  & 0 &  0 &  11.13  &  0 &  0 \\
       0  & 0 &   4.75  &  0 &   0 &   0 &  14.96 &  0 \\
      0  &  0 &  0  &  5.73 &   0 &  0  & 0 &  11.13] \\
\end{bmatrix} \\
P_2 & {=} \begin{bmatrix} 5.26  & \minus 2.24 \\ 
     \minus 2.24  &  7.54 \\
\end{bmatrix} \\
P_3 & {=} \begin{bmatrix} 6.47 &  0 &  0 &  0 &  \minus 2.36  &   0 &   0 &  0 \\
     0  &  6.47 &  0 &  0 &   0 &  \minus 2.36 &   0 &   0 \\
     0 &  0 &   6.47 &  0 &   0 &   0 &  \minus 2.36 &   0 \\
     0 &  0 &  0 &   6.47 &   0 &   0 &  0 &  \minus 2.36 \\
     \minus 2.35 &   0 &   0 &   0 &   7.90 &  0 &  0 &  0 \\
      0 &  \minus 2.35 &   0 &   0 &  0 &   7.90 &  0 &  0 \\
      0  &  0 &  \minus 2.35 &   0 &  0 &  0 &   7.90 &  0 \\
      0  &  0 &   0 &  \minus 2.35 &  0 &  0 &  0 &   7.90 \\
\end{bmatrix}
\end{aligned}
\end{equation}
\vspace{-10mm}
\end{figure*}}
\normalsize
\hspace{-3mm} and $\epsilon = 1.607$ \ifbool{conf}{with suitable matrices $P_1$, $P_2$, and $P_3$ }{}satisfy conditions (\ref{cond:phil_1}) and (\ref{thrm_cond1}) in Theorem \ref{thrm1} with $\bar{\kappa}_1 = 9.78$, $\kappa_1 = 31.44$, $\bar{\kappa}_2 = 1$, and $\alpha_2 = 18.923$. Figure \ref{fig:1} shows the trajectories of $e_i - e_k$, $\varepsilon_{a_i}$ for components $i \in \{1,2,3,4,5\}$ of a solution $\phi$ for the case where $\sigma = \sigma^*$ with initial conditions $\phi_{e}(0,0) = (1,-1,2,-2,0)$, $\phi_{\eta}(0,0) = (0,-3,1,-4,-1)$, and clock rates $a_i$ in the range $(0.85,1.15)$. \ifbool{two_col}{\hspace{-1mm}}{The bottom plot in Figure \ref{fig:1} depicts the Lyapunov trajectory $V$ evaluated along the solution $\phi$ with the upper bound given in (\ref{thrm1_bound}) projected onto the regular time domain. Observe that the exponential bound provided in (\ref{thrm1_bound}) jumps along the solution, validating our theoretical results on the exponential stability of the system.}\ifbool{rep}{\footnote{Code at github.com/HybridSystemsLab/HybridClockSync}}{\footnote{See \cite{12} for more simulations under different scenarios including a larger simulation with $N = 100$ nodes. Code at github.com/HybridSystemsLab/HybridClockSync}}

\begin{remark}  Theorem 4.5 not only assures that the proposed algorithm guarantees global exponential stability of the set $\A$ defined in (\ref{set:A}), but also that such a property is robust to perturbations — see Section VI for details.  It should be noted that the property is global in all of its variables, in the sense that regardless of the initial condition for the state $x$ of $\HS$, in particular, convergence (in distance) to $\A$ is assured.  These properties are not evident in other algorithms in the literature — in particular, the initial conditions for the variables $\tau_{new}$ and $\tau_{old}$ in \cite[Theorem 1]{7} need to be properly chosen to avoid unboundedness of the update law.  
\end{remark}

\begin{remark}
Observe that condition (\ref{cond:phil_1}) may be difficult to satisfy numerically as it may not be convex in $\gamma$ and $P_1$. The authors in \cite{97} utilize a polytopic embedding strategy to arrive at a linear matrix inequality in which one needs to find some matrices $X_i$ such that the exponential matrix is an element in the convex hull of the $X_i$ matrices. Such an algorithm can be adapted to our setting.
\end{remark}

\ifbool{rep}{
\begin{figure}
\centering
\includegraphics[trim={0mm 0mm 0mm 0mm},clip,width=1 \textwidth]{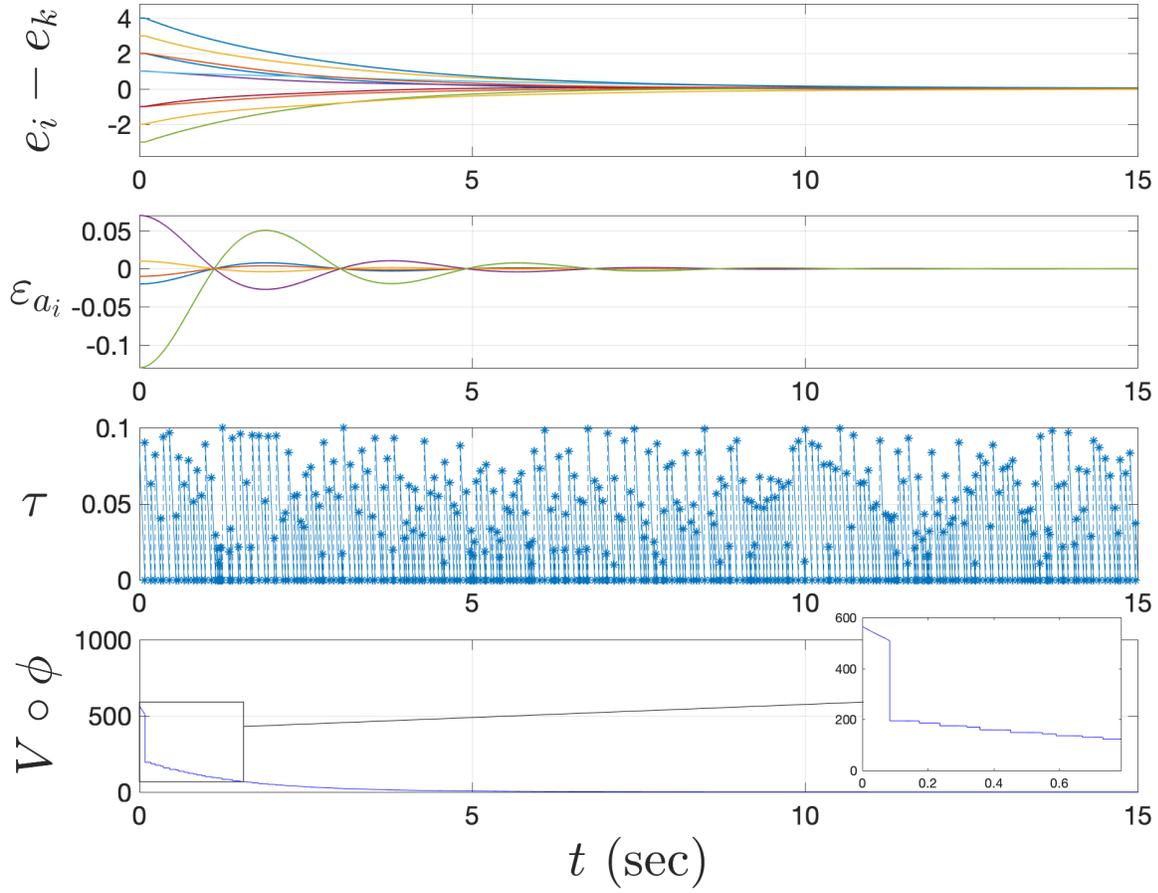}
\caption{\label{fig:1} (top) The trajectories of the state component errors $e_i - e_k$, $\varepsilon_{a_i}$, and $\tau$ for $i \in \{1,2,3,4,5\}$ of the solution $\phi$ for the case where $\sigma = \sigma^*$. (bottom) Plot of $V(\chi_{\varepsilon})$ evaluated along the solution $\phi$ with the associated bound (denoted $V_b$) given in (\ref{thrm1_bound}) projected onto the regular time domain.}
\end{figure}

\begin{figure}
\centering
\includegraphics[trim={0mm 0mm 0mm 0mm},clip,width=1 \textwidth]{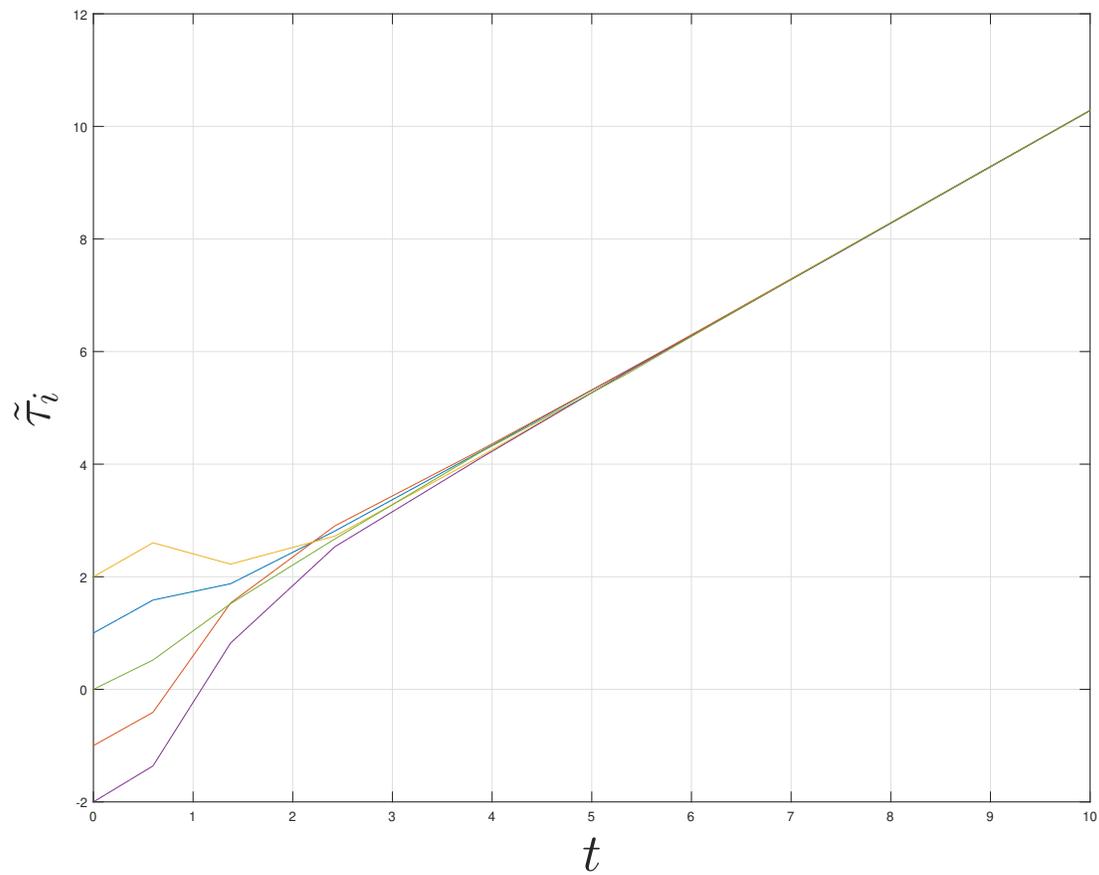}
\caption{\label{fig:2} The trajectories of the virtual clocks $\tilde{\tau}_i$ for $i \in \{1,2,3,4,5\}$ of the solution $\phi$ for the case where $\sigma = \sigma^*$.}
\end{figure}}{
\begin{figure}
\centering
\includegraphics[trim={0mm 7mm 0mm 5mm},clip,width=0.5 \textwidth]{../Figures/Multi_fig2_sigS_v10.eps}
\vspace{-0mm} \caption{\label{fig:1} The trajectories of the solution $\phi$ for state component errors $e_i {-} e_k$, $\varepsilon_{a_i}$, and $\tau$. Plot of $V$ evaluated along the solution $\phi$ projected onto the regular time domain (bottom).} \vspace{-0mm}
\end{figure}
}

\ifbool{two_col}{
\section{Key Properties of $\HS$ and \\ Proof of the Main Result} \label{sec:cl_properties}
}{
\section{Key Properties of the Nominal Closed-Loop System} \label{sec:cl_properties}}

\subsection{Reduced Model -- First Pass}

\color{black}

In this section, we recast the hybrid system $\HS$ into a reduced model obtained by setting $u = \eta  - \hat{a} + \sigma^*\textbf{1}_n$. This reduced model enables assessing asymptotic stability of $\A$. It is given in error coordinates for the parameter estimation of the internal clock rate and also the error of the internal clock state. We let $\varepsilon_a = a - \hat{a}$ denote the estimation error of the internal clock rate and $\varepsilon_{\tau} = \hat{\tau} - \tau^*$ represent the estimation error of the internal clock state. The state of the reduced model is given by $x_{\varepsilon} := (e,\eta, \varepsilon_a, \varepsilon_{\tau}, \tau) \in \reals^n \times \reals^n \times \reals^n \times \reals^n \times [0,T_2] =: \mathcal{X}_{\varepsilon}$ with dynamics defined by the data
\ifbool{two_col}{ 
\begin{equation} \label{eqn:hy_eps}
\begin{aligned}
\hspace{-2mm} f_{\varepsilon}(x_{\varepsilon}) & := (\eta + \varepsilon_a, h \eta, \mu \varepsilon_{\tau}, - \varepsilon_{\tau} - \varepsilon_a, -1) & \forall x_{\varepsilon} \in C_{\varepsilon}, \\
\hspace{-2mm} G_{\varepsilon}(x_{\varepsilon}) & := (e, \minus \gamma \mathcal{L} e, \varepsilon_a, \varepsilon_{\tau},  [T_1,T_2]) & \forall x_{\varepsilon} \in D_{\varepsilon}
\end{aligned}
\end{equation} 
\noindent
where 
}{
\begin{equation} \label{eqn:hy_eps}
\begin{aligned}
& f_{\varepsilon}(x_{\varepsilon}) := \begin{bmatrix}
\eta + \varepsilon_a  \\ h \eta \\ \mu \varepsilon_{\tau} \\ - \varepsilon_{\tau} - \varepsilon_a \\ -1 \\
\end{bmatrix} \hspace{1mm} \forall x_{\varepsilon} \in C_{\varepsilon},
& G_{\varepsilon}(x_{\varepsilon}) := \begin{bmatrix}
e \\ \minus \gamma \mathcal{L} e \\ \varepsilon_a \\ \varepsilon_{\tau} \\  [T_1,T_2]
\end{bmatrix} \hspace{1mm} \forall x_{\varepsilon} \in D_{\varepsilon}
\end{aligned}
\end{equation}
\noindent
where
}
$C_{\varepsilon} := \mathcal{X}_{\varepsilon}$ and $D_{\varepsilon} := \{ x_{\varepsilon} \in \mathcal{X}_{\varepsilon} : \tau = 0 \}$. This system is denoted $\HS_{\varepsilon} = (C_{\varepsilon},f_{\varepsilon},D_{\varepsilon},G_{\varepsilon})$. Note that the construction $u = \eta - \hat{a} + \sigma^* \textbf{1}_n$, which holds along all solutions after the first jump, leads to $\dot{e} = \eta + \varepsilon_a$. 

\ifbool{conf}{
Observe that the state of $\HS_{\varepsilon}$ utilizes most of the state components of $\HS$ except for the control input $u$ as it's value is captured by the states $\eta$ and $\varepsilon$ directly in the dynamics of $e$. 
}{To relate the properties of the reduced model to those of the hybrid system $\HS$, we establish a result showing an equivalency between the solutions of $\HS$ in (\ref{eqn:cl_hysys}) and $\HS_{\varepsilon}$ defined above. The result shows that after the first jump, each solution $\phi$ to $\HS$ is equivalent to a solution $\phi^{\varepsilon}$ to $\HS_{\varepsilon}$ when the trajectories of the timer variable $\tau$ for both solutions are equal. To facilitate such a result, we define the function $M: \mathcal{X} \to \mathcal{X}_{\varepsilon}$ given by
\begin{equation} \label{eqn:m}
M(x) := (e, \eta, a - \hat{a}, \hat{\tau} - \tau^*, \tau)
\end{equation}
\noindent
where $x = (e, u, \eta, \tau^*, \hat{a}, \hat{\tau}, \tau)$, as defined in (\ref{eqn:state}), and the function $\widetilde{M}: \mathcal{X}_{\varepsilon} \times \reals^n_{\geq 0} \times \reals^n_{\geq 0} \to \mathcal{X}$ given by
\ifbool{two_col}{
\begin{equation} \label{eqn:m}
\widetilde{M} (x_\varepsilon, \hat{\tau}, \tau^*) := (e, \eta - (a - \varepsilon_a) + \sigma^* \textbf{1}_n, \eta, \hat{\tau} - \varepsilon_{\tau}, a - \varepsilon_a, \varepsilon_{\tau} +  \tau^*, \tau)
\end{equation}
}{
\begin{equation} \label{eqn:m}
\widetilde{M} (x_\varepsilon, \hat{\tau}, \tau^*) := \begin{bmatrix}
e \\ \eta - (a - \varepsilon_a) + \sigma^* \textbf{1}_n \\ \eta \\ \hat{\tau} - \varepsilon_{\tau} \\ a - \varepsilon_a \\ \varepsilon_{\tau} +  \tau^* \\ \tau 
\end{bmatrix}
\end{equation}
}


\ifbool{two_col}{
\begin{lemma} \label{lem:equiv} 
Let $T_2 \geq T_1 > 0$, digraph $\DG$, and hybrid systems $\HS$ and $\HS_{\varepsilon}$ be given as in (\ref{eqn:cl_hysys}) and (\ref{eqn:hy_eps}), respectively. For each $\phi \in \mathcal{S}_{\HS}$ and each\footnote{Note that for a given solution $\phi^{\varepsilon}(t,j)$ to $\HS_{\varepsilon}$, the solution components are given by $\phi^{\varepsilon}(t,j) = \big (\phi_e^{\varepsilon}(t,j), \phi_\eta^{\varepsilon}(t,j), \phi^{\varepsilon}_{\varepsilon_a}(t,j), \phi^{\varepsilon}_{\varepsilon_\tau}(t,j), \phi^{\varepsilon}_{\tau}(t,j) \big )$.} $\phi^{\varepsilon} \in \mathcal{S}_{\HS_{\varepsilon}}$ such that $\phi (0,0) = \widetilde{M} \big (\phi^{\varepsilon}(0,0),  \phi_{\hat{\tau}}(0,0), \phi_{\tau^*}(0,0) \big )$ and timer components $\phi_{\tau}(t,j) = \phi^{\varepsilon}_{\tau}(t,j)$ for all $(t,j) \in \mbox{\rm dom } \phi$, it follows that $\mbox{dom } \phi = \mbox{dom } \phi^{\varepsilon}$ and $\phi (t,j) = \widetilde{M} \big (\phi^{\varepsilon}(t,j), \phi_{\hat{\tau}}(t,j), \phi_{\tau^*}(t,j) \big )$ for each $(t,j) \in \mbox{\rm dom } \phi$.
\end{lemma}
}{
\begin{lemma} \label{lem:equiv} 
Let $T_2 \geq T_1 > 0$, digraph $\DG$, and hybrid systems $\HS$ and $\HS_{\varepsilon}$ be given as in (\ref{eqn:cl_hysys}) and (\ref{eqn:hy_eps}), respectively. For each $\phi \in \mathcal{S}_{\HS}$ and each\footnote{Note that for a given solution $\phi^{\varepsilon}(t,j)$ to $\HS_{\varepsilon}$, the solution components are given by $$\phi^{\varepsilon}(t,j) = \big (\phi_e^{\varepsilon}(t,j), \phi_\eta^{\varepsilon}(t,j), \phi^{\varepsilon}_{\varepsilon_a}(t,j), \phi^{\varepsilon}_{\varepsilon_\tau}(t,j), \phi^{\varepsilon}_{\tau}(t,j) \big )$$} $\phi^{\varepsilon} \in \mathcal{S}_{\HS_{\varepsilon}}$ such that $\phi (0,0) = \widetilde{M} \big (\phi^{\varepsilon}(0,0),  \phi_{\hat{\tau}}(0,0), \phi_{\tau^*}(0,0) \big )$ and timer components $\phi_{\tau}(t,j) = \phi^{\varepsilon}_{\tau}(t,j)$ for all $(t,j) \in \mbox{\rm dom } \phi$, it follows that $\mbox{dom } \phi = \mbox{dom } \phi^{\varepsilon}$ and
\begin{equation} \label{eqn:equiv_sol_M}
\phi (t,j) = \widetilde{M} \big (\phi^{\varepsilon}(t,j), \phi_{\hat{\tau}}(t,j), \phi_{\tau^*}(t,j) \big ) \hspace{10mm} \forall (t,j) \in \mbox{\rm dom } \phi
\end{equation}
\end{lemma}

\noindent See the appendix for proof.
}
}

With the reduced model $\HS_{\varepsilon}$ in place, we consider the following set to globally exponentially stabilize for $\HS_{\varepsilon}$:
\begin{equation} \label{set:A_eps}
\A_{\varepsilon} {:=} \{ x_{\varepsilon} \in \mathcal{X}_{\varepsilon} : e_i {=} e_k, \eta_i {=} 0 \hspace{1mm} \forall i, k \in \mathcal{V}, \hspace{1pt} \varepsilon_{a} {=} 0, \varepsilon_{\tau} {=} 0 \}
\end{equation}
\noindent
This set is equivalent to $\A$ in the sense that the point-to-set distance metrics $|x|_{\A}$ and $|x|_{\A_{\varepsilon}}$ are equivalent when the map $\widetilde{M}$ is applied, as demonstrated in the results that follow. \ifbool{conf}{}{


\ifbool{two_col}{
\begin{lemma} \label{lem:set_dists}
Given sets $\A$ and $\A_{\varepsilon}$ as in (\ref{set:A}) and (\ref{set:A_eps}), respectively, for each $x =  (e, u, \eta, \tau^*, \hat{a}, \hat{\tau}, \tau)$, $x_{\varepsilon}$, $\hat{\tau}$, and $\tau^*$ such that $x \in \mathcal{X}$, $(x_{\varepsilon}, \hat{\tau}, \tau^* ) \in \mathcal{X}$, and $u = \eta - \hat{a} + \sigma^* \textbf{1}_n$ then $|x|_{\A} = |x_{\varepsilon}|_{\A_{\varepsilon}}$ and $|\widetilde{M}(x_{\varepsilon}, \hat{\tau}, \tau^*) |_{\A} = |x|_{\A}$.
\end{lemma}
}{
\begin{lemma} \label{lem:set_dists}
Given sets $\A$ and $\A_{\varepsilon}$ as in (\ref{set:A}) and (\ref{set:A_eps}), respectively, for each $x =  (e, u, \eta, \tau^*, \hat{a}, \hat{\tau}, \tau)$, $x_{\varepsilon}$, $\hat{\tau}$, and $\tau^*$ such that $x \in \mathcal{X}$, $(x_{\varepsilon}, \hat{\tau}, \tau^* ) \in \mathcal{X}$, and $u = \eta - \hat{a} + \sigma^* \textbf{1}_n$ then 
\begin{equation} \label{eqn:set_dists_1}
|x|_{\A} = |x_{\varepsilon}|_{\A_{\varepsilon}}
\end{equation} and
\begin{equation} \label{eqn:set_dists_2}
|\widetilde{M}(x_{\varepsilon}, \hat{\tau}, \tau^*) |_{\A} = |x|_{\A}
\end{equation}
\end{lemma}

With the stabilization set defined for $\HS_{\varepsilon}$, we have the following result that shows that if the set $\A_{\varepsilon}$ is globally exponentially stable for $\HS_{\varepsilon}$ then the set $\A$ is also globally exponentially stable for $\HS$.}


\begin{lemma} \label{lem:GES1}
Given $T_2 \geq T_1 > 0$ and a strongly connected digraph $\DG$, the set $\A$ in (\ref{set:A}) is GES for the hybrid system $\HS$ if $\A_{\varepsilon}$ in (\ref{set:A_eps}) is GES for the hybrid system $\HS_{\varepsilon}$.
\end{lemma}

\ifbool{two_col}{For proofs, see  \cite{12} and the Appendix.}
{
\noindent See the appendix for proof.
 }
 }

\subsection{Reduced Model -- Second Pass} 

Global exponential stability of $\A_{\varepsilon}$ for $\HS_{\varepsilon}$ is established by performing a Lyapunov analysis on a version of $\HS_{\varepsilon}$ obtained after an appropriate change of coordinates, one where the flow and jump dynamics are linearized. The model is obtained by exploiting an important property of the eigenvalues of the Laplacian matrix for strongly connected digraphs. 

\ifbool{two_col}{To this end, let $\DG$ be a strongly connected digraph. By exploiting certain known properties of strongly connected graphs as outlined in\cite[Lemma 2.1]{12} and \cite[Lemma 2.2]{12}, \EndNew one has that zero is a simple eigenvalue of the Laplacian matrix $\mathcal{L}$ with an associated eigenvector $v_1 = \frac{1}{\sqrt{N}} \textbf{1}_N$. Furthermore, there exists a nonsingular matrix}{To this end, let $\DG$ be a strongly connected digraph. By Lemma \ref{lem:olfati} and Lemma \ref{lem:godsil}, one has that zero is a simple eigenvalue of the Laplacian matrix $\mathcal{L}$ with an associated eigenvector $v_1 = \frac{1}{\sqrt{N}} \textbf{1}_N$. Furthermore, there exists a nonsingular matrix
}
\begin{equation} \label{eqn:T_mat}
\mathcal{T} = [v_1, \mathcal{T}_1]
\end{equation} 
where $\mathcal{T}_1 \in \reals^{N \times N \minus 1}$ is a matrix whose columns are the remaining eigenvectors of $\mathcal{L}$, i.e., $[v_2, \ldots, v_N]$, such that $\mathcal{T}^{\minus 1} \mathcal{L} \mathcal{T} = \begin{bmatrix} 0 & 0 \\ 0 & \bar{\mathcal{L}} \end{bmatrix}$, where $\mathcal{L}$ is the graph Laplacian of $\DG$ and $\bar{\mathcal{L}}$ is a diagonal matrix with the nonnegative eigenvalues of $\mathcal{L}$ as the diagonal elements given by $( \lambda_2, \lambda_3, \ldots, \lambda_N)$,  see \cite{5}, \cite{fax2004TAC}, and \cite{6} for more details.

To perform the said change of coordinates, we use $\Tmat$ to first perform the following transformations: $\bar{e} = \Tmat^{-1}e$, $\bar{\eta} = \Tmat^{-1}\eta$, $\bar{\varepsilon}_a = \Tmat^{-1}\varepsilon_a$ and $\bar{\varepsilon}_{\tau} = \Tmat^{-1}\varepsilon_{\tau}$. Then, we define vectors $\bar{z} = (\bar{z}_1, \bar{z}_2)$ and $\bar{w} = (\bar{w}_1, \bar{w}_2)$, where $\bar{z}_1 := (\bar{e}_1, \bar{\eta}_1)$, $\bar{z}_2 := (\bar{e}_2,\ldots,\bar{e}_N, \bar{\eta}_2, \ldots, \bar{\eta}_N)$, $\bar{w}_1 = (\bar{\varepsilon}_{a_1},\bar{\varepsilon}_{\tau_1})$, and $\bar{w}_2 = (\bar{\varepsilon}_{a_2},\ldots, \bar{\varepsilon}_{a_n},\bar{\varepsilon}_{\tau_2},\ldots, \bar{\varepsilon}_{\tau_n})$. Finally, we define $\chi_{\varepsilon} := (\bar{z}_1, \bar{z}_2, \bar{w}_1, \bar{w}_2, \tau) \in \reals^2 \times \reals^{2(n-1)} \times \reals^2 \times \reals^{2(n-1)} \times [0,T_2] =: \mathcal{X}_{\varepsilon}$ as the state of the new version of $\HS_{\varepsilon}$, which is denoted $\widetilde{\HS}_{\varepsilon}$ and has data given by
\ifbool{two_col}{
\begin{equation} \label{eqn:H_vareps}
\begin{aligned}
\widetilde{f}_{\varepsilon}(\chi_{\varepsilon}) {:=} \begin{bmatrix}
A_{f_1} \bar{z}_1 \\ A_{f_2} \bar{z}_2 \\ A_{f_3} \bar{w}_1 \\ A_{f_4} \bar{w}_2 \\  -1
\end{bmatrix} {+} \begin{bmatrix}
B_{f_1} \bar{w}_1 \\ B_{f_2} \bar{w}_2 \\ 0 \\ 0 \\ 0
\end{bmatrix}, \hspace{1mm}
\widetilde{G}_{\varepsilon}(\chi_{\varepsilon}) {:=} \begin{bmatrix}
A_{g_1} \bar{z}_1 \\ A_{g_2} \bar{z}_2 \\ \bar{w}_1 \\ \bar{w}_2 \\ [T_1,T_2]
\end{bmatrix}
\end{aligned}
\end{equation} 
}{
\begin{equation} \label{eqn:H_vareps}
\begin{aligned}
\widetilde{f}_{\varepsilon}(\chi_{\varepsilon}) {:=} \begin{bmatrix}
A_{f_1} \bar{z}_1 \\ A_{f_2} \bar{z}_2 \\ A_{f_3} \bar{w}_1 \\ A_{f_4} \bar{w}_2 \\  -1
\end{bmatrix} {+} \begin{bmatrix}
B_{f_1} \bar{w}_1 \\ B_{f_2} \bar{w}_2 \\ 0 \\ 0 \\ 0
\end{bmatrix} \hspace{2mm} \forall \chi_{\varepsilon} \in \tilde{C}_{\varepsilon}, \hspace{3mm}
\widetilde{G}_{\varepsilon}(\chi_{\varepsilon}) {:=} \begin{bmatrix}
A_{g_1} \bar{z}_1 \\ A_{g_2} \bar{z}_2 \\ \bar{w}_1 \\ \bar{w}_2 \\ [T_1,T_2]
\end{bmatrix} \hspace{2mm} \forall \chi_{\varepsilon} \in \tilde{D}_{\varepsilon}
\end{aligned}
\end{equation}
}
\normalsize
\noindent
 for each $\chi_{\varepsilon}$ in $\widetilde{C}_{\varepsilon} := \mathcal{X}_{\varepsilon}$ and in $\widetilde{D}_{\varepsilon} := \{ \chi_{\varepsilon} \in \mathcal{X}_{\varepsilon} : \tau = 0 \}$, respectively, with
\begin{equation} \label{eqn:matrices}
\begin{aligned}
& A_{f_1} {=} \begin{bmatrix}
0 & 1  \\ 0 & h
\end{bmatrix},
& & A_{f_2} {=} \begin{bmatrix} 0 & I_{m} \\ 0 & hI_{m}
\end{bmatrix}, 
& & A_{f_3} {=} \begin{bmatrix}
0 & \mu  \\ \minus 1 & \minus 1
\end{bmatrix} \\
& A_{f_4} {=} \begin{bmatrix} 0 & \mu I_{m} \\ \minus I_{m} & \minus I_{m}
\end{bmatrix}, 
& & B_{f_1} {=} \begin{bmatrix} 1 & 0 \\ 0 & 0
\end{bmatrix}, 
& & B_{f_2} {=} \begin{bmatrix} I_{m} & 0 \\ 0 & 0
\end{bmatrix} \\
& A_{g_1} {=} \begin{bmatrix} 1 & 0  \\ 0 & 0 
\end{bmatrix},
& & A_{g_2} {=} \begin{bmatrix} I_{m} & 0  \\ \minus \gamma \bar{\mathcal{L}} & 0 
\end{bmatrix}
\end{aligned}
\end{equation}
\normalsize
\noindent
and $m = N-1$. Then, $\widetilde{\HS}_{\varepsilon} = ( \widetilde{C}_{\varepsilon}, \widetilde{f}_{\varepsilon}, \widetilde{D}_{\varepsilon}, \widetilde{G}_{\varepsilon})$ denotes the new version of $\HS_{\varepsilon}$. The set $\A_{\varepsilon}$ to stabilize in the new coordinates for this hybrid system is given by 
\begin{equation} \label{set:A_tilde}
\begin{aligned}
\hspace{-1mm} \widetilde{\A}_{\varepsilon} & := \{ \chi_{\varepsilon} \in \mathcal{X}_{\varepsilon} : \bar{z}_1 {=} (e^* \hspace{-1mm} ,0), \bar{z}_2 {=} 0, \bar{w}_1 {=} 0, \bar{w}_2 {=} 0, e^* \hspace{-1mm} \in \reals \}
\end{aligned}
\end{equation}

\ifbool{conf}{
In the next result, we demonstrate how global exponential stability of $\tilde{\A}_{\varepsilon}$ for $\widetilde{\HS}_{\varepsilon}$ implies global exponential stability of $\A_{\varepsilon}$ for $\HS_{\varepsilon}$. This is accomplished by exploiting the relationship that exists between the two systems through the transformation matrix $\mathcal{T}$ on the sets $\tilde{\A}_{\varepsilon}$ and $\A_{\varepsilon}$.
}{
In the following two results, we first demonstrate the relationship between the sets $\tilde{\A}_{\varepsilon}$ for $\widetilde{\HS}_{\varepsilon}$ and $\A_{\varepsilon}$ for $\HS_{\varepsilon}$ so as to solve Problem \ref{prob:1}. Then, similar to Lemma \ref{lem:GES1}, we show that global exponential stability of $\tilde{\A}_{\varepsilon}$ for $\widetilde{\HS}_{\varepsilon}$ implies global exponential stability of $\A_{\varepsilon}$ for $\HS_{\varepsilon}$. See the appendix for proofs.


\begin{lemma} \label{lem:equiv_transform} 
Let $T_2 \geq T_1 > 0$, digraph $\DG$, and hybrid systems $\HS_{\varepsilon}$ and $\tilde{\HS}_{\varepsilon}$ be given as in (\ref{eqn:hy_eps}) and (\ref{eqn:H_vareps}), respectively. For each solution $\phi \in \mathcal{S}_{\HS_{\varepsilon}}$ there exists a solution $\tilde{\phi} \in \mathcal{S}_{\tilde{\HS}_{\varepsilon}}$ such that $\phi(t,j) = \Gamma \tilde{\phi}(t,j)$ for each $(t,j) \in \mbox{dom } \phi$ if and only if for each solutions $\tilde{\phi} \in \mathcal{S}_{\tilde{\HS}_{\varepsilon}}$ there exists a solution $\phi \in \mathcal{S}_{\HS_{\varepsilon}}$ such that $\tilde{\phi}(t,j) = \Gamma^{\minus 1} \phi(t,j)$ for each $(t,j) \in \mbox{dom } \tilde{\phi}$, where $\Gamma = {\rm diag}(\mathcal{T}, \mathcal{T}, \mathcal{T}, \mathcal{T}, 1)$.
\end{lemma}

\ifbool{two_col}{}{
\noindent See the appendix for proof.
}
}


\ifbool{two_col}{
\begin{lemma} \label{lem:set_equiv}
Given $0 < T_1 \leq T_2$ and a strongly connected digraph $\DG$, $\xi \in \A_{\varepsilon}$ if and only if $\chi_{\varepsilon} := \Gamma^{\minus 1} \xi \in \tilde{\A}_{\varepsilon}$, where $\Gamma^{\minus 1} = \mbox{\rm diag}(\mathcal{T}^{\minus 1}, \mathcal{T}^{\minus 1}, \mathcal{T}^{\minus 1}, \mathcal{T}^{\minus 1},  1)$ and $\mathcal{T}$ is given in (\ref{eqn:T_mat}). Moreover, for each $x_{\varepsilon} \in \mathcal{X}_{\varepsilon}$ and each $\chi_{\varepsilon} \in \mathcal{X}_{\varepsilon}$, $|\chi_{\varepsilon}|_{\tilde{\A}_{\varepsilon}} \leq |\Gamma^{\minus 1}| |x_{\varepsilon}|_{\A_{\varepsilon}}$ and $|x_{\varepsilon}|_{\A_{\varepsilon}} \leq |\Gamma| |\chi_{\varepsilon}|_{\tilde{\A}_{\varepsilon}}$.
\end{lemma}}{
\begin{lemma} \label{lem:set_equiv}
Given $0 < T_1 \leq T_2$ and a strongly connected digraph $\DG$, $\xi \in \A_{\varepsilon}$ if and only if $\chi_{\varepsilon} := \Gamma^{\minus 1} \xi \in \tilde{\A}_{\varepsilon}$, where $\Gamma^{\minus 1} = \mbox{\rm diag}(\mathcal{T}^{\minus 1}, \mathcal{T}^{\minus 1}, \mathcal{T}^{\minus 1}, \mathcal{T}^{\minus 1},  1)$ and $\mathcal{T}$ is given in (\ref{eqn:T_mat}). Moreover, for each $x_{\varepsilon} \in \mathcal{X}_{\varepsilon}$ and each $\chi_{\varepsilon} \in \mathcal{X}_{\varepsilon}$
\begin{equation} \label{eqn:set_bound_a}
|\chi_{\varepsilon}|_{\tilde{\A}_{\varepsilon}} \leq |\Gamma^{\minus 1}| |x_{\varepsilon}|_{\A_{\varepsilon}}
\end{equation}
and
\begin{equation} \label{eqn:set_bound_b}
|x_{\varepsilon}|_{\A_{\varepsilon}} \leq |\Gamma| |\chi_{\varepsilon}|_{\tilde{\A}_{\varepsilon}}
\end{equation}
\end{lemma}}

\ifbool{two_col}{}{
\noindent See the appendix for proof.
}


\begin{lemma} \label{lem:GES2}
Given $0 < T_1 \leq T_2$ and a strongly connected digraph $\DG$, the set $\tilde{\A}_{\varepsilon}$ is GES for the hybrid system $\widetilde{\HS}_{\varepsilon}$ if and only if $\A_{\varepsilon}$ is GES for the hybrid system $\HS_{\varepsilon}$.
\end{lemma}

\ifbool{two_col}{}{
\noindent See the appendix for proof.
}

\subsection{Parameter Estimator}

Exponential stability of the set $\tilde{\A}_{\varepsilon}$ for $\widetilde{\HS}_{\varepsilon}$ hinges upon the convergence of the estimate $\hat{a}$ to $a$. We present a result establishing convergence of $\hat{a}$ to $a$ by considering a model reduction of $\widetilde{\HS}_{\varepsilon}$. To this end, consider the state $\chi_{\varepsilon_r} := (\bar{w}_1, \bar{w}_2, \tau) \in \reals^2 \times \reals^{2(n-1)} \times [0,T_2] =: \mathcal{X}_{\varepsilon_r}$. Its dynamics are given by the system $\widetilde{\HS}_{\varepsilon_r} = ( \widetilde{C}_{\varepsilon_r}, \widetilde{f}_{\varepsilon_r},  \widetilde{D}_{\varepsilon_r}, \widetilde{G}_{\varepsilon_r})$ with data
\ifbool{two_col}{
$\widetilde{f}_{\varepsilon}(\chi_{\varepsilon_r})$ for each $\chi_{\varepsilon_r} \in \widetilde{C}_{\varepsilon_r}$ and $\widetilde{G}_{\varepsilon_r}(\chi_{\varepsilon_r})$ for each $\chi_{\varepsilon_r} \in \widetilde{D}_{\varepsilon_r}$  where
\ifbool{two_col}{
\begin{equation} \label{eqn:hy_param}
\begin{aligned}
\widetilde{f}_{\varepsilon_r}(\chi_{\varepsilon_r}) & := \big ( A_{f_3} \bar{w}_1, A_{f_4} \bar{w}_2,  \minus 1 \big ), \\
\widetilde{G}_{\varepsilon_r}(\chi_{\varepsilon_r}) & := \big ( \bar{w}_1, \bar{w}_2, [T_1,T_2] \big )
\end{aligned}
\end{equation}
}{
\begin{equation} \label{eqn:hy_param}
\begin{aligned}
\widetilde{f}_{\varepsilon_r}(\chi_{\varepsilon_r}) := \begin{bmatrix} A_{f_3} \bar{w}_1 \\ A_{f_4} \bar{w}_2 \\  \minus 1
\end{bmatrix}, \hspace{2mm} 
\widetilde{G}_{\varepsilon_r}(\chi_{\varepsilon_r}) := \begin{bmatrix}
\bar{w}_1 \\ \bar{w}_2 \\ [T_1,T_2]
\end{bmatrix}
\end{aligned}
\end{equation}
}
\noindent
where
$\widetilde{C}_{\varepsilon_r} := \mathcal{X}_{\varepsilon_r}$ and $\widetilde{D}_{\varepsilon_r} := \{ \chi_{\varepsilon_r} {\in} \mathcal{X}_{\varepsilon_r} : \tau = 0 \}$. 
}{
\begin{equation*}
\begin{aligned}
\tilde{f}_{\varepsilon_r}(\chi_{\varepsilon_r}) := \begin{bmatrix} A_{f_3} \bar{w}_1 \\ A_{f_4} \bar{w}_2 \\  \minus 1
\end{bmatrix} \hspace{2mm}  \forall \chi_{\varepsilon_r} \in \tilde{C}_{\varepsilon_r} := \mathcal{X}_{\varepsilon_r},  & &
\tilde{G}_{\varepsilon_r}(\chi_{\varepsilon_r}) := \begin{bmatrix}
\bar{w}_1 \\ \bar{w}_2 \\ [T_1,T_2]
\end{bmatrix} \hspace{2mm}  \forall \chi_{\varepsilon_r} \in \tilde{D}_{\varepsilon_r} := \{ \chi_{\varepsilon_r} {\in} \mathcal{X}_{\varepsilon_r} : \tau {=} 0 \} 
\end{aligned}
\end{equation*}
}
\noindent
For this system, the set to exponentially stabilize is given by 
\begin{equation} \label{set:est_param}
\widetilde{\A}_{\varepsilon_r} := \{0\} \times \{0\} \times [0,T_2]
\end{equation}

In the next result, we show global exponential stability of the set $\tilde{\A}_{\varepsilon_r}$ for $\widetilde{\HS}_{\varepsilon_r}$ through the satisfaction of matrix inequalities. See the appendix for proof.


\ifbool{conf}{
\begin{proposition} \label{prop:est_param} Let $0 < T_1 \leq T_2$ be given. If there exist a positive scalar $\mu$ and positive definite symmetric matrices $P_2$, $P_3$ such that 
\begin{equation} \label{eqn:lyap1}
P_2 A_{f_3} + A_{f_3}^{\top} P_2 \prec 0
\end{equation}
\begin{equation} \label{eqn:lyap2}
P_3 A_{f_4} + A_{f_4}^{\top} P_3 \prec 0
\end{equation}
hold, then the set $\widetilde{\A}_{\varepsilon_r}$ is globally exponentially stable for the hybrid system $\widetilde{\HS}_{\varepsilon_r}$. Furthermore, every solution $\phi$ to $\widetilde{\HS}_{\varepsilon_r}$ satisfies 
\begin{equation} \label{eqn:est_param_bnd}
|\phi(t,j)|_{\widetilde{\A}_{\varepsilon_r}} \leq \sqrt{\frac{\alpha_2}{\alpha_1}} \exp{ \Big (\minus \frac{\tilde{\beta}}{2 \alpha_2} t \Big ) } |\phi(0,0)|_{\widetilde{\A}_{\varepsilon_r}}
\end{equation}
\noindent
for each $(t,j) \in \mbox{dom } \phi$, with {\small $\alpha_1 = \mathrm{min} \{ \lambda_{min} (P_2), \lambda_{min} (P_3) \}$ $\alpha_2 = \mathrm{max} \{ \lambda_{max} (P_2), \lambda_{max} (P_3) \}$}
and $\tilde{\beta} > 0$.
\end{proposition}
}{
\begin{proposition} \label{prop:est_param}If there exists a positive scalar $\mu$ and positive definite symmetric matrices $P_2$, $P_3$ such that,  with $A_{f_3}$ and $A_{f_4}$ as in (\ref{eqn:matrices}), 
\ifbool{two_col}{the conditions in (\ref{eqn:lyap1})}{
\begin{equation} \label{eqn:lyap1b}
P_2 A_{f_3} + A_{f_3}^{\top} P_2 \prec 0
\end{equation}
\begin{equation} \label{eqn:lyap2b}
P_3 A_{f_4} + A_{f_4}^{\top} P_3 \prec 0
\end{equation}}
hold, then the set $\widetilde{\A}_{\varepsilon_r}$ is globally exponentially stable for the hybrid system $\widetilde{\HS}_{\varepsilon_r}$. Furthermore, every solution $\tilde{\phi}$ to $\widetilde{\HS}_{\varepsilon_r}$ satisfies 
\begin{equation} \label{eqn:est_param_bnd}
| \tilde{\phi}  (t,j)|_{\widetilde{\A}_{\varepsilon_r}} \leq \sqrt{\frac{\alpha_{\bar{w}_2}}{\alpha_{\bar{w}_1}}} \exp{ \Big ( \minus \frac{ \bar{\gamma}   \tilde{\beta}}{2 \alpha_{\bar{w}_2}}  (t + j)  \Big )} | \tilde{\phi}  (0,0)|_{\widetilde{\A}_{\varepsilon_r}}
\end{equation} 
for each $(t,j) \in \mbox{\rm dom }  \tilde{\phi}  $, with $\alpha_{\bar{w}_1} =$ $\mathrm{min} \{ \lambda_{min} (P_2),$ $\lambda_{min} (P_3) \}$, $\alpha_{\bar{w}_2} =$ $\mathrm{max} \{ \lambda_{max} (P_2),$ $\lambda_{max} (P_3) \}$, $\tilde{\beta} > 0$, and $\bar{\gamma} = \min \{1-\gamma, \gamma T_1\}$.
\end{proposition}

\ifbool{two_col}{}{
\noindent See the appendix for proof.
}
}

\subsection{Proof of Theorem \ref{thrm1}} \label{sec:thm_proof}

\ifbool{two_col}{
Consider the following Lyapunov function candidate  for $\widetilde{\HS}_{\varepsilon}$  
\begin{equation} \label{eqn:lyap}
V(\chi_{\varepsilon}) := V_1(\chi_{\varepsilon}) + V_2(\chi_{\varepsilon}) + V_{\varepsilon_r}(\chi_{\varepsilon}) \hspace{3mm} \forall \chi_{\varepsilon} \in \mathcal{X}_{\varepsilon}
\end{equation}
\noindent \sloppy
where $V_1(\chi_{\varepsilon}) {=} \exp{(2 h \tau)} \bar{\eta}_1^2$, $$V_2(\chi_{\varepsilon}) {=}\bar{z}^{\top}_2 \exp{(A^{\top}_{f_2} \tau )} P_1 \exp{(A_{f_2} \tau )} \bar{z}_2$$ and $V_{\varepsilon_r}(\chi_{\varepsilon}) =$ $\bar{w}_1^{\top} P_2 \bar{w}_1$ $+ \bar{w}_2^{\top} P_3 \bar{w}_2$. Note that there exist two positive scalars $\alpha_1$, $\alpha_2$ such that 
\begin{equation} \label{eqn:lyap_def}
\alpha_1 |\chi_{\varepsilon}|_{\tilde{\A}_{\varepsilon}}^2 \leq V(\chi_{\varepsilon}) \leq \alpha_2 |\chi_{\varepsilon}|_{\tilde{\A}_{\varepsilon}}^2 \hspace{5mm} \forall \chi_{\varepsilon} \in \tilde{C}_{\varepsilon} \cup \tilde{D}_{\varepsilon}
\end{equation} 
\noindent
With $P_1$ positive definite and noting the nonsingularity of $\exp{(A_{f_2} \tau)}$ for every $\tau$, we have
\ifbool{two_col}{
$\alpha_1 = \underset{\nu \in [0,T_2]}{\mathrm{min}} \Big \{ \exp{(2h \nu)},$ $\lambda_{min} \big (\exp{(A^{\top}_{f_2} \nu )} P_1 \exp{(A_{f_2} \nu )} \big ), \lambda_{min}(P_2), \lambda_{min}(P_3) \Big \}$
}{
\begin{equation*} \small
\begin{aligned}
\alpha_1 & = \underset{\nu \in [0,T_2]}{\mathrm{min}} \Big \{ \exp{(2h \nu)}, \lambda_{min} \big (\exp{(A^{\top}_{f_2} \nu )} P_1 \exp{(A_{f_2} \nu )} \big ), \\ 
& \hspace{53mm} \lambda_{min}(P_2), \lambda_{min}(P_3) \Big \}
\end{aligned}
\end{equation*}}
 and $\alpha_2$ as in (\ref{eqn:thm_consts}).  For each $\chi_{\varepsilon} \in \tilde{C}_{\varepsilon}$, one has
\begin{equation} \label{eqn:change_in_v}
\begin{aligned}
\langle \nabla V(\chi_{\varepsilon}), \tilde{f}_{\varepsilon}(\chi_{\varepsilon}) \rangle & = 2 \bar{z}_2^{\top} \big ( \exp{(A^{\top}_{f_2} \tau )} P_1 \exp{(A_{f_2} \tau)} \big ) B_{f_2}  \bar{w}_2 \\
& \hspace{10mm} + \bar{w}_1^{\top} (P_2 A_{f_3} + A_{f_3}^{\top} P_2) \bar{w}_1 \\
& \hspace{15mm} + \bar{w}_2^{\top} (P_3 A_{f_4} + A_{f_4}^{\top} P_3) \bar{w}_2 \\ 
\end{aligned}
\end{equation}
\noindent
Now, by noting the conditions in (\ref{eqn:lyap1}), with $\beta_1 > 0$ and $\beta_2 > 0$ such that $P_2 A_{f_3} + A_{f_3}^{\top} P_2 \leq -\beta_1 I$, and $P_3 A_{f_4} + A_{f_4}^{\top} P_3 \leq -\beta_2 I$ then one has
\begin{equation}
\begin{aligned}
\langle \nabla V(\chi_{\varepsilon}), \tilde{f}_{\varepsilon}(\chi_{\varepsilon}) \rangle & \leq \kappa_1 | \bar{z}_2 | | \bar{w}_2 | - \beta_1 | \bar{w}_1 |^2 - \beta_2 | \bar{w}_2 |^2 \\ 
\end{aligned}
\end{equation}
\noindent
where 
\ifbool{two_col}{ 
$\kappa_1$ is as given in (\ref{eqn:thm_consts}). Applying Young's inequality to $\kappa_1 |\bar{z}_2| |\bar{w}_2|$,\footnote{In particular, we are utilizing the relation $ab \leq \frac{a^2}{2 \epsilon} + \frac{\epsilon b^2}{2}$ where $a,b \in \reals$ and $\epsilon > 0$.} we obtain
}{
\begin{equation*}
\kappa_1 = 2 \underset{\nu \in [0,T_2]}{\mathrm{max}} \big | \exp{(A^{\top}_{f_2} \nu )} P_1 \exp{(A_{f_2} \nu )} \big | |B_{f_2}|
\end{equation*}
\noindent
Applying Young's inequality to $\kappa_1 |\bar{z}_2| |\bar{w}_2|$,\footnote{In particular, we are utilizing the relation $ab \leq \frac{a^2}{2 \epsilon} + \frac{\epsilon b^2}{2}$ where $a,b \in \reals$ and $\epsilon > 0$.} we obtain
}
\begin{equation}
\begin{aligned}
\langle \nabla V(\chi_{\varepsilon}), \tilde{f}_{\varepsilon}(\chi_{\varepsilon}) \rangle  
& \leq \frac{\kappa_1}{2 \epsilon} | \bar{z}_2 
|^2 \minus \beta_1 | \bar{w}_1|^2 + \Big ( \frac{\kappa_1 \epsilon}{2} - \beta_2 \Big ) | \bar{w}_2|^2
\end{aligned}
\end{equation}
\noindent
where $\epsilon > 0$. We then upper bound the inequality by picking the largest coefficient, i.e, $\bar{\kappa}_1 = \max \Big \{ \frac{\kappa_1}{2 \epsilon}, \Big ( \frac{\kappa_1 \epsilon}{2} - \beta_2 \Big )  \Big \}$, leading to
\begin{equation}
\langle \nabla V(\chi_{\varepsilon}), \tilde{f}_{\varepsilon}(\chi_{\varepsilon}) \rangle \leq \frac{\bar{\kappa}_1}{\alpha_2} V(\chi_{\varepsilon}) 
\end{equation}
\noindent
Now, for the analysis across jumps, note that for all $\chi_{\varepsilon} \in \tilde{D}_{\varepsilon}$, $\tau = 0$. At jumps, $\tau$ is  mapped to some point $\nu \in [T_1 , T_2]$. Then, at jumps, for each $g \in \tilde{G}_{\varepsilon}$ one has
\begin{equation*}
\begin{aligned}
V(g) {\minus} V(\chi_{\varepsilon}) & = \\ 
& \hspace{-5mm} \minus \bar{\eta}_1^2 {+} \bar{z}_2^{\top} \big ( A_{g_2}^{\top} \exp{(A_{f_2}^{\top} \nu)} P_1 \exp{(A_{f_2} \nu)} A_{g_2} {\minus} P_1 \big ) \bar{z}_2 \\
& \leq \minus \bar{\kappa}_2 \big ( | \bar{\eta}_1 |^2 + | \bar{z}_2 |^2 \big )
\end{aligned}
\end{equation*}
\normalsize
\ifbool{two_col}{
where $\bar{\kappa}_2$ and $\kappa_2$ are as given in (\ref{eqn:thm_consts}), from where we have 
}{
\noindent
where $\bar{\kappa}_2 = \max \{1 , \kappa_2 \}$ and, by continuity of condition (\ref{cond:phil_1}), $\kappa_2 > 0$ such that
\small
\begin{equation*}
\kappa_2 {\in} \Big ( 0,  \minus \underset{v {\in} [T_1,T_2]}{\min}  {\lambda_{min}} (A^{\top}_{g_2} \exp{(A^{\top}_{f_2} v )} P_1 \exp{(A_{f_2} v )} A_{g_2} {\minus} P_1 ) \Big )
\end{equation*} 
\normalsize
\noindent
 from where we have 
}
\begin{equation} \label{eqn:jump_bnd}
\begin{aligned}
V(g) - V(\chi_{\varepsilon}) \leq \minus \bar{\kappa}_2 \big ( | \bar{\eta}_1 |^2 + | \bar{z}_2 |^2 \big )
\end{aligned}
\end{equation}
\noindent
 Utilizing the upper bound $\alpha_2$ from the definition of $V$ in (\ref{eqn:lyap_def}),  for all $\chi_{\varepsilon} \in \tilde{D}_{\varepsilon}$, one has $V(\chi_{\varepsilon}) \leq \alpha_2 \big (| \bar{\eta}_1 |^2 + |\bar{z}_2|^2 + |\bar{w}|^2 \big )$. Dividing by $\alpha_2$ and rearranging terms, one has
\begin{equation} \label{thrm_zdef}
\begin{aligned}
\minus (| \bar{\eta}_1 |^2 + |\bar{z}_2|^2) \leq - \frac{1}{\alpha_2} V(\chi_{\varepsilon}) + |\bar{w}|^2
\end{aligned}
\end{equation}
\noindent
Then, by inserting (\ref{thrm_zdef}) into (\ref{eqn:jump_bnd}), we obtain
\begin{equation} \label{eqn:v_jumps}
V(g) - V(\chi_{\varepsilon}) \leq \bar{\kappa}_2 \Big (- \frac{1}{\alpha_2} V(\chi_{\varepsilon}) + |\bar{w}|^2 \Big )
\end{equation}
\noindent
Now, by noting that $\langle \nabla V(\chi_{\varepsilon}), \tilde{f}(\chi_{\varepsilon}) \rangle  \leq \frac{\bar{\kappa}_1}{\alpha_2} V(\chi_{\varepsilon})$ and by (\ref{eqn:v_jumps}), pick a solution  $\tilde{\phi}$ to $\widetilde{\HS}_{\varepsilon}$  with initial condition $\tilde{\phi}(0,0) \in \widetilde{C}_{\varepsilon} \cup \widetilde{D}_{\varepsilon}$.  Let the jumps of $\tilde{\phi}$ occur at times $(t_j,j) \in \{j' : \exists t' : (t',j') \in \mbox{dom } \phi\}$.  For each $(t,j) \in [0,t_1] \times \{0\}$ one has
\begin{equation*}
V(  \tilde{\phi} (  t,0 )) \leq \exp \Big ( \frac{\bar{\kappa}_1}{\alpha_2} \hspace{0.5mm} t_1 \Big ) V( \tilde{\phi} (  0,0)) 
\end{equation*}
\noindent
At $(t_1,1)$, one has
\begin{equation*}
\begin{aligned}
V(  \tilde{\phi} (  t_1,1)) & {\leq} \Big ( 1 \minus \frac{\bar{\kappa}_2}{\alpha_2} \Big ) \exp \Big ( \frac{\bar{\kappa}_1}{\alpha_2} \hspace{0.5mm} t_1 \Big ) V(  \tilde{\phi} (  0,0)) {+} \bar{\kappa}_2 |\bar{w}(t_1,0)|^2 \\
\end{aligned}
\end{equation*}
\noindent
Then, for each $(t,j) \in [t_1,t_2] \times \{1\}$
\begin{equation*}
\begin{aligned}
V(  \tilde{\phi} (  t,1)) & = \exp \Big ( \frac{\bar{\kappa}_1}{\alpha_2} \hspace{0.5mm} t_2 \Big ) \Big ( 1 - \frac{\bar{\kappa}_2}{\alpha_2} \Big ) V(  \tilde{\phi} (  0,0)) \\ 
& \hspace{1mm} + \exp \big ( \frac{\bar{\kappa}_1}{\alpha_2} (t_2 - t_1) \big ) \bar{\kappa}_2 |\bar{w}(t_1,0)|^2
\end{aligned}
\end{equation*}
\normalsize
\noindent
At $(t_2,2)$, one has
\begin{equation*}
\begin{aligned}
V(  \tilde{\phi} (  t_2,2)) & \leq \exp \Big ( \frac{\bar{\kappa}_1}{\alpha_2} \hspace{0.5mm} t_2 \Big ) \Big ( 1 - \frac{\bar{\kappa}_2}{\alpha_2} \Big )^2 V(  \tilde{\phi} (  0,0)) \\
& \hspace{-2mm} + \bar{\kappa}_2 \Big [ \exp \Big ( \frac{\bar{\kappa}_1}{\alpha_2} (t_2 - t_1) \Big ) |\bar{w}(t_1,0)|^2 + |\bar{w}(t_2,1)|^2 \Big ]
\end{aligned}
\end{equation*}
\normalsize
\noindent
A general form of the bound is given by
\begin{equation} \label{eqn:gen_bound_main}
\begin{aligned}
V(  \tilde{\phi} (  t,j)) & \leq \exp \Big ( \frac{\bar{\kappa}_1}{\alpha_2} \hspace{0.5mm} t_j \Big ) \Big ( 1 {-} \frac{\bar{\kappa}_2}{\tilde{\alpha}_1} \Big )^j V(  \tilde{\phi} (  0,0)) \\
& \hspace{2mm} + \bar{\kappa}_2 \Big( \sum_{k=1}^{j} \exp \Big ( \frac{\bar{\kappa}_1}{\alpha_2} (t_{k + 1} \minus t_k) \Big ) |\bar{w}(t_{k},k \minus 1)|^2 \Big )
\end{aligned}
\end{equation}
\normalsize
\noindent
Noting that $t_{j+1} - t_j \leq T_2$ and $\frac{\bar{\kappa}_1}{\alpha_2} > 0$,  the latter term can be further bounded as
\begin{equation*}
\begin{aligned}
& \bar{\kappa}_2 \Big( \sum_{k=1}^{j} {\exp} \big ( \frac{\bar{\kappa}_1}{\alpha_2} (t_{k+1} \minus t_k) \big ) |\bar{w}(t_{k},k \minus 1)|^2 \Big ) \\
& \hspace{30mm} \leq  \bar{\kappa}_2 {\exp} \big ( \frac{\bar{\kappa}_1}{\alpha_2} \hspace{0.5mm} T_2 \big ) \mathrm{sup}_{(t,j) \in \mbox{\footnotesize dom} \tilde{\phi}} |\bar{w}(t,j)|^2
\end{aligned}
\end{equation*}
\noindent
Moreover, since $t_j \leq T_2 (j+1)$ and $\frac{\bar{\kappa}_1}{\alpha_2} > 0$, we can also put a stricter bound on the first term in (\ref{eqn:gen_bound_main}) as follows:
\begin{equation*}
\begin{aligned}
\hspace{0mm} {\exp} \Big ( \frac{\bar{\kappa}_1}{\alpha_2} \hspace{0.5mm} t_{j} \Big ) \Big ( 1 {\minus} \frac{\bar{\kappa}_2}{\alpha_2} \Big )^j V(  \tilde{\phi} (  0,0)) & \\
& \hspace{-35mm} \leq {\exp} \Big ( \frac{\bar{\kappa}_1}{\alpha_2} \hspace{0.5mm} T_2 \Big ) \Big ( \exp \Big ( \frac{\bar{\kappa}_1}{\alpha_2} \hspace{0.5mm} T_2 \Big ) \Big ( 1 {-} \frac{\bar{\kappa}_2}{\alpha_2} \Big ) \Big )^j V(  \tilde{\phi} (  0,0))
\end{aligned}
\end{equation*}
\normalsize
\noindent 
Thus 
\begin{equation*} \label{thrm1_bound_a}
\begin{aligned}
V(  \tilde{\phi} (  t,j)) & \leq \\ 
& \exp \Big ( \frac{\bar{\kappa}_1}{\alpha_2} \hspace{0.5mm} T_2 \Big ) \Big ( \exp \Big ( \frac{\bar{\kappa}_1}{\alpha_2} \hspace{0.5mm} T_2 \Big ) \Big ( 1 {-} \frac{\bar{\kappa}_2}{\alpha_2} \Big ) \Big )^j V(  \tilde{\phi} (  0,0)) \\
& \hspace{13mm} + \bar{\kappa}_2 \exp \Big ( \frac{\bar{\kappa}_1}{\alpha_2} \hspace{0.5mm} T_2 \Big ) \mathrm{sup}_{(t,j) \in \mbox{\footnotesize dom} \tilde{\phi}}|\bar{w}(t,j)|^2
\end{aligned}
\end{equation*}
\normalsize
\noindent
Then, from the result of Proposition \ref{prop:est_param}, we have (\ref{eqn:est_param_bnd})
with $\alpha_{\bar{w}_1} =$ $\mathrm{min} \{ \lambda_{min} (P_2), \lambda_{min} (P_3) \}$ and $\alpha_{\bar{w}_2} =$ $\mathrm{max} \{ \lambda_{max} (P_2), \lambda_{max} (P_3) \}$.  Now, to improve readability, we have omitted including the use of the notation $V(\tilde{\phi}(t,j))$ when evaluating $V$ along the trajectory for the solution $\tilde{\phi}$ opting instead for the use of the state components of $\chi_{\varepsilon}$ directly. In particular, we remind the reader that the notation $\bar{w}(t,j)$ corresponds to the $\bar{w}$ component of a solution, i.e.,  $\phi_{\bar{w}}(t,j)$. Now, combining the inequality with (\ref{eqn:lyap_def}) and noting $V(\tilde{\phi}(0,0)) \leq \alpha_2 |\tilde{\phi}(0,0) |_{\tilde{\A}_{\varepsilon}}^2$ one has for each $(t,j) \in \mbox{dom } \phi$
\begin{equation} \label{thrm1_bound}
\begin{aligned} 
|\phi(t,j)|_{\tilde{\A}_{\varepsilon}} & \leq \\
& \hspace{-15mm} \sqrt{\frac{\alpha_2}{\alpha_1}} |\phi(0,0) |_{\tilde{\A}_{\varepsilon}} \exp \Big ( \frac{\bar{\kappa}_1}{2 \alpha_2} \hspace{0.5mm} T_2 \Big ) \Big ( \exp \Big ( \frac{\bar{\kappa}_1}{2 \alpha_2} \hspace{0.5mm} T_2 \Big ) \Big ( 1 {-} \frac{\bar{\kappa}_2}{2 \alpha_2} \Big ) \Big )^j \\
& \hspace{-15mm} {+} \sqrt{\bar{\kappa}_2} \exp \Big ( \frac{\bar{\kappa}_1}{2 \alpha_2} \hspace{0.5mm} T_2 \Big ) \sqrt{ \frac{\alpha_{\bar{w}_2}}{\alpha_{\bar{w}_1}} \exp{ \Big ( \frac{ \minus \bar{\gamma} \tilde{\beta}}{2 \alpha_{\bar{w}_2}} (t {+} j) \Big )^2 } |\phi_{\bar{w}}(0,0)|^2_{\tilde{\A}_{\varepsilon_r}} }
\end{aligned}
\end{equation}
\normalsize
By the given conditions, the set $\tilde{\A}_{\varepsilon}$ is globally exponentially stable and attractive for $\widetilde{\HS}_{\varepsilon}$.  Now, by utilizing Lemmas \ref{lem:equiv_transform} - \ref{lem:GES2}, we can establish global exponential stability to the set $\A_{\varepsilon}$ for $\HS_{\varepsilon}$. In particular, Lemma \ref{lem:equiv_transform} establishes the relation between  $\widetilde{\HS}_{\varepsilon}$ and $\HS_{\varepsilon}$. In turn, we can then make use of Lemmas  \ref{lem:equiv} - \ref{lem:GES1}, where Lemma \ref{lem:equiv} establishes the reduction from $\HS$ to $\HS_{\varepsilon}$, to show that the set $\A$ is globally exponentially stable and attractive for $\HS$ in (\ref{eqn:cl_hysys}). 
}{
Consider the following Lyapunov function candidate  for $\widetilde{\HS}_{\varepsilon}$  
\begin{equation} \label{eqn:lyap}
V(\chi_{\varepsilon}) := V_1(\chi_{\varepsilon}) + V_2(\chi_{\varepsilon}) + V_{\varepsilon_r}(\chi_{\varepsilon}) \hspace{3mm} \forall \chi_{\varepsilon} \in \mathcal{X}_{\varepsilon}
\end{equation}
\noindent
where 
\begin{align*}
& V_1(\chi_{\varepsilon}) = \exp{(2 h \tau)} \bar{\eta}_1^2 \\
& V_2(\chi_{\varepsilon}) = \bar{z}^{\top}_2 \exp{(A^{\top}_{f_2} \tau )} P_1 \exp{(A_{f_2} \tau )} \bar{z}_2 \\
& V_{\varepsilon_r}(\chi_{\varepsilon}) = \bar{w}_1^{\top} P_2 \bar{w}_1 + \bar{w}_2^{\top} P_3 \bar{w}_2
\end{align*} Note that there exist two positive scalars $\alpha_1$, $\alpha_2$ such that 
\begin{equation} \label{eqn:lyap_def}
\alpha_1 |\chi_{\varepsilon}|_{\tilde{\A}_{\varepsilon}}^2 \leq V(\chi_{\varepsilon}) \leq \alpha_2 |\chi_{\varepsilon}|_{\tilde{\A}_{\varepsilon}}^2 \hspace{1cm} \forall \chi_{\varepsilon} \in \tilde{C}_{\varepsilon} \cup \tilde{D}_{\varepsilon}
\end{equation} 
\noindent
With $P_1$ positive definite and noting the nonsingularity of $\exp{(A_{f_2} \tau)}$ for every $\tau$, we have

\begin{equation*} \small
\begin{aligned}
\alpha_1 & = \underset{\nu \in [0,T_2]}{\mathrm{min}} \Big \{ \exp{(2h \nu)}, \lambda_{min} \big (\exp{(A^{\top}_{f_2} \nu )} P_1 \exp{(A_{f_2} \nu )} \big ), \\ 
& \hspace{53mm} \lambda_{min}(P_2), \lambda_{min}(P_3) \Big \}
\end{aligned}
\end{equation*}
 and $\alpha_2$ as in (\ref{eqn:thm_consts}).  For each $\chi_{\varepsilon} \in \tilde{C}_{\varepsilon}$, one has
\begin{equation}
\begin{aligned}
\langle \nabla V(\chi_{\varepsilon}), \tilde{f}_{\varepsilon}(\chi_{\varepsilon}) \rangle & = 2 \bar{z}_2^{\top} \big ( \exp{(A^{\top}_{f_2} \tau )} P_1 \exp{(A_{f_2} \tau)} \big ) B_{f_2}  \bar{w}_2 \\
& \hspace{10mm} + \bar{w}_1^{\top} (P_2 A_{f_3} + A_{f_3}^{\top} P_2) \bar{w}_1 \\
& \hspace{15mm} + \bar{w}_2^{\top} (P_3 A_{f_4} + A_{f_4}^{\top} P_3) \bar{w}_2 \\ 
\end{aligned}
\end{equation}
\noindent
Now, by noting (\ref{eqn:lyap1}) and (\ref{eqn:lyap2}), with $\beta_1 > 0$ and $\beta_2 > 0$ such that $P_2 A_{f_3} + A_{f_3}^{\top} P_2 \leq -\beta_1 I$, and $P_3 A_{f_4} + A_{f_4}^{\top} P_3 \leq -\beta_2 I$ then one has
\begin{equation}
\begin{aligned}
\langle \nabla V(\chi_{\varepsilon}), \tilde{f}_{\varepsilon}(\chi_{\varepsilon}) \rangle & \leq \kappa_1 | \bar{z}_2 | | \bar{w}_2 | - \beta_1 | \bar{w}_1 |^2 - \beta_2 | \bar{w}_2 |^2 \\ 
\end{aligned}
\end{equation}
\noindent
where
\begin{equation*}
\kappa_1 = 2 \underset{\nu \in [0,T_2]}{\mathrm{max}} \big | \exp{(A^{\top}_{f_2} \nu )} P_1 \exp{(A_{f_2} \nu )} \big | |B_{f_2}|
\end{equation*}
\noindent
Applying Young's inequality to $\kappa_1 |\bar{z}_2| |\bar{w}_2|$, \footnote{In particular, we are utilizing the relation $ab \leq \frac{a^2}{2 \epsilon} + \frac{\epsilon b^2}{2}$ where $a,b \in \reals$ and $\epsilon > 0$.} we obtain
\begin{equation}
\begin{aligned}
\langle \nabla V(\chi_{\varepsilon}), \tilde{f}_{\varepsilon}(\chi_{\varepsilon}) \rangle 
& \leq \frac{\kappa_1}{2 \epsilon} | \bar{z}_2 |^2 + \frac{\kappa_1 \epsilon}{2} | \bar{w}_2 |^2 \minus \beta_1 | \bar{w}_1|^2 \minus \beta_2 | \bar{w}_2|^2 \\ 
& \leq \frac{\kappa_1}{2 \epsilon} | \bar{z}_2 
|^2 \minus \beta_1 | \bar{w}_1|^2 + \Big ( \frac{\kappa_1 \epsilon}{2} - \beta_2 \Big ) | \bar{w}_2|^2
\end{aligned}
\end{equation}
\noindent
where $\epsilon > 0$, we then upper bound the inequality by picking the largest coefficient, i.e, $\bar{\kappa}_1 = \max \Big \{ \frac{\kappa_1}{2 \epsilon}, \Big ( \frac{\kappa_1 \epsilon}{2} - \beta_2 \Big )  \Big \}$, leading to
\begin{equation}
\begin{aligned}
\langle \nabla V(\chi_{\varepsilon}), \tilde{f}_{\varepsilon}(\chi_{\varepsilon}) \rangle & \leq \bar{\kappa}_1 \big ( | \bar{z}_2 |^2 + | \bar{w}_1|^2 + | \bar{w}_2 |^2 \big ) \\
& \leq \bar{\kappa}_1 \big ( |\chi_{\varepsilon}|^2_{\tilde{\A}_{\varepsilon}} \big ) \\
& \leq \bar{\kappa}_1 \Big ( \frac{1}{\alpha_2} V(\chi_{\varepsilon}) \Big ) \\
& \leq \frac{\bar{\kappa}_1}{\alpha_2} V(\chi_{\varepsilon}) 
\end{aligned}
\end{equation}
\noindent
Now, for the analysis across jumps, note that for all $\chi_{\varepsilon} \in \tilde{D}_{\varepsilon}$, $\tau = 0$. At jumps, $\tau$ is  mapped to some point $\nu \in [T_1 , T_2]$. Then, at jumps, for each $g \in \tilde{G}_{\varepsilon}$ one has
\begin{equation}
\begin{aligned}
V(g) {\minus} V(\chi_{\varepsilon}) & = \minus \bar{\eta}_1^2 -  \bar{z}_2^{\top} P_1 \bar{z}_2 \\
& \hspace{5mm} {+} (A_{g_2} \bar{z}_2)^{\top} \exp{(A^{\top}_{f_2} \nu )} P_1 \exp{(A_{f_2} \nu)} (A_{g_2} \bar{z}_2) \\
& = \minus \bar{\eta}_1^2 \\
& \hspace{5mm} {+} \bar{z}_2^{\top} \big ( A_{g_2}^{\top} \exp{(A_{f_2}^{\top} \nu)} P_1 \exp{(A_{f_2} \nu)} A_{g_2} {\minus} P_1 \big ) \bar{z}_2 \\
& \leq \minus | \bar{\eta}_1 |^2 \minus \kappa_2 | \bar{z}_2 |^2 \\
& \leq \minus \bar{\kappa}_2 \big ( | \bar{\eta}_1 |^2 + | \bar{z}_2 |^2 \big )
\end{aligned}
\end{equation}
\normalsize
\noindent
where $\bar{\kappa}_2 = \max \{1 , \kappa_2 \}$ and, by continuity of condition (\ref{cond:phil_1}), $\kappa_2 > 0$ such that
\small
\begin{equation*}
\kappa_2 {\in} \Big ( 0,  \minus \underset{v {\in} [T_1,T_2]}{\min}  {\lambda_{min}} (A^{\top}_{g_2} \exp{(A^{\top}_{f_2} v )} P_1 \exp{(A_{f_2} v )} A_{g_2} {\minus} P_1 ) \Big )
\end{equation*} 
\normalsize
\noindent
 for where we have 
\begin{equation} \label{eqn:jump_bnd}
\begin{aligned}
V(g) - V(\chi_{\varepsilon}) \leq \minus \bar{\kappa}_2 \big ( | \bar{\eta}_1 |^2 + | \bar{z}_2 |^2 \big )
\end{aligned}
\end{equation}
\noindent
 Utilizing the upper bound $\alpha_2$ from the definition of $V$ in (\ref{eqn:lyap_def}),  for all $\chi_{\varepsilon} \in \tilde{D}_{\varepsilon}$, one has
\begin{equation}
\begin{aligned}
V(\chi_{\varepsilon}) & \leq \alpha_2 \big (| \bar{\eta}_1 |^2 + |\bar{z}_2|^2 + |\bar{w}|^2 \big ) \\
\end{aligned}
\end{equation}
\noindent
Dividing by $\alpha_2$ and rearranging terms, one has
\begin{equation} \label{thrm_zdef}
\begin{aligned}
\minus (| \bar{\eta}_1 |^2 + |\bar{z}_2|^2) \leq - \frac{1}{\alpha_2} V(\chi_{\varepsilon}) + |\bar{w}|^2 \\
\end{aligned}
\end{equation}
\noindent
Then, by inserting (\ref{thrm_zdef}) into (\ref{eqn:jump_bnd}),
\begin{equation} \label{eqn:v_jumps}
\begin{aligned}
V(g) - V(\chi_{\varepsilon}) & \leq \minus \bar{\kappa}_2 \big ( | \bar{\eta}_1 |^2 + | \bar{z}_2 |^2 \big ) \\
V(g) - V(\chi_{\varepsilon}) & \leq \bar{\kappa}_2 \Big (- \frac{1}{\alpha_2} V(\chi_{\varepsilon}) + |\bar{w}|^2 \Big ) \\
V(g) & \leq - \frac{\bar{\kappa}_2}{\alpha_2} V(\chi_{\varepsilon}) + \bar{\kappa}_2 |\bar{w}|^2 + V(\chi_{\varepsilon}) \\
V(g) & \leq \Big (1 - \frac{\bar{\kappa}_2}{\alpha_2} \Big ) V(\chi_{\varepsilon}) + \bar{\kappa}_2  |\bar{w}|^2 \\
\end{aligned}
\end{equation}
\noindent
Now, by noting that $\langle \nabla V(\chi_{\varepsilon}), \tilde{f}(\chi_{\varepsilon}) \rangle  \leq \frac{\bar{\kappa}_1}{\alpha_2} V(\chi_{\varepsilon})$ and by (\ref{eqn:v_jumps}), pick a solution  $\tilde{\phi}$ to $\widetilde{\HS}_{\varepsilon}$  with initial condition $\tilde{\phi}(0,0) \in \widetilde{C}_{\varepsilon} \cup \widetilde{D}_{\varepsilon}$.  Let the jumps of $\tilde{\phi}$ occur at times $(t_j,j) \in \{j' : \exists t' : (t',j') \in \mbox{dom } \phi\}$.  For each $(t,j) \in [0,t_1] \times \{0\}$ one has
\begin{equation}
V(t,0) \leq \exp \Big ( \frac{\bar{\kappa}_1}{\alpha_2} \hspace{0.5mm} t_1 \Big ) V(0,0) 
\end{equation}
\noindent
At $(t_1,1)$ 
\begin{equation*}
\begin{aligned}
V(t_1,1) & \leq \Big ( 1 - \frac{\bar{\kappa}_2}{\alpha_2} \Big ) V(t_1,0) + \bar{\kappa}_2 |\bar{w}(t_1,0)|^2 \\ 
& \leq \Big ( 1 \minus \frac{\bar{\kappa}_2}{\alpha_2} \Big ) \exp \Big ( \frac{\bar{\kappa}_1}{\alpha_2} \hspace{0.5mm} t_1 \Big ) V(0,0) + \bar{\kappa}_2 |\bar{w}(t_1,0)|^2 \\
\end{aligned}
\end{equation*}
\noindent
Then, for each $(t,j) \in [t_1,t_2] \times \{1\}$
\begin{equation*}
\begin{aligned}
V(t,1) & \leq \exp \Big ( \frac{\bar{\kappa}_1}{\alpha_2} (t_2 - t_1) \Big ) V(t_1,1) \\
& \leq \exp \Big ( \frac{\bar{\kappa}_1}{\alpha_2} (t_2 - t_1) \Big ) \Big [ \Big ( 1 - \frac{\bar{\kappa}_2}{\alpha_2} \Big ) \exp \big ( \bar{\kappa}_1 \hspace{0.5mm} t_1 \big ) V(0,0) \\
& \hspace{5mm} + \bar{\kappa}_2 |\bar{w}(t_1,0)|^2 \Big ] \\
& \leq \exp \Big ( \frac{\bar{\kappa}_1}{\alpha_2} (t_2 - t_1) \Big ) \Big ( 1 \minus  \frac{\bar{\kappa}_2}{\tilde{\alpha}_1} \Big ) \exp \Big ( \frac{\bar{\kappa}_1}{\alpha_2} \hspace{0.5mm} t_1 \Big ) V(0,0) \\
& \hspace{5mm} + \exp \Big ( \frac{\bar{\kappa}}{\alpha_2} (t_2 - t_1) \Big ) \bar{\kappa}_2 |\bar{w}(t_1,0)|^2 \\
& = \exp \Big ( \frac{\bar{\kappa}_1}{\alpha_2} \hspace{0.5mm} t_2 \Big ) \Big ( 1 - \frac{\bar{\kappa}_2}{\alpha_2} \Big ) V(0,0) \\ 
& \hspace{5mm} + \exp \big ( \frac{\bar{\kappa}_1}{\alpha_2} (t_2 - t_1) \big ) \bar{\kappa}_2 |\bar{w}(t_1,0)|^2 \\
\end{aligned}
\end{equation*}
\normalsize
\noindent
At $(t_2,2)$
\begin{equation*}
\begin{aligned}
V(t_2,2) & \leq \Big ( 1 - \frac{\bar{\kappa}_2}{\alpha_2} \Big ) V(t_2,1) + \bar{\kappa}_2 |\bar{w}(t_2,1)|^2 \\ 
& \leq \Big ( 1 - \frac{\bar{\kappa}_2}{\alpha_2} \Big ) \exp \Big ( \frac{\bar{\kappa}_1}{\alpha_2} \hspace{0.5mm} t_2 \Big ) \Big ( 1 - \frac{\bar{\kappa}_2}{\alpha_2} \Big ) V(0,0) \\
& \hspace{5mm} + \exp \Big ( \frac{\bar{\kappa}_1}{\alpha_2} (t_2 - t_1) \Big ) \bar{\kappa}_2  |\bar{w}(t_1,0)|^2 {+} \bar{\kappa}_2 |\bar{w}(t_2,1)|^2 \\
& \leq \exp \Big ( \frac{\bar{\kappa}_1}{\alpha_2} \hspace{0.5mm} t_2 \Big ) \Big ( 1 - \frac{\bar{\kappa}_2}{\alpha_2} \Big )^2 V(0,0) \\
& \hspace{5mm} + \bar{\kappa}_2 \Big [ \exp \Big ( \frac{\bar{\kappa}_1}{\alpha_2} (t_2 - t_1) \Big ) |\bar{w}(t_1,0)|^2 + |\bar{w}(t_2,1)|^2 \Big ] \\
\end{aligned}
\end{equation*}
\normalsize
\noindent
A general form of the bound is given by
\begin{equation} \label{eqn:gen_bound_main}
\begin{aligned}
V(t,j) & \leq \exp \Big ( \frac{\bar{\kappa}_1}{\alpha_2} \hspace{0.5mm} t_j \Big ) \Big ( 1 {-} \frac{\bar{\kappa}_2}{\tilde{\alpha}_1} \Big )^j V(0,0) \\
& \hspace{5mm} + \bar{\kappa}_2 \Big( \sum_{k=1}^{j} \exp \Big ( \frac{\bar{\kappa}_1}{\alpha_2} (t_{k + 1} \minus t_k) \Big ) |\bar{w}(t_{k},k \minus 1)|^2 \Big )
\end{aligned}
\end{equation}
\normalsize
\noindent
Noting that $t_{j+1} - t_j \leq T_2$ and $\frac{\bar{\kappa}_1}{\alpha_2} > 0$,  the latter term can be further bounded as
\begin{equation*}
\begin{aligned}
& \bar{\kappa}_2 \Big( \sum_{k=1}^{j} {\exp} \big ( \frac{\bar{\kappa}_1}{\alpha_2} (t_{k+1} \minus t_k) \big ) |\bar{w}(t_{k},k \minus 1)|^2 \Big ) \\
& \hspace{30mm} \leq  \bar{\kappa}_2 {\exp} \big ( \frac{\bar{\kappa}_1}{\alpha_2} \hspace{0.5mm} T_2 \big ) \mathrm{sup}_{(t,j) \in \mbox{\footnotesize dom} \tilde{\phi}} |\bar{w}(t,j)|^2
\end{aligned}
\end{equation*}
\noindent
Moreover, since $t_j \leq T_2 (j+1)$ and $\frac{\bar{\kappa}_1}{\alpha_2} > 0$, we can also put a stricter bound on the first term in (\ref{eqn:gen_bound_main}) as follows:
\begin{equation*}
\begin{aligned}
\hspace{-1mm} {\exp} \Big ( \frac{\bar{\kappa}_1}{\alpha_2} \hspace{0.5mm} t_{j} \Big ) \Big ( 1 {\minus} \frac{\bar{\kappa}_2}{\alpha_2} \Big )^j V(0,0) & \\
& \hspace{-25mm} \leq {\exp} \Big ( \frac{\bar{\kappa}_1}{\alpha_2} \hspace{0.5mm} T_2 (j{+}1) \Big ) \Big ( 1 {\minus} \frac{\bar{\kappa}_2}{\alpha_2} \Big )^j V(0,0) \\
& \hspace{-25mm} \leq {\exp} \Big ( \frac{\bar{\kappa}_1}{\alpha_2} \hspace{0.5mm} T_2 \Big ) \Big ( \exp \Big ( \frac{\bar{\kappa}_1}{\alpha_2} \hspace{0.5mm} T_2 \Big ) \Big ( 1 {-} \frac{\bar{\kappa}_2}{\alpha_2} \Big ) \Big )^j V(0,0)
\end{aligned}
\end{equation*}
\normalsize
\noindent 
Thus 
\begin{equation} \label{thrm1_bound_a}
\begin{aligned}
V(t,j) & \leq \exp \Big ( \frac{\bar{\kappa}_1}{\alpha_2} \hspace{0.5mm} T_2 \Big ) \Big ( \exp \Big ( \frac{\bar{\kappa}_1}{\alpha_2} \hspace{0.5mm} T_2 \Big ) \Big ( 1 - \frac{\bar{\kappa}_2}{\alpha_2} \Big ) \Big )^j V(0,0) \\
& \hspace{5mm} + \bar{\kappa}_2 \exp \Big ( \frac{\bar{\kappa}_1}{\alpha_2} \hspace{0.5mm} T_2 \Big ) \mathrm{sup}_{(t,j) \in \mbox{\footnotesize dom} \tilde{\phi}}|\bar{w}(t,j)|^2
\end{aligned}
\end{equation}
\normalsize
\noindent
Then, from the result of Proposition \ref{prop:est_param}, we have 
\begin{equation*}
|\tilde{\phi}_{\bar{w}}(t,j)| \leq \sqrt{\frac{\alpha_{\bar{w}_2}}{\alpha_{\bar{w}_1}}} \exp{ \Big ( \minus \frac{\tilde{\beta}}{2 \alpha_{\bar{w}_2}} t \Big ) } |\tilde{\phi}_{\bar{w}}(0,0)|_{\tilde{\A}_{\varepsilon_r}}
\end{equation*} 
with $\alpha_{\bar{w}_1} = \mathrm{min} \{ \lambda_{min} (P_2), \lambda_{min} (P_3) \}$ and $\alpha_{\bar{w}_2} = \mathrm{max} \{ \lambda_{max} (P_2), \lambda_{max} (P_3) \}$.  Now, to improve readability, we have omitted including the use of the notation $V(\tilde{\phi}(t,j))$ when evaluating $V$ along the trajectory for the solution $\tilde{\phi}$ opting instead for the use of the state components of $\chi_{\varepsilon}$ directly. In particular, we remind the reader that the notation $\bar{w}(t,j)$ corresponds to the $\bar{w}$ component of a solution, i.e.,  $\phi_{\bar{w}}(t,j)$. Thus, we have   
\ifbool{two_col}{
\begin{equation} \label{thrm1_bound}
\begin{aligned} 
V(t,j) & \leq \exp \Big ( \frac{\bar{\kappa}_1}{\alpha_2} \hspace{0.5mm} T_2 \Big ) \Big ( \exp \Big ( \frac{\bar{\kappa}_1}{\alpha_2} \hspace{0.5mm} T_2 \Big ) \Big ( 1 {-} \frac{\bar{\kappa}_2}{\alpha_2} \Big ) \Big )^j V(0,0) \\
& \hspace{-5mm} + \bar{\kappa}_2 \exp \Big ( \frac{\bar{\kappa}_1}{\alpha_2} \hspace{0.5mm} T_2 \Big )  \frac{\alpha_{\bar{w}_2}}{\alpha_{\bar{w}_1}} \exp{ \Big ( \frac{ \minus \bar{\gamma} \tilde{\beta}}{2 \alpha_{\bar{w}_2}} (t {+} j) \Big )^2 } |\phi_{\bar{w}}(0,0)|^2_{\tilde{\A}_{\varepsilon_r}}  \\
& \hspace{50mm} \forall (t,j) \in \mbox{dom } \tilde{\phi}
\end{aligned}
\end{equation}
}{
\begin{equation} \label{thrm1_bound}
\begin{aligned} 
V(t,j) & \leq \exp \Big ( \frac{\bar{\kappa}_1}{\alpha_2} \hspace{0.5mm} T_2 \Big ) \Big ( \exp \Big ( \frac{\bar{\kappa}_1}{\alpha_2} \hspace{0.5mm} T_2 \Big ) \Big ( 1 {-} \frac{\bar{\kappa}_2}{\alpha_2} \Big ) \Big )^j V(0,0) \\
& \hspace{5mm} + \bar{\kappa}_2 \exp \Big ( \frac{\bar{\kappa}_1}{\alpha_2} \hspace{0.5mm} T_2 \Big )  \frac{\alpha_{\bar{w}_2}}{\alpha_{\bar{w}_1}} \exp{ \Big ( \minus \frac{\bar{\gamma} \tilde{\beta}}{2 \alpha_{\bar{w}_2}} (t + j) \Big )^2 } |\phi_{\bar{w}}(0,0)|^2_{\tilde{\A}_{\varepsilon_r}}  \\
& \hspace{90mm} \forall (t,j) \in \mbox{dom } \tilde{\phi}
\end{aligned}
\end{equation}
}
\noindent
Now, combining the inequality with (\ref{eqn:lyap_def}) and noting $V(\phi(0,0)) \leq \alpha_2 |\phi(0,0) |_{\tilde{\A}_{\varepsilon}}^2$ one has
\ifbool{two_col}{
\begin{equation}
\begin{aligned} 
|\phi(t,j)|^2_{\tilde{\A}_{\varepsilon}} & \leq \\
& \hspace{-15mm} \frac{\alpha_2}{\alpha_1} |\phi(0,0) |_{\tilde{\A}_{\varepsilon}}^2 \exp \Big ( \frac{\bar{\kappa}_1}{\alpha_2} \hspace{0.5mm} T_2 \Big ) \Big ( \exp \Big ( \frac{\bar{\kappa}_1}{\alpha_2} \hspace{0.5mm} T_2 \Big ) \Big ( 1 {-} \frac{\bar{\kappa}_2}{\alpha_2} \Big ) \Big )^j \\
& \hspace{-12mm} {+} \bar{\kappa}_2 \exp \Big ( \frac{\bar{\kappa}_1}{\alpha_2} \hspace{0.5mm} T_2 \Big )  \frac{\alpha_{\bar{w}_2}}{\alpha_{\bar{w}_1}} \exp{ \Big ( \minus \frac{\bar{\gamma} \tilde{\beta}}{2 \alpha_{\bar{w}_2}} (t {+} j) \Big )^2 } |\phi_{\bar{w}}(0,0)|^2_{\tilde{\A}_{\varepsilon_r}}  \\
& \hspace{50mm} \forall (t,j) \in \mbox{dom } \phi
\end{aligned}
\end{equation}
}{
\begin{equation}
\begin{aligned} 
|\phi(t,j)|^2_{\tilde{\A}_{\varepsilon}} & \leq \alpha_1^{\minus 1} \Big ( \alpha_2 |\phi(0,0) |_{\tilde{\A}_{\varepsilon}}^2 \Big ) \exp \Big ( \frac{\bar{\kappa}_1}{\alpha_2} \hspace{0.5mm} T_2 \Big ) \Big ( \exp \Big ( \frac{\bar{\kappa}_1}{\alpha_2} \hspace{0.5mm} T_2 \Big ) \Big ( 1 {-} \frac{\bar{\kappa}_2}{\alpha_2} \Big ) \Big )^j \\
& \hspace{5mm} + \bar{\kappa}_2 \exp \Big ( \frac{\bar{\kappa}_1}{\alpha_2} \hspace{0.5mm} T_2 \Big )  \frac{\alpha_{\bar{w}_2}}{\alpha_{\bar{w}_1}} \exp{ \Big ( \frac{ \minus \bar{\gamma} \tilde{\beta}}{2 \alpha_{\bar{w}_2}} (t + j) \Big )^2 } |\phi_{\bar{w}}(0,0)|^2_{\tilde{\A}_{\varepsilon_r}}  \\
& \hspace{90mm} \forall (t,j) \in \mbox{dom } \phi
\end{aligned}
\end{equation}
}
\noindent
Then, taking the square root on both sides, one has
\ifbool{two_col}{
\begin{equation}
\begin{aligned} 
|\phi(t,j)|_{\tilde{\A}_{\varepsilon}} & \leq \\
& \hspace{-15mm} \sqrt{\frac{\alpha_2}{\alpha_1}} |\phi(0,0) |_{\tilde{\A}_{\varepsilon}} \exp \Big ( \frac{\bar{\kappa}_1}{2 \alpha_2} \hspace{0.5mm} T_2 \Big ) \Big ( \exp \Big ( \frac{\bar{\kappa}_1}{2 \alpha_2} \hspace{0.5mm} T_2 \Big ) \Big ( 1 {-} \frac{\bar{\kappa}_2}{2 \alpha_2} \Big ) \Big )^j \\
& \hspace{-15mm} {+} \sqrt{\bar{\kappa}_2} \exp \Big ( \frac{\bar{\kappa}_1}{2 \alpha_2} \hspace{0.5mm} T_2 \Big ) \sqrt{ \frac{\alpha_{\bar{w}_2}}{\alpha_{\bar{w}_1}} \exp{ \Big ( \frac{ \minus \bar{\gamma} \tilde{\beta}}{2 \alpha_{\bar{w}_2}} (t {+} j) \Big )^2 } |\phi_{\bar{w}}(0,0)|^2_{\tilde{\A}_{\varepsilon_r}} } \\
& \hspace{50mm} \forall (t,j) \in \mbox{dom } \phi
\end{aligned}
\end{equation}
}{
\begin{equation}
\begin{aligned} 
|\phi(t,j)|_{\tilde{\A}_{\varepsilon}} & \leq \sqrt{\frac{\alpha_2}{\alpha_1}} |\phi(0,0) |_{\tilde{\A}_{\varepsilon}} \exp \Big ( \frac{\bar{\kappa}_1}{2 \alpha_2} \hspace{0.5mm} T_2 \Big ) \Big ( \exp \Big ( \frac{\bar{\kappa}_1}{2 \alpha_2} \hspace{0.5mm} T_2 \Big ) \Big ( 1 {-} \frac{\bar{\kappa}_2}{2 \alpha_2} \Big ) \Big )^j \\
& \hspace{5mm} + \sqrt{\bar{\kappa}_2} \exp \Big ( \frac{\bar{\kappa}_1}{2 \alpha_2} \hspace{0.5mm} T_2 \Big ) \sqrt{ \frac{\alpha_{\bar{w}_2}}{\alpha_{\bar{w}_1}} \exp{ \Big ( \minus \frac{\bar{\gamma} \tilde{\beta}}{2 \alpha_{\bar{w}_2}} (t + j) \Big )^2 } |\phi_{\bar{w}}(0,0)|^2_{\tilde{\A}_{\varepsilon_r}} } \\
& \hspace{90mm} \forall (t,j) \in \mbox{dom } \phi
\end{aligned}
\end{equation}
}
\normalsize
By the given conditions, the set $\tilde{\A}_{\varepsilon}$ is globally exponentially stable and attractive for $\widetilde{\HS}_{\varepsilon}$.  Now, by utilizing Lemmas \ref{lem:equiv_transform} - \ref{lem:GES2}, we can establish global exponential stability to the set $\A_{\varepsilon}$ for $\HS_{\varepsilon}$, in turn we can then make use of Lemmas  \ref{lem:equiv} - \ref{lem:GES1} to then show that the set $\A$ is globally exponentially stable and attractive for $\HS$ in (\ref{eqn:cl_hysys}). }

\section{Robustness to Communication Noise, Clock Drift Perturbations, and Error on ${\bf \sigma}$} \label{sec:robustness}

\ifbool{two_col}{
Under a realistic scenario, it is often the case that the system is subjected to disturbances. In this section, we present results on input-to-state stability (ISS) of the system when it is affected by different types of disturbances. First, we present an ISS result that considers communication noise. Then, we present an ISS result on the parameter estimation sub-system when it is subjected to noise on the internal clock output. Finally, we present an ISS result on noise introduced to the desired clock rate reference $\sigma^*$.\ifbool{conf}{Due to space constraints, the proofs of these results can be found in \cite{12}.}{} 
}{
Under a realistic scenario, it is often the case that the system is subjected to various noise disturbances. Environmental factors can affect the internal clock dynamics and introduce noise to the communication medium in the form of communication delay. In this section we present results on input-to-state stability (ISS) of the system when it is affected by such sources of noise. We will first present an ISS result on the parameter estimation sub-system when it is subjected to noise on the internal clock output, we will then present an ISS result that considers communication noise, last but not least, we will present an ISS result on noise introduced to the desired clock rate reference $\sigma^*$.\ifbool{conf}{Due to space constraints, the proofs of these results can be found in \cite{12}.}{} We will henceforth refer to the following notion of ISS for Hybrid Systems in the presentation of these results, defined as follows:}

\begin{definition} (Input-to-state stability) A hybrid system $\HS$ with input $m$ is input-to-state stable with respect to  a set $\A \subset \reals^n$  if there exist $\beta \in \mathcal{KL}$ and $\kappa \in \mathcal{K}$ such that each solution pair $(\phi, m)$ to $\HS$ satisfies $|\phi(t,j)|_{\A} \leq \mathrm{max} \{\beta(|\phi(0,0)|_{\A},t+j), \kappa(|m|_{\infty}) \}$ for each $(t,j) \in \mathrm{ dom }$ $\phi$.
\end{definition}

\subsection{Robustness to Communication Noise}

We consider the case when the measurements of the timer $\tilde{\tau}_i$ is affected by noise  $m_{e_i} \in \reals$, $i \in \nodes$.  As a result,  the output of each agent is given by  $\tilde{\tau}_i + m_{e_i}$.  In the presence of this noise,  the update law to $\eta_i^+$ in the hybrid controller in  (\ref{protocol}) becomes \ifbool{two_col}{
\small
\begin{equation*}
\begin{aligned}
\eta_i^+ = -\gamma \sum_{k \in \mathcal{N}(i)} (\tilde{\tau}_i - \tilde{\tau}_k) -\gamma \sum_{k \in \mathcal{N}(i)} (m_{e_i} - m_{e_k})
\end{aligned}
\end{equation*}}{
\begin{equation*}
\begin{aligned}
\eta_i^+ & = -\gamma \sum_{k \in \mathcal{N}(i)} (y_i - y_k) \\ 
& = -\gamma \sum_{k \in \mathcal{N}(i)} (\tilde{\tau}_i - \tilde{\tau}_k) -\gamma \sum_{k \in \mathcal{N}(i)} (m_{e_i} - m_{e_k})
\end{aligned}
\end{equation*}}
\normalsize
\noindent
Performing the same change of coordinates, as in the proof of Theorem \ref{thrm1}, we show that $\widetilde{\HS}_{\varepsilon}$ is ISS to communication noise  $m_e := (m_{e_1}, m_{e_2}, \ldots, m_{e_n} ) \in \reals^n$.  Recalling the change of coordinates $\bar{e} = \Tmat^{-1}e$ and $\bar{\eta} = \Tmat^{-1}\eta$, let $\bar{m}_{e} = \Tmat^{-1} m_{e}$. The update law $\bar{\eta}^+$, is given by $\bar{\eta}^+ = (0, -\gamma \bar{\mathcal{L}} \bar{e} -\gamma \bar{\mathcal{L}} \bar{m}_e)$ with $\bar{\eta}_1$ unaffected by the communication noise.

 Using the update law for $\bar{\eta}$ under the effect of $\bar{m}_e$,  we define the perturbed hybrid system $\widetilde{\HS}_m$ with state vector $\chi_m := (\bar{z}_1, \bar{z}_2, \bar{w}_1, \bar{w}_2, \tau) \in \mathcal{X}_{\varepsilon}$, where, again $\bar{z}_1 = (\bar{e}_1, \bar{\eta}_1)$, $\bar{z}_2 = (\bar{e}_2,\ldots,\bar{e}_N, \bar{\eta}_2, \ldots, \bar{\eta}_N)$, $\bar{w}_1 = (\bar{\varepsilon}_{a_1},\bar{\varepsilon}_{\tau_1})$, and $\bar{w}_2 = (\bar{\varepsilon}_{a_2},\ldots, \bar{\varepsilon}_{a_n},\bar{\varepsilon}_{\tau_2},\ldots, \bar{\varepsilon}_{\tau_n})$. Moreover, let $\bar{m}_{\bar{z}_2} = (0, \bar{m}_e)$. The data $( \tilde{C}_m, \tilde{f}_m,  \tilde{D}_m, \tilde{G}_m)$ for the new system $\widetilde{\HS}_m$ is given by
\ifbool{two_col}{$\widetilde{f}_m(\chi_m) := \widetilde{f}_{\varepsilon}(\chi_m)$ for each $\chi_m \in \widetilde{C}_m$ and $\widetilde{G}_m(\chi_m, \bar{m}_{\varepsilon} ) := \widetilde{G}_{\varepsilon}(\chi_m) - 
(0, B_g \bar{m}_{\bar{z}_2}, 0, 0, 0 )$ for each $\chi_m \in \widetilde{D}_m$}{
\begin{equation*} 
\begin{aligned}
\widetilde{f}_m(\chi_m) & := \widetilde{f}_{\varepsilon}(\chi_m) & & \forall \chi_m \in \widetilde{C}_m \\
\widetilde{G}_m(\chi_m, \bar{m}_{\varepsilon} ) & := \widetilde{G}_{\varepsilon}(\chi_m) - \begin{bmatrix}
0 \\ B_g \bar{m}_{\bar{z}_2} \\ 0 \\ 0 \\ 0
\end{bmatrix} & & \forall \chi_m \in \widetilde{D}_m \\
\end{aligned}
\end{equation*}}
\normalsize
\noindent
where $\widetilde{C}_m := \mathcal{X}_{\varepsilon}$, $\widetilde{D}_m := \{ \chi_m \in \mathcal{X}_m : \tau = 0 \}$, and $B_g = \begin{bmatrix}
0 & \gamma \bar{\mathcal{L}}
\end{bmatrix}^{\top}$.

\begin{theorem} \label{thm:iss2}
Given a strongly connected digraph $\DG$, if the parameters $T_2 \geq T_1 > 0$, $\mu > 0$, $h \in \reals$, $\gamma > 0$, and positive definite symmetric matrices $P_1$, $P_2$, and $P_3$ are such that (\ref{cond:phil_1}) and (\ref{thrm_cond1}) hold, the hybrid system $\widetilde{\HS}_m$ with input $\bar{m}_e$ is ISS with respect to $\widetilde{\A}_{\varepsilon}$ in (\ref{set:A_tilde}). \ifbool{rep}{
Furthermore, for each $(t,j) \in \mbox{dom } \phi$ every solution $\phi$ to $\widetilde{\HS}_m$ satisfies 
\begin{equation*}
\begin{aligned} 
|\phi(t,j)|_{\tilde{\A}_{\varepsilon}} \leq & \sqrt{\frac{\alpha_2}{\alpha_1}} |\phi(0,0) |_{\tilde{\A}_{\varepsilon}} \exp \Big ( \frac{\bar{\kappa}_1}{2 \alpha_2} \hspace{0.5mm} T_2 \Big ) \Big ( \exp \Big ( \frac{\bar{\kappa}_1}{2 \alpha_2} \hspace{0.5mm} T_2 \Big ) \Big ( 1 {-} \frac{\bar{\kappa}_2}{2 \alpha_2} \Big ) \Big )^j \\
& \hspace{0mm} {+} \sqrt{\bar{\kappa}_2} \exp \Big ( \frac{\bar{\kappa}_1}{2 \alpha_2} \hspace{0.5mm} T_2 \Big ) \sqrt{ \frac{\alpha_{\bar{w}_2}}{\alpha_{\bar{w}_1}} \exp{ \Big ( \frac{ \minus \bar{\gamma} \tilde{\beta}}{2 \alpha_{\bar{w}_2}} (t {+} j) \Big )^2 } |\phi_{\bar{w}}(0,0)|^2_{\tilde{\A}_{\varepsilon_r}} } \\
& \hspace{7mm} + \sqrt{\tilde{\kappa}_{\bar{m}_2}} \exp \Big (\frac{\kappa}{4 \epsilon_2} T_2 \Big ) \sqrt{\underset{(t,j) \in \mbox{\footnotesize \rm dom} \hspace{1mm} \phi}{\mathrm{sup}} |\bar{m}_{\bar{z}_2}(t,j)|^2}
\end{aligned}
\end{equation*} 
\EndNew
}{
Furthermore, for each $(t,j) \in \mbox{dom } \phi$ every solution $\phi$ to $\widetilde{\HS}_m$ satisfies $|\phi(t,j)|_{\tilde{\A}_{\varepsilon}} \leq$
\begin{equation*}
\begin{aligned} 
& \sqrt{\frac{\alpha_2}{\alpha_1}} |\phi(0,0) |_{\tilde{\A}_{\varepsilon}} \exp \Big ( \frac{\bar{\kappa}_1}{2 \alpha_2} \hspace{0.5mm} T_2 \Big ) \Big ( \exp \Big ( \frac{\bar{\kappa}_1}{2 \alpha_2} \hspace{0.5mm} T_2 \Big ) \Big ( 1 {-} \frac{\bar{\kappa}_2}{2 \alpha_2} \Big ) \Big )^j \\
& \hspace{0mm} {+} \sqrt{\bar{\kappa}_2} \exp \Big ( \frac{\bar{\kappa}_1}{2 \alpha_2} \hspace{0.5mm} T_2 \Big ) \sqrt{ \frac{\alpha_{\bar{w}_2}}{\alpha_{\bar{w}_1}} \exp{ \Big ( \frac{ \minus \bar{\gamma} \tilde{\beta}}{2 \alpha_{\bar{w}_2}} (t {+} j) \Big )^2 } |\phi_{\bar{w}}(0,0)|^2_{\tilde{\A}_{\varepsilon_r}} } \\
& \hspace{7mm} + \sqrt{\tilde{\kappa}_{\bar{m}_2}} \exp \Big (\frac{\kappa}{4 \epsilon_2} T_2 \Big ) \sqrt{\underset{(t,j) \in \mbox{\footnotesize \rm dom} \hspace{1mm} \phi}{\mathrm{sup}} |\bar{m}_{\bar{z}_2}(t,j)|^2}
\end{aligned}
\end{equation*} 
}
\end{theorem}


\ifbool{two_col}{The proof of this result utilizes a Lyapunov analysis using the function candidate $V$ in (\ref{eqn:lyap}). Since the disturbance is present during jumps, we show that $V$ can be upper bounded resulting in a bounded disturbance in $V$ when evaluated along a given solution to $\widetilde{\HS}_{m_{\sigma}}$; see \cite{12} for more details.}
{
\begin{proof}
Consider the same Lyapunov function candidate $V(\chi_m) = V_1(\chi_m) + V_2(\chi_m) + V_{\varepsilon_r} (\chi_m)$ from the proof of Theorem \ref{thrm1}. During flows, there is no contribution from the perturbation thus the derivative of $V$ is unchanged from the proof of Theorem \ref{thrm1}. Thus, one has
\ifbool{rep}{
\begin{equation*}
\begin{aligned}
\langle \nabla V(\chi_m), \tilde{f}(\chi_m) \rangle & \leq 2 \bar{z}_2^{\top} \big ( \exp{A^{\top}_{f_2} \tau } P \exp{A_{f_2} \tau} \big ) B_{f_2}  \bar{w}_2 + \bar{w}_1^{\top} (P_1 A_{f_3} + A_{f_3}^{\top} P_1) \bar{w}_1 + \bar{w}_2^{\top} (P_2 A_{f_4} + A_{f_4}^{\top} P_2) \bar{w}_2 \\ 
\end{aligned}
\end{equation*}}{
\begin{equation*}
\begin{aligned}
\langle \nabla V(\chi_m), \tilde{f}(\chi_m) \rangle & \leq 2 \bar{z}_2^{\top} \big ( \exp{A^{\top}_{f_2} \tau } P \exp{A_{f_2} \tau} \big ) B_{f_2}  \bar{w}_2 \\ 
& \hspace{3mm} + \bar{w}_1^{\top} (P_1 A_{f_3} + A_{f_3}^{\top} P_1) \bar{w}_1 \\
& \hspace{5mm} + \bar{w}_2^{\top} (P_2 A_{f_4} + A_{f_4}^{\top} P_2) \bar{w}_2
\end{aligned}
\end{equation*}}
\noindent
then by following the same notions of the proof in Theorem \ref{thrm1}, one has
\begin{equation*}
\begin{aligned}
\langle \nabla V(\chi_m), \tilde{f}(\chi_m) \rangle 
& \leq \frac{\kappa_1}{2 \epsilon} | \bar{z}_2 
|^2 - \beta_1 | \bar{w}_1|^2 + \Big ( \frac{\kappa_1 \epsilon}{2} - \beta_2 \Big ) | \bar{w}_2|^2 \\ 
& \leq \bar{\kappa}_1 \big (|\bar{z}_2|^2 + |\bar{w}_1|^2 + |\bar{w}_2|^2 \big ) \\
& \leq \bar{\kappa}_1 V(x)
\end{aligned}
\end{equation*}
\noindent
where $\bar{\kappa}_1 = \max \Big \{ \frac{\kappa_1}{2 \epsilon} , \Big ( \frac{\kappa_1 \epsilon}{2} - \beta_2 \Big ) \Big \}$ and $\varepsilon > 0$. At jumps, triggered when $\tau = 0$, one has, for each $\chi_m \in \tilde{D}_m \setminus \tilde{\A}_{\varepsilon}$ and $g \in \tilde{G}_m(\chi_m)$
\ifbool{rep}{
\begin{equation} \label{thrm_iss_jumps}
\begin{aligned}
V(g) {-} V(\chi_m) & \leq {-} \bar{\eta}_1^2 + (A_{g_2} \bar{z}_2 - B_g \bar{m}_{\bar{z}_2})^{\top} Q (A_{g_2} \bar{z}_2 - B_g \bar{m}_{\bar{z}_2}) - \bar{z}_2^{\top} P_1 \bar{z}_2 \\
& \leq {-} \bar{\eta}_1^2 + (A_{g_2} \bar{z}_2)^{\top} \exp{A^{\top}_{f_2} \tau } P_1 \exp{A_{f_2} \tau} (A_{g_2} \bar{z}_2) - 2 (B_g \bar{m}_{\bar{z}_2})^{\top} \exp{A^{\top}_{f_2} \tau } P_1 \exp{A_{f_2} \tau} (A_{g_2} \bar{z}_2) \\ 
& \hspace{60mm} + (B_g \bar{m}_{\bar{z}_2})^{\top} \exp{A^{\top}_{f_2} \tau } P_1 \exp{A_{f_2} \tau} (B_g \bar{m}_{\bar{z}_2}) {-}  \bar{z}_2^{\top} P_1 \bar{z}_2 \\ 
\end{aligned}
\end{equation}}{
\begin{equation} \label{thrm_iss_jumps}
\begin{aligned}
V(g) {\minus} V(\chi_m) & \leq {\minus} \bar{\eta}_1^2 {+} (A_{g_2} \bar{z}_2 {-} B_g \bar{m}_{\bar{z}_2})^{\top} Q (A_{g_2} \bar{z}_2 {-} B_g \bar{m}_{\bar{z}_2}) \\
& \hspace{5mm} - \bar{z}_2^{\top} P_1 \bar{z}_2 \\
& \leq {\minus} \bar{\eta}_1^2 + (A_{g_2} \bar{z}_2)^{\top} \exp{A^{\top}_{f_2} \tau } P_1 \exp{A_{f_2} \tau} (A_{g_2} \bar{z}_2) \\
& \hspace{5mm} {-} 2 (B_g \bar{m}_{\bar{z}_2})^{\top} \exp{A^{\top}_{f_2} \tau } P_1 \exp{A_{f_2} \tau} (A_{g_2} \bar{z}_2) \\ 
& \hspace{10mm} + (B_g \bar{m}_{\bar{z}_2})^{\top} \exp{A^{\top}_{f_2} \tau } P_1 \exp{A_{f_2} \tau} (B_g \bar{m}_{\bar{z}_2}) \\
& \hspace{15mm} {-}  \bar{z}_2^{\top} P_1 \bar{z}_2 \\ 
\end{aligned}
\end{equation}}
\noindent
From \ref{cond:phil_1} and the proof in Theorem \ref{thrm1}, there exists a scalar $\kappa_2$ such that $$\bar{z}_2^{\top} (A_{g_2}^{\top} \exp{A^{\top}_{f_2} v } P_1 \exp{A_{f_2} v} A_{g_2} - P_1) \bar{z}_2 \leq - \kappa_2 \bar{z}_2^{\top} \bar{z}_2$$
\noindent
leading to
\ifbool{rep}{
\begin{equation} \label{thrm_iss_jumps1}
\begin{aligned}
V(g) {-} V(\chi_m) & \leq - \bar{\eta}_1^2 -\kappa_2 \bar{z}_2^{\top} \bar{z}_2 - 2 (B_g \bar{m}_{\bar{z}_2})^{\top} \exp{A^{\top}_{f_2} \tau } P_1 \exp{A_{f_2} \tau} (A_{g_2} \bar{z}_2) + (B_g \bar{m}_{\bar{z}_2})^{\top} \exp{A^{\top}_{f_2} \tau } P_1 \exp{A_{f_2} \tau} (B_g \bar{m}_{\bar{z}_2})
\end{aligned}
\end{equation}}{
\begin{equation} \label{thrm_iss_jumps1}
\begin{aligned}
V(g) {-} V(\chi_m) & \leq - \bar{\eta}_1^2 -\kappa_2 \bar{z}_2^{\top} \bar{z}_2 \\
& \hspace{5mm} - 2 (B_g \bar{m}_{\bar{z}_2})^{\top} \exp{A^{\top}_{f_2} \tau } P_1 \exp{A_{f_2} \tau} (A_{g_2} \bar{z}_2) \\
& \hspace{10mm} + (B_g \bar{m}_{\bar{z}_2})^{\top} \exp{A^{\top}_{f_2} \tau } P_1 \exp{A_{f_2} \tau} (B_g \bar{m}_{\bar{z}_2})
\end{aligned}
\end{equation}}
\noindent
Let $Q = \exp{A^{\top}_{f_2} \tau } P_1 \exp{A_{f_2} \tau}$, then applying Young's inequality on the third term such that 
\ifbool{rep}{
\begin{equation*}
\begin{aligned}
\bar{m}_{\bar{z}_2}^{\top} B_g^{\top} Q A_{g_2} \bar{z}_2 & \leq \frac{1}{2 \epsilon_2} \Big ( \bar{m}_{\bar{z}_2}^{\top} B_g^{\top} Q A_{g_2} \Big )^{\top} \Big ( \bar{m}_{\bar{z}_2}^{\top} B_g^{\top} Q A_{g_2} \Big ) + \frac{\epsilon_2}{2} \bar{z}_2^{\top} \bar{z}_2 \\
& \leq \frac{1}{2 \epsilon_2} \Big | \big ( B_g^{\top} Q A_{g_2} \big ) \big ( B_g^{\top} Q A_{g_2} \big )^{\top} \Big | \bar{m}_{\bar{z}_2}^{\top} \bar{m}_{\bar{z}_2} + \frac{\epsilon_2}{2} \bar{z}_2^{\top} \bar{z}_2 \\
\end{aligned}
\end{equation*}}{
\begin{equation*}
\begin{aligned}
\bar{m}_{\bar{z}_2}^{\top} B_g^{\top} Q A_{g_2} \bar{z}_2 & \leq \frac{1}{2 \epsilon_2} \Big ( \bar{m}_{\bar{z}_2}^{\top} B_g^{\top} Q A_{g_2} \Big )^{\top} \Big ( \bar{m}_{\bar{z}_2}^{\top} B_g^{\top} Q A_{g_2} \Big ) \\ 
& \hspace{5mm} + \frac{\epsilon_2}{2} \bar{z}_2^{\top} \bar{z}_2 \\
& \leq \frac{1}{2 \epsilon_2} \Big | \big ( B_g^{\top} Q A_{g_2} \big ) \big ( B_g^{\top} Q A_{g_2} \big )^{\top} \Big | \bar{m}_{\bar{z}_2}^{\top} \bar{m}_{\bar{z}_2} \\
& \hspace{5mm} + \frac{\epsilon_2}{2} \bar{z}_2^{\top} \bar{z}_2 \\
\end{aligned}
\end{equation*}}
\noindent
where $\epsilon_2 > 0$, we then have
\ifbool{rep}{
\begin{equation} \label{thrm_iss_jumps2}
\begin{aligned}
V(g) \minus V(\chi_m) & \leq - \bar{\eta}_1^2 -\kappa_2 \bar{z}_2^{\top} \bar{z}_2 {-} \Big ( {\small \frac{1}{2 \epsilon_2}} \big | \big ( B_g^{\top} Q A_{g_2} \big ) \big ( B_g^{\top} Q A_{g_2} \big )^{\top} \big | \bar{m}_{\bar{z}_2}^{\top} \bar{m}_{\bar{z}_2} + \frac{\epsilon_2}{2} \bar{z}_2^{\top} \bar{z}_2 \Big ) \\ 
& \hspace{15mm} + \bar{m}_{\bar{z}_2} B_g^{\top} Q B_g \bar{m}_{\bar{z}_2} \\
& \leq - \bar{\eta}_1^2 - \big ( \kappa_2 + \frac{\epsilon_2}{2} \big ) \bar{z}_2^{\top} \bar{z}_2 {+} \big ( | B_g^{\top} Q B_g | \minus \frac{1}{2 \epsilon_2} | ( B_g^{\top} Q A_{g_2} )( B_g^{\top} Q A_{g_2} )^{\top} | \big ) \bar{m}_{\bar{z}_2}^{\top}  \bar{m}_{\bar{z}_2} 
\end{aligned}
\end{equation}}{
\begin{equation} \label{thrm_iss_jumps2}
\begin{aligned}
V(g) \minus V(\chi_m) & \leq - \bar{\eta}_1^2 -\kappa_2 \bar{z}_2^{\top} \bar{z}_2 \\ & \hspace{5mm} {-} \Big ( {\small \frac{1}{2 \epsilon_2}} \big | \big ( B_g^{\top} Q A_{g_2} \big ) \big ( B_g^{\top} Q A_{g_2} \big )^{\top} \big | \bar{m}_{\bar{z}_2}^{\top} \bar{m}_{\bar{z}_2}
\\ & \hspace{10mm} + \frac{\epsilon_2}{2} \bar{z}_2^{\top} \bar{z}_2 \Big ) + \bar{m}_{\bar{z}_2} B_g^{\top} Q B_g \bar{m}_{\bar{z}_2} \\
& \leq - \bar{\eta}_1^2 - \big ( \kappa_2 + \frac{\epsilon_2}{2} \big ) \bar{z}_2^{\top} \bar{z}_2 {+} \big ( | B_g^{\top} Q B_g | \\
& \hspace{5mm} \minus \frac{1}{2 \epsilon_2} | ( B_g^{\top} Q A_{g_2} )( B_g^{\top} Q A_{g_2} )^{\top} | \big ) \bar{m}_{\bar{z}_2}^{\top}  \bar{m}_{\bar{z}_2} 
\end{aligned}
\end{equation}}
\noindent
by noting $|A_{g_2}|, |B_g| \leq \gamma \lambda_{max}({\bar{\mathcal{L}}})$ let $$\kappa_{\bar{m}_2} = \big ( \lambda_{max}({\bar{\mathcal{L}}}) \big )^2 \underset{v \in [0,T_2]}{\mathrm{max}} \Big \{ \lambda_{max} \big (\exp{A^{\top}_{f_2} v } P_1 \exp{A_{f_2} v } \big ) \Big \}$$ then we let $\epsilon_2 = \kappa_2$ and
\ifbool{rep}{
\begin{equation*} 
\begin{aligned}
V(g) - V(\chi_m) & \leq - \bar{\eta}_1^2 - \big ( \kappa_2 + \frac{\kappa_2}{2} \big ) \bar{z}_2^{\top} \bar{z}_2 +  \big ( \gamma^2 \kappa_{\bar{m}_2} - \frac{1}{2 \kappa_2} \gamma^4 \kappa_{\bar{m}_2}^2 \big ) \bar{m}_{\bar{z}_2}^{\top}  \bar{m}_{\bar{z}_2} 
\end{aligned}
\end{equation*}}{\begin{equation*} 
\begin{aligned}
V(g) - V(\chi_m) & \leq - \bar{\eta}_1^2 - \big ( \kappa_2 + \frac{\kappa_2}{2} \big ) \bar{z}_2^{\top} \bar{z}_2 \\
& \hspace{5mm} +  \big ( \gamma^2 \kappa_{\bar{m}_2} - \frac{1}{2 \kappa_2} \gamma^4 \kappa_{\bar{m}_2}^2 \big ) \bar{m}_{\bar{z}_2}^{\top}  \bar{m}_{\bar{z}_2} 
\end{aligned}
\end{equation*}}
\noindent
now let $\tilde{\kappa}_{\bar{m}_2} = \big ( \gamma^2 \kappa_{\bar{m}_2} - \frac{1}{2 \kappa_2} \gamma^4 \kappa_{\bar{m}_2}^2 \big )$ then at jumps one has
\begin{equation} \label{thrm_iss_jumps3}
\begin{aligned}
V(g) - V(\chi_m) & \leq - \bar{\kappa}_2 (|\bar{\eta}_1|^2 + |\bar{z}_2|^2) + \tilde{\kappa}_{\bar{m}_2} |\bar{m}_{\bar{z}_2}|^2 
\end{aligned}
\end{equation}
\noindent
where $\bar{\kappa}_2 = \max \big \{1, \frac{3 \kappa_2}{2} \big \}$.   Now, recall from (\ref{thrm_zdef}) in the proof of Theorem \ref{thrm1},
\begin{equation}
\begin{aligned}
-(|\bar{\eta}_1|^2 + |\bar{z}_2|^2) \leq - \frac{1}{\alpha_2} V(\chi_{\varepsilon}) + |\bar{w}|^2
\end{aligned}
\end{equation}
\noindent
by then plugging (\ref{thrm_zdef}) in to (\ref{thrm_iss_jumps3}) one has
\begin{equation*}
\begin{aligned}
V(g) - V(\chi_m) & \leq \frac{3 \kappa_2}{2} \Big ( \minus \frac{1}{\alpha_2} V(\chi_{\varepsilon}) + |\bar{w}|^2 \Big ) + \tilde{\kappa}_{\bar{m}_2} |\bar{m}_{\bar{z}_2}|^2 \\
& \leq \minus \frac{3 \kappa_2}{2 \alpha_2} V(\chi_{\varepsilon}) + \frac{3 \kappa_2}{2} |\bar{w}|^2 + \tilde{\kappa}_{\bar{m}_2} |\bar{m}_{\bar{z}_2}|^2 
\end{aligned}
\end{equation*}
\noindent
then at jumps one has
\begin{equation*}
\begin{aligned}
V(g) & \leq \Big ( 1 - \frac{3 \kappa_2}{2 \alpha_2} \Big ) V(\chi_{\varepsilon}) + \frac{3 \kappa_2}{2} |\bar{w}|^2 + \tilde{\kappa}_{\bar{m}_2} |\bar{m}_{\bar{z}_2}|^2 
\end{aligned}
\end{equation*}
\noindent
Noting $\langle \nabla V(\chi_{\varepsilon}), \tilde{f}(\chi_{\varepsilon}) \rangle  \leq \bar{\kappa} V(\chi_{\varepsilon})$, one can then pick a solution with initial conditions $\phi(0,0) \in \tilde{C}_m \cup \tilde{D}_m$ and find the trajectory of $V(t,j)$ is bounded as follows
\begin{equation*}
\begin{aligned}
V(t,j) & \leq \exp \big ( \bar{\kappa} \hspace{0.5mm} T_2 \big ) \Big ( \exp \big ( \bar{\kappa} \hspace{0.5mm} T_2 \big ) \Big ( 1 - \frac{3 \kappa_2}{2 \alpha_2} \Big ) \Big )^j V(0,0) \\
& \hspace{5mm} + \frac{3 \kappa_2}{2} \exp \big ( \bar{\kappa} \hspace{0.5mm} T_2 \big ) \mathrm{sup}_{(t,j) \in \mbox{\footnotesize dom} \phi} |\bar{w}(t,j)|^2 \\
& \hspace{10mm} + \tilde{\kappa}_{\bar{m}_2} \exp \Big (\frac{\kappa}{2 \epsilon_2} T_2 \Big ) \mathrm{sup}_{(t,j) \in \mbox{\footnotesize dom} \phi} |\bar{m}_{\bar{z}_2}|^2
\end{aligned}
\end{equation*}
\end{proof}}

\subsection{Robustness to Perturbations on Internal Clock Drift}

In this section, we consider a disturbance $m_{\tau^*_i} \in \reals$, $i \in \nodes$ added to the output of the internal clock. Let $y^{\tau^*}_i := \tau^*_i + m_{\tau^*_i}$, $i \in \nodes$, define the perturbed internal clock output. Then the dynamics of the original estimation system in (\ref{eqn:cl_hysys}) under this disturbance becomes 
\begin{equation}
\begin{aligned}
& \dot{\hat{\tau}}_i = \hat{a}_i {-} (\hat{\tau}_i {-} y^{\tau^*}_i), & & \dot{\hat{a}}_i = \minus \mu ( \hat{\tau}_i {-} y^{\tau^*}_i ) & & \tau \in [0, T_2] \\
& \hat{\tau}^+_i = \hat{\tau}_i, &  &\hat{a}^+_i = \hat{a}_i & & \tau = 0
\end{aligned}
\end{equation}
\noindent
In error coordinates $\varepsilon_{\hat{a}_i} = a_i - \hat{a}_i$, $\varepsilon_{\tau_i} = \hat{\tau}_i - \tau_i^*$, this leads to 
\begin{align*}
& \dot{\varepsilon}_{\tau_i} = - \varepsilon_{\tau_i} - \varepsilon_{a_i} + m_{\tau^*_i}, \hspace{1mm} \dot{\varepsilon}_{\hat{a}_i} = \mu \varepsilon_{\tau_i} - \mu m_{\tau^*_i}  & \tau & \in [0, T_2] \\
& \varepsilon_{\tau_i}^+ = \varepsilon_{\tau_i}, \hspace{23mm} \varepsilon_{a_i}^+ = \varepsilon_{a_i} \hspace{15mm} & \tau & = 0
\end{align*}
Similar to the result presented in Proposition \ref{prop:est_param}, for the estimation sub-system we will consider the same reduction $\widetilde{\HS}_{\varepsilon_r}$ that now captures the perturbation. Recall the coordinate transformations $\bar{\varepsilon}_a = \Tmat^{-1}\varepsilon_a$ and $\bar{\varepsilon}_{\tau} = \Tmat^{-1}\varepsilon_{\tau}$ for the respective internal clock and parameter estimation errors. Moreover, recall $\bar{w} = (\bar{w}_1,\bar{w}_2)$ where $\bar{w}_1 = (\bar{\varepsilon}_{a_1},\bar{\varepsilon}_{\tau_1})$ and $\bar{w}_2 = (\bar{\varepsilon}_{a_2},\ldots, \bar{\varepsilon}_{a_n},\bar{\varepsilon}_{\tau_2},\ldots, \bar{\varepsilon}_{\tau_n})$. Let $\bar{m}_{\tau^*} = \Tmat^{-1}m_{\tau^*}$ and $\bar{q} = (\bar{q}_1, \bar{q}_2)$ where $\bar{q}_1 = (\bar{m}_{\tau^*_1}, \bar{m}_{\tau^*_1})$ and $\bar{q}_2 = (\bar{m}_{\tau^*_2}, \ldots, \bar{m}_{\tau^*_n}, \bar{m}_{\tau^*_2}, \ldots, \bar{m}_{\tau^*_n})$. Now, consider the reduced coordinates $\chi_{m_r} := (\bar{w}_1, \bar{w}_2, \tau) \in \reals^n \times \reals^n \times [0,T_2] =: \mathcal{X}_{\varepsilon}$. The data of this reduced system is given by $\widetilde{\HS}_{m_r} = ( \widetilde{C}_{\varepsilon}, \widetilde{f}_{\varepsilon},  \widetilde{D}_{\varepsilon}, \widetilde{G}_{\varepsilon})$ where
\ifbool{two_col}{$\widetilde{f}_{m_r}(\chi_{m_r}, \bar{q}) := \widetilde{f}_{\varepsilon_r}(\chi_{m_r}) + (B_{m_1} \bar{q}_1, B_{m_2} \bar{q}_2, 0)$ for each $\chi_{m_r} \in \widetilde{C}_{m_r}$ and $\widetilde{G}_{m_r}(\chi_{m_r}) := ( \bar{w}_1, \bar{w}_2, [T_1,T_2])$ for each $\chi_{m_r} \in \widetilde{D}_{m_r}$}{
\begin{equation*}
\begin{aligned}
\widetilde{f}_{m r}(\chi_{m r}) & := \begin{bmatrix} A_{f_3} \bar{w}_1 \\ A_{f_4} \bar{w}_2 \\  -1
\end{bmatrix} + \begin{bmatrix}
B_{m_1} \bar{q}_1 \\ B_{m_2} \bar{q}_2 \\ 0
\end{bmatrix} & & \forall \chi_{m r} \in \widetilde{C}_{m r} \\
\widetilde{G}_{m r}(\chi_{m r}) & := \begin{bmatrix}
\bar{w}_1 \\ \bar{w}_2 \\ [T_1,T_2]
\end{bmatrix} & & \forall \chi_{m r} \in \widetilde{D}_{m r} \\
\end{aligned}
\end{equation*}
}
\normalsize
\noindent
where $\tilde{C}_{m_r} := \mathcal{X}_{\varepsilon}$, $\tilde{D}_{m_r} := \{ \chi_m \in \mathcal{X}_{\varepsilon} : \tau = 0 \}$, and 
\ifbool{two_col}{
$B_{m_1} = (
[\mu, 0], [0, 1] )$, $B_{m_2} = ([
\mu I, 0], [0, I])$.}{
$$B_{m_1} = \begin{bmatrix}
\mu & 0 \\ 0 & 1 
\end{bmatrix}, \hspace{1mm} B_{m_2} = \begin{bmatrix}
\mu I & 0 \\ 0 & I 
\end{bmatrix}$$}

\ifbool{two_col}{
\begin{theorem}
If there exists a positive scalar $\mu$ and positive definite symmetric matrices $P_2$, $P_3$ such that the conditions in (\ref{eqn:lyap1}) hold, the hybrid system $\widetilde{\HS}_{m r}$ with input $\bar{m}_{\tau^*}$ is ISS with respect to $\widetilde{\A}_{\varepsilon_r}$  given in (\ref{set:est_param}). 
\end{theorem}

The proof of this result is established by picking a solution to the model reduction $\widetilde{\HS}_{\varepsilon_r}$, integrating the disturbance that is treated as an input to the system,  and then bounding the integral.  A proof of this result can be found in \cite{12}. 
}{
\begin{theorem}
If there exists a positive scalar $\mu$ and positive definite symmetric matrices $P_2$, $P_3$ such that (\ref{eqn:lyap1}) and (\ref{eqn:lyap2}) hold, the hybrid system $\widetilde{\HS}_{m r}$ with input $\bar{m}_{\tau^*}$ is ISS with respect to $\widetilde{\A}_{\varepsilon_r}$.
\end{theorem}}

\ifbool{two_col}{}{\begin{proof}
Since the matrices $A_{f_3}$ and $A_{f_4}$ are Hurwitz and the states $\bar{w}_1$ and $\bar{w}_2$ do not jump, we can estimation system as a continuous time system and write the solution explicitly for the states $\bar{w}_1$ and $\bar{w}_2$. 
\begin{equation}
\begin{aligned}
\phi_{\bar{w}_1} (t,j) & = \exp(A_{f_3}(t - 0)) \phi_{\bar{w}_1}(0,0) \\ 
& \hspace{5mm} + \int_0^{t} \exp(A_{f_3}(t - s)) B_{m_1} \bar{q}_1(s) ds 
\end{aligned}
\end{equation}
\noindent
and
\begin{equation}
\begin{aligned}
\phi_{\bar{w}_2} (t,j) & = \exp(A_{f_4}(t - 0)) \phi_{\bar{w}_2}(0,0) \\ 
& \hspace{5mm} + \int_0^{t} \exp(A_{f_4}(t - s)) B_{m_2} \bar{q}_2(s) ds
\end{aligned} 
\end{equation}
\noindent
then by bounding $|\exp(A_{f_3}(t - 0))| {\leq} \rho_1 \exp{-\lambda_1(t-0)}$ and $|\exp(A_{f_4}(t - 0))| {\leq} \rho_2 \exp{-\lambda_2(t-0)}$ one has
\begin{equation}
\begin{aligned}
|\phi_{\bar{w}_1} (t,j)| & \leq \rho_1 \exp{-\lambda_1(t-0)} |\phi_{\bar{w}_1} (0,0)| \\ 
& \hspace{5mm} + \int_0^{t} \rho_1 \exp{-\lambda_1(t-s)} |B_{m_1}| |\bar{q}_1(s)| ds \\
& \leq \rho_1 \exp{-\lambda_1(t-0)} |\phi_{\bar{w}_1} (0,0)| + \frac{\rho_1 |B_{m_1}|}{\lambda_1} \underset{0 \leq s \leq t}{\mathrm{sup}} |\bar{q}_2(s)|
\end{aligned}
\end{equation}
\noindent
and
\begin{equation}
\begin{aligned}
|\phi_{\bar{w}_2} (t,j)| & \leq \rho_2 \exp{-\lambda_2(t-0)} |\phi_{\bar{w}_2} (0,0)| \\ 
& \hspace{5mm} + \int_0^{t} \rho_2 \exp{-\lambda_2(t-s)} |B_{m_2}| |\bar{q}_2(s)| ds \\
& \leq \rho_2 \exp{-\lambda_2(t-0)} |\phi_{\bar{w}_2} (0,0)| + \frac{\rho_2 |B_{m_2}|}{\lambda_2} \underset{0 \leq s \leq t}{\mathrm{sup}} |\bar{q}_2(s)|
\end{aligned}
\end{equation}
\end{proof}}


\subsection{Robustness to Error on ${\bf \sigma}$}

 In this section, we consider a disturbance on $\sigma^*$ to capture the scenario where $\sigma^*$ is not precisely known, i.e., $\sigma_i \neq \sigma^*$. Let $\varepsilon_{\sigma_i} = \sigma_i - \sigma^*$ represent the error between the injected and the ideal clock rate. 
 \ifbool{conf}{}{
 Treating $\varepsilon_{\sigma}$ as a perturbation to the system $\HS_{\varepsilon}$, one has
\ifbool{two_col}{$\dot{x}_{\varepsilon} = f_{\varepsilon}(x_{\varepsilon}) + 
(\varepsilon_{\sigma}, 0, 0, 0, 0)$ for each $x_{\varepsilon} \in C_{\varepsilon}$ and $x_{\varepsilon}^+ \in \big (e, -\gamma \mathcal{L} e, \varepsilon_a, \varepsilon_{\tau}, [T_1,T_2] \big )$ for each $x_{\varepsilon} \in D_{\varepsilon}$.
}{
\begin{equation*}
\begin{aligned}
& \dot{x}_{\varepsilon} = \begin{bmatrix}
\eta + \varepsilon_a  \\ h \eta \\ \mu \varepsilon_{\tau} \\ - \varepsilon_{\tau} - \varepsilon_a \\ -1 \\
\end{bmatrix} + \begin{bmatrix}
\varepsilon_{\sigma} \\ 0 \\ 0 \\ 0 \\ 0
\end{bmatrix} & & \forall x_{\varepsilon} \in C_{\varepsilon} \\
& x_{\varepsilon}^+ \in \big (e, -\gamma \mathcal{L} e, \varepsilon_a, \varepsilon_{\tau}, [T_1,T_2] \big ) & & \forall x_{\varepsilon} \in D_{\varepsilon}
\end{aligned}
\end{equation*}
}
} 
To show how the perturbation affects $\widetilde{\HS}_{\varepsilon}$, let $\bar{\varepsilon}_{\sigma} = \Tmat^{-1} \varepsilon_{\sigma}$, then let $\bar{m}_{\sigma} = (\bar{m}_{\sigma_1},\bar{m}_{\sigma_2})$ where $\bar{m}_{\sigma_1} = \bar{\varepsilon}_{\sigma_1}$ and $\bar{m}_{\sigma_2} = (\bar{\varepsilon}_{\sigma_2}, \ldots , \bar{\varepsilon}_{\sigma_n})$.

We define this perturbed hybrid system $\widetilde{\HS}_{m_{\sigma}}$ with state vector $\chi_{m_{\sigma}} := (\bar{z}_1, \bar{z}_2, \bar{w}_1, \bar{w}_2, \tau) \in \mathcal{X}_{\varepsilon}$.  Its dynamics are given by the new system $\widetilde{\HS}_{m_{\sigma}} = ( \tilde{C}_{m_{\sigma}}, \tilde{f}_{m_{\sigma}},  \tilde{D}_{m_{\sigma}}, \tilde{G}_{m_{\sigma}})$ with data $\tilde{f}_{m_{\sigma}}(\chi_{m_{\sigma}})$ for each $\chi_{m_{\sigma}} \in \tilde{C}_{m_{\sigma}}$ $:= \mathcal{X}_{\varepsilon}$ and $\tilde{G}_{m_{\sigma}}(\chi_{m_{\sigma}})$ for each $\chi_{m_{\sigma}} \in \tilde{D}_{m_{\sigma}}$ $:= \{ \chi_{m_{\sigma}} \in \mathcal{X}_{\varepsilon} : \tau = 0 \}$ where 
\ifbool{two_col}{ $\widetilde{f}_{m_{\sigma}}(\chi_{m_{\sigma}}, \bar{m}_{\sigma}) :=$ $\widetilde{f}_{\varepsilon}(\chi_{m_{\sigma}})$ $+ \big ( \bar{m}_{\sigma_1}, \bar{m}_{\sigma_2}, 0, 0, 0 \big )$ and $\widetilde{G}_{m_{\sigma}}(\chi_{m_{\sigma}}) :=$ $\widetilde{G}_{\varepsilon}(\chi_{m_{\sigma}})$ leading to the following result.
}{
\small
\begin{equation*}
\begin{aligned}
\hspace{-1mm} \widetilde{f}_{m_{\sigma}}(\chi_{m_{\sigma}}) & := \begin{bmatrix}
A_{f_1} \bar{z}_1 + B_{f_1} \bar{w}_1 \\ A_{f_2} \bar{z}_2 + B_{f_2} \bar{w}_2 \\ A_{f_3} \bar{w}_1 \\ A_{f_4} \bar{w}_2 \\  -1
\end{bmatrix} + \begin{bmatrix}
\bar{m}_{\sigma_1} \\ \bar{m}_{\sigma_2} \\ 0 \\ 0 \\ 0
\end{bmatrix}  \\
\hspace{-1mm} \widetilde{G}_{m_{\sigma}}(\chi_{m_{\sigma}}) & := \begin{bmatrix}
[A_{g_1} \bar{z}_1]^\top, & [A_{g_2}^\top \bar{z}_2]^\top, & \bar{w}_1^\top, & \bar{w}_2^\top, & [T_1,T_2]
\end{bmatrix}^\top
\end{aligned}
\end{equation*}
\normalsize
}

\begin{theorem}
Given a strongly connected digraph $\DG$, if the parameters $T_2 \geq T_1 > 0$, $\mu > 0$, $h \in \reals$, $\gamma > 0$, and positive definite symmetric matrices $P_1$, $P_2$, and $P_3$ are such that (\ref{cond:phil_1}) and (\ref{thrm_cond1}) hold, the hybrid system $\widetilde{\HS}_{m_{\sigma}}$ with input $\bar{m}_{\sigma}$ is ISS with respect to $\widetilde{\A}_{\varepsilon}$ given in (\ref{set:A_tilde}).
\end{theorem} 

The proof of this result largely follows the same approach used in the proof of Theorem \ref{thm:iss2}, namely, a Lyapunov analysis using the function candidate $V$ in (\ref{eqn:lyap}). Since the disturbance is present during flows, we show that the derivative of $V$ can be upper bounded resulting in a bounded disturbance in $V$ when evaluated along a given solution to $\widetilde{\HS}_{m_{\sigma}}$; see \cite{12} for more details.

\ifbool{two_col}{}{
\begin{proof}
Consider the same Lyapunov function candidate from the proof of Theorem \ref{thrm1} expressed for $\chi_{m_{\sigma}}$  $$V(\chi_{m_{\sigma}}) = V_1(\chi_{m_{\sigma}}) + V_2(\chi_{m_{\sigma}}) + V_{\varepsilon r} (\chi_{m_{\sigma}})$$ The contribution from the perturbation only affects the system during flows. For each $\chi_{m_{\sigma}} \in \tilde{C}_{m_{\sigma}}$ the change in $V$ is given by
\ifbool{rep}{
\begin{equation*}
\begin{aligned}
\langle \nabla V(\chi_{m_{\sigma}}), \tilde{f}_{m_{\sigma}}(\chi_{m_{\sigma}}) \rangle & \leq 2 \bar{z}_2^{\top} \exp{A^{\top}_{f_2} \tau } P_1 \exp{A_{f_2} \tau} (B_{f_2} \bar{w}_2 + \bar{m}_{\sigma_2} ) \\
& \hspace{5mm} + \bar{w}_1^{\top} (P_1 A_{f_3} + A_{f_3}^{\top} P_1) \bar{w}_1 + \bar{w}_2^{\top} (P_2 A_{f_4} + A_{f_4}^{\top} P_2) \bar{w}_2
\end{aligned}
\end{equation*}}{
\begin{equation*}
\begin{aligned}
\langle \nabla V(\chi_{m_{\sigma}}), \tilde{f}_{m_{\sigma}}(\chi_{m_{\sigma}}) \rangle & \leq 2 \bar{z}_2^{\top} \exp{A^{\top}_{f_2} \tau } P_1 \exp{A_{f_2} \tau} (B_{f_2} \bar{w}_2 + \bar{m}_{\sigma_2} ) \\
& \hspace{5mm} + \bar{w}_1^{\top} (P_1 A_{f_3} + A_{f_3}^{\top} P_1) \bar{w}_1 \\
& \hspace{10mm} + \bar{w}_2^{\top} (P_2 A_{f_4} + A_{f_4}^{\top} P_2) \bar{w}_2
\end{aligned}
\end{equation*}}
\noindent
From conditions (\ref{eqn:lyap1}) and (\ref{eqn:lyap2}), let $P_2 A_{f_3} + A_{f_3}^{\top} P_2 < -\beta_1 I$ and $P_3 A_{f_4} + A_{f_4}^{\top} P_3 < -\beta_2 I$ then one has
\ifbool{rep}{
\begin{equation}
\begin{aligned}
\langle \nabla V(\chi_{m_{\sigma}}), \tilde{f}_{m_{\sigma}}(\chi_{m_{\sigma}}) \rangle & \leq \kappa_1 | \bar{z}_2 | | \bar{w}_2 | + \frac{\kappa_1}{|B_{f_2}|} |\bar{z}_2| |\bar{m}_{\sigma_2}| - \beta_1 | \bar{w}_1 |^2 - \beta_2 | \bar{w}_2 |^2
\end{aligned}
\end{equation}}{
\begin{equation}
\begin{aligned}
\langle \nabla V(\chi_{m_{\sigma}}), \tilde{f}_{m_{\sigma}}(\chi_{m_{\sigma}}) \rangle & \leq \kappa_1 | \bar{z}_2 | | \bar{w}_2 | + \frac{\kappa_1}{|B_{f_2}|} |\bar{z}_2| |\bar{m}_{\sigma_2}| \\ 
& \hspace{5mm} - \beta_1 | \bar{w}_1 |^2 - \beta_2 | \bar{w}_2 |^2
\end{aligned}
\end{equation}}
\noindent
then applying Young's equality to the first and second terms one has
\ifbool{rep}{
\begin{equation}
\begin{aligned}
\langle \nabla V(\chi_{m_{\sigma}}), \tilde{f}_{m_{\sigma}}(\chi_{m_{\sigma}}) \rangle & \leq \frac{\kappa_1}{2 \epsilon} | \bar{z}_2 |^2 + \frac{\kappa_1 \epsilon}{2}  | \bar{w}_2 |^2  + \frac{\kappa_1}{2 \rho |B_{f_2}|} | \bar{z}_2 |^2 + \frac{\kappa_1 \rho}{2 |B_{f_2}|} | \bar{m}_{\sigma_2} |^2 - \beta_1 | \bar{w}_1|^2 - \beta_2 | \bar{w}_2|^2 \\
& \leq \Big ( \frac{\kappa_1}{2 \epsilon} + \frac{\kappa_1}{2 \rho |B_{f_2}|} \Big ) | \bar{z}_2 |^2 + \Big ( \frac{\kappa_1 \epsilon}{2} - \beta_2 \Big ) | \bar{w}_2|^2 - \beta_1 | \bar{w}_1|^2 + \frac{\kappa_1 \rho}{2 |B_{f_2}|} | \bar{m}_{\sigma_2} |^2
\end{aligned}
\end{equation}}{
\begin{equation*}
\begin{aligned}
\langle \nabla V(\chi_{m_{\sigma}}), \tilde{f}_{m_{\sigma}}(\chi_{m_{\sigma}}) \rangle & \leq \frac{\kappa_1}{2 \epsilon} | \bar{z}_2 |^2 + \frac{\kappa_1 \epsilon}{2}  | \bar{w}_2 |^2  + \frac{\kappa_1}{2 \rho |B_{f_2}|} | \bar{z}_2 |^2 \\
& \hspace{5mm} + \frac{\kappa_1 \rho}{2 |B_{f_2}|} | \bar{m}_{\sigma_2} |^2 - \beta_1 | \bar{w}_1|^2 - \beta_2 | \bar{w}_2|^2 \\
& \leq \Big ( \frac{\kappa_1}{2 \epsilon} + \frac{\kappa_1}{2 \rho |B_{f_2}|} \Big ) | \bar{z}_2 |^2 \\
& \hspace{5mm} + \Big ( \frac{\kappa_1 \epsilon}{2} - \beta_2 \Big ) | \bar{w}_2|^2 - \beta_1 | \bar{w}_1|^2 \\
& \hspace{10mm} + \frac{\kappa_1 \rho}{2 |B_{f_2}|} | \bar{m}_{\sigma_2} |^2
\end{aligned}
\end{equation*}}
\noindent
Since $|B_{f_2}| = 1$ then
\begin{equation}
\begin{aligned}
\langle \nabla V(\chi_{m_{\sigma}}), \tilde{f}_{m_{\sigma}}(\chi_{m_{\sigma}}) \rangle & \leq \tilde{\kappa} \big ( |\bar{z}_2|^2 + |\bar{w}_1|^2 + |\bar{w}_2|^2 \big ) \\ & \hspace{5mm} + \frac{\kappa_1 \rho}{2} | \bar{m}_{\sigma_2} |^2 \\
& \leq \tilde{\kappa} V(\chi_{\varepsilon}) + \frac{\kappa_1 \rho}{2} | \bar{m}_{\sigma_2} |^2
\end{aligned}
\end{equation}
\noindent
where $\tilde{\kappa} = \max \Big \{ \frac{\kappa_1}{2 \epsilon} + \frac{\kappa_1}{2 \rho}, \Big ( \frac{\kappa_1 \epsilon}{2} - \beta_2 \Big )  \Big \}$ and $\epsilon, \rho > 0$. Since the perturbation does not affect the system at jumps then, recall from the proof of Theorem \ref{thrm1} that, across jumps for each $\chi_{m_{\sigma}} \in \tilde{D}_{m_{\sigma}}$ and $g \in \tilde{G}_{m_{\sigma}}$ one has
\begin{equation*}
\begin{aligned}
V(g) - V(\chi_{m_{\sigma}}) \leq \minus \bar{\eta}_1^2 {+} \bar{z}_2^{\top} \Big ( A_{g_2}^{\top} \exp{A_{f_2}^{\top} v} P_1 \exp{A_{f_2} v} A_{g_2} {-} P_1 \Big ) \bar{z}_2 
\end{aligned}
\end{equation*}
\noindent
leading to the following bound
\begin{equation*}
\begin{aligned}
V(g) & \leq \Big (1 - \frac{\bar{\kappa}_2}{\alpha_2} \Big ) V(\chi_{\varepsilon}) + \bar{\kappa}_2  |\bar{w}|^2 \\
\end{aligned}
\end{equation*}
\noindent
from (\ref{eqn:v_jumps}). Then a general bound for the Lyapunov trajectory is given by
\ifbool{rep}{
\begin{equation}
\begin{aligned}
V(t,j) & \leq \exp \big ( \tilde{\kappa} \hspace{0.5mm} T_2 \big ) \Big ( \exp \big ( \tilde{\kappa} \hspace{0.5mm} T_2 \big ) \Big ( 1 - \frac{\bar{\kappa}_2}{\alpha_2} \Big ) \Big )^j V(0,0) + \bar{\kappa}_2 \exp \big ( \bar{\kappa} \hspace{0.5mm} T_2 \big ) \mathrm{sup}_{(t,j) \in \mbox{\footnotesize dom} \phi}|\bar{w}(t,j)|^2 \\ & \hspace{20mm} + \frac{\kappa_1 \rho}{2} \int_{0}^{t} \exp \big ( \tilde{\kappa} (t - \tau) \big ) | \bar{m}_{\sigma_2} |^2
\end{aligned}
\end{equation}}{
\begin{equation*}
\begin{aligned}
V(t,j) & \leq \exp \big ( \tilde{\kappa} \hspace{0.5mm} T_2 \big ) \Big ( \exp \big ( \tilde{\kappa} \hspace{0.5mm} T_2 \big ) \Big ( 1 - \frac{\bar{\kappa}_2}{\alpha_2} \Big ) \Big )^j V(0,0) \\
& \hspace{5mm} + \bar{\kappa}_2 \exp \big ( \bar{\kappa} \hspace{0.5mm} T_2 \big ) \mathrm{sup}_{(t,j) \in \mbox{\footnotesize dom} \phi}|\bar{w}(t,j)|^2 \\ & \hspace{10mm} + \frac{\kappa_1 \rho}{2} \int_{0}^{t} \exp \big ( \tilde{\kappa} (t - \tau) \big ) | \bar{m}_{\sigma_2} |^2
\end{aligned}
\end{equation*}}

\end{proof}}

\ifbool{two_col}{}{

\subsection{Noise on the communication and clock rate reference $\sigma^*$ with aperiodic communication events}

\ifbool{conf}{
\begin{example} In this example we demonstrate by simulation the system's robustness to noise on the communication channel and the clock rate reference $\sigma^*$. Consider the same system presented in Example \ref{ex:1}. Figure \ref{fig:3} shows input-to-state stability for the trajectories of the errors $e_i - e_k$ for the components $i \in \{1,2,3,4,5\}$ of a solution $\phi$ for the case where the system is subjected to communication noise $m_{e_i}(t,j) \in (0, 0.1)$ and noise on the clock rate reference $m_{\sigma^*_i}(t,j) \in (0.85, 1.15)$ for all $(t,j) \in \mbox{dom } \phi$, respectively. Moreover, after the respective transient period, the norm of the relative error $|e_i - e_k|$  converges to an average value of $0.0229$ when subjected to noise $m_{\sigma^*_i}$ and $0.0549$ for noise $m_{e_i}$.
\end{example}
}{
\begin{example} In this example we demonstrate the system $\HS$ robustness to noise on the communication channel and the clock rate reference $\sigma^*$. Consider the same system presented in Example \ref{ex:1}. 
Figure \ref{fig:3} shows ISS for the trajectories of the errors $e_i - e_k$ for the components $i \in \{1,2,3,4,5\}$ of a solution $\phi$ for the case where the system is subjected to communication noise $m_{e_i}(t,j) \in (0, 1)$ and noise on the clock rate reference $m_{\sigma^*_i}(t,j) \in (0.85, 1.15)$ for all $(t,j) \in  {\rm dom } \phi$, respectively. Moreover, after the respective transient period for each case, the norm of the relative error $|e_i - e_k|$  for each solution converges to an average value of $0.0229$ when subjected to noise $m_{\sigma^*_i}$ and $0.0549$ for noise $m_{e_i}$.
\end{example}}

\ifbool{rep}{
\begin{figure}[H]
\centering
\includegraphics[trim={0mm 0mm 0mm 0mm},clip,width=1 \textwidth]{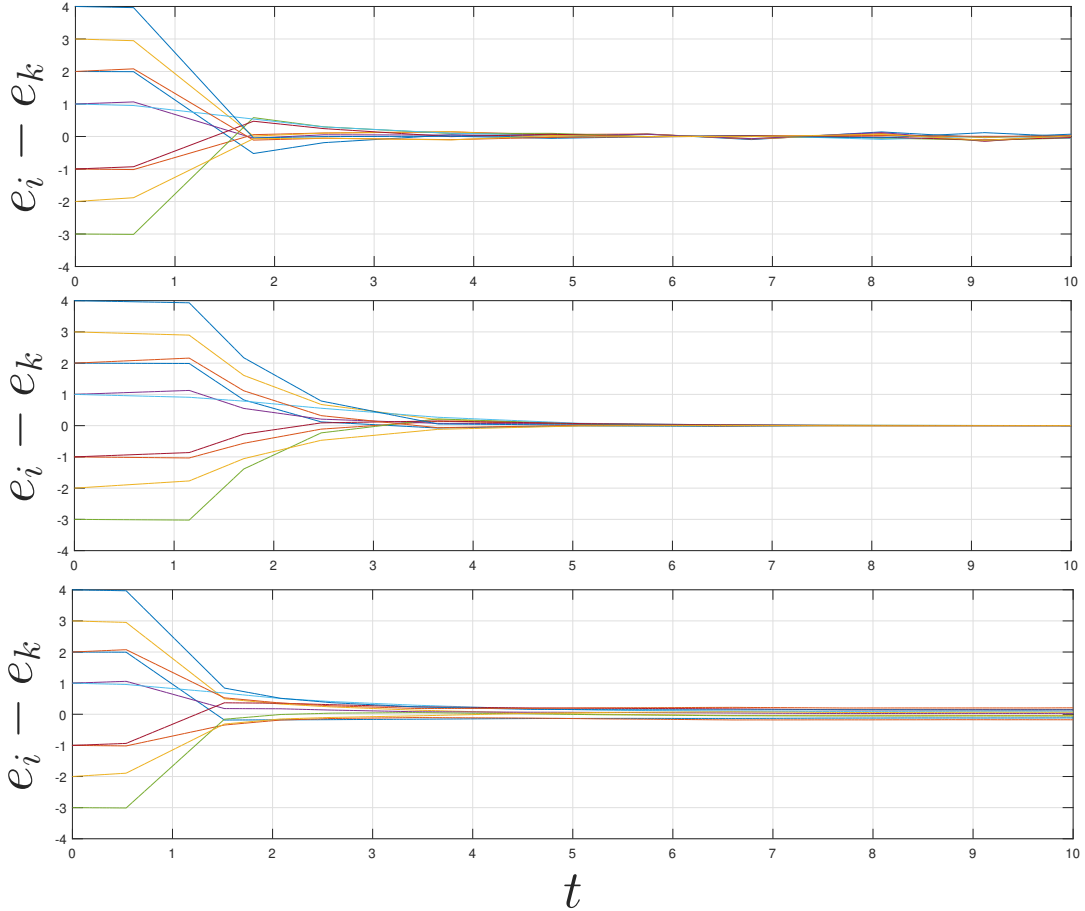}
\caption{\label{fig:3} (top) The trajectories of the errors $e_i - e_k$ for the components $i \in \{1,2,3,4,5\}$ of a solution $\phi$ for the case where the system is subjected to communication noise $m_{e_i}$ (top), noise on the internal clock output $m_{\tau_i}$ (middle), and noise on the clock rate reference $m_{\sigma^*_i}$ (bottom).}
\end{figure}
}{
\begin{figure}[H]
\centering
\includegraphics[trim={0mm 7mm 0mm 0mm},clip,width=0.5 \textwidth]{../Figures/error_all_noise_multi_v5.eps} \vspace{-5mm}
\caption{\label{fig:3} (top) The trajectories of the errors $e_i - e_k$ for the components $i \in \{1,2,3,4,5\}$ of a solution $\phi$ for the case where the system is subjected to communication noise $m_{e_i}$ (top) and noise on the clock rate reference $m_{\sigma^*_i}$ (bottom).} \vspace{-1mm}
\end{figure}}
}

\ifbool{rep}{\subsection{Robustness to Communication Event Delays}}{
\subsection{\color{subsectioncolor} Robustness to Communication Event Delays}}

\EndNew
In this section, we consider communication delays and show that the asymptotic stability of $\HS$ is robust (in a practical sense) to them. To model communication delays, we define a new hybrid system, denoted $\HS_{\delta}$, that is an augmentation of $\widetilde{\HS}_{\varepsilon}$ 
implementing the following mechanism: when a communication event is triggered due to $\tau$ reaching zero, 
the information to communicate is stored in a memory state and only available to the agents
after $T^\delta$ seconds after. \ifbool{rep}{}{Due to space constraints, the details behind the construction of $\HS_{\delta}$ are in \cite{12}. 
}

\ifbool{rep}{
Inspired by the results in \cite{213}, we analyze communication delays by defining a  higher-order model 
that implements a mechanism that memorizes the information at the original event and updates the
state of the system after the delay. To construct such higher-order model, the state vector $x$ of $\HS$ in (\ref{eqn:cl_hysys}) is decomposed into $x_1$ and $x_2$, where $x_1$ is composed of the states that are independent of a delay and $x_2$ represents the state components that are affected by the delay.

To this end, let $x^* = (x_1, x_2) \in \mathcal{X}^* := \mathcal{X}_1 \times \mathcal{X}_2$ where $x_1 = (e, \tau^*, \hat{a}, \hat{\tau}, \tau) \in \mathcal{X}_1$, where $\mathcal{X}_1 := \reals^n \times \reals^n_{\geq 0} \times \reals^n \times \reals^n_{\geq 0}$, and $x_2 = (u, \eta) \in \mathcal{X}_2$, where $\mathcal{X}_2 := \reals^n \times \reals^n$. We rewrite the hybrid model $\HS$ in (\ref{eqn:cl_hysys}) for the new state vector $x^*$ by defining the hybrid model
\begin{equation} \label{eqn:hy_recast}
\HS^* := (C^*,f^*,D^*,G^*) 
\end{equation} 
which is equivalently to $\HS$ by with slightly different notation. 
\noindent
In fact, the data $(C^*,f^*,D^*,G^*)$ is given by $(\dot{e}, \dot{\tau}^*, \dot{\hat{a}}, \dot{\hat{\tau}}, \dot{\tau}, \dot{u}, \dot{\eta} )$ $= (a + u - \sigma^* \textbf{1}_n, a, -\mu ( \hat{\tau} - \tau^* ), \hat{a} - (\hat{\tau} - \tau^*), -1, h \eta -\mu ( \hat{\tau} - \tau^* ), h \eta) =: f^*(x^*)$ for each $x^* \in C^*$ and
$(e^+, {\tau^*}^+, \hat{a}^+, \hat{\tau}^+, \tau^+, u^+, \eta^+)$ $= (e, \tau^*, \hat{a}, \hat{\tau}, [T_1,T_2], -\gamma \mathcal{L} e - \hat{a} + \sigma^* \textbf{1}_n, -\gamma \mathcal{L} e) =: G^*(x^*)$ for each $x^* \in D^*$, where $C^* := \mathcal{X}^*$ and $D^* := \{x^* \in \mathcal{X}^* : \tau = 0 \}$.
In this setting, the set of interest is
\begin{equation}
\begin{aligned}
\A^* := \{x^* \in \mathcal{X}^* : \hspace{1mm} & e_i = e_k, \eta_i = 0, \hat{a}_i = a_i, \hat{\tau}_i = \tau^*_i, \\
& u_i = \eta_i - \hat{a}_i + \sigma^* \hspace{9mm} \forall i ,k \in \mathcal{V} \}
\end{aligned}
\end{equation}

Given that $\HS^*$ is a rewritten version of $\HS$, the following two results are immediate. 

\begin{lemma}
The hybrid system $\HS^*$ satisfies the hybrid basic conditions defined in \cite[Assumption 6.5]{4}.
\end{lemma}

\begin{lemma}
The set $\A^*$ is asymptotically stable for the hybrid system $\HS^*$.
\end{lemma}

With the formulation of $\HS^*$, we are now ready to define the augmented model that allows for the analysis of delays on communication events.
Given the partition of  $x^*$ into non-delayed and delayed components $x_1$ and $x_2$, respectively, let $\widehat{G} : \mathcal{X}^* \rightrightarrows \mathcal{X}_1 \times \mathcal{X}_2 \times \mathcal{X}_2$ be a set-valued mapping such that 
$$\widehat{G}(x^*) = P(G^*(x^*) \times \{ x_2 \} )$$ 
for every $x^* = (x_1, x_2) \in \mathcal{X}^*$ 
where
$$P(y_1,y_2,x_2) = (y_1,x_2,y_2) \hspace{1cm} \forall (y_1, y_2, x_2) \in \mathcal{X}_1 \times \mathcal{X}_2 \times \mathcal{X}_2$$

Given a class-$\mathcal{K}_{\infty}$ function $\alpha$ and a continuous function $\rho : \reals^n \to \reals_{\geq 0}$ for every $T_{\delta} \geq 0$ let
\begin{align*}
C_T & := \{ x^* \in \mathcal{X}^*  : (x^* + \alpha(T_{\delta}) \rho(x^*) \mathbb{B}) \cap C^* \} \\
f_T(x^*) & := f^*((x^* + \alpha(T_{\delta}) \rho(x^*) \mathbb{B}) \cap C^*) + \alpha(T_{\delta}) \rho(x^*) \mathbb{B}
\end{align*}
\noindent
We can now define the delayed hybrid system augmented from $\HS^*$, let 
\begin{equation} \label{eqn:hy_delta}
\HS_{\delta} := (C_{\delta}, f_{\delta}, D_{\delta}, G_{\delta})
\end{equation}
\noindent
with state $x_{\delta} = (x^*, \mu, \tau_{\delta}) \in \mathcal{X}_{T^\delta} := \mathcal{X}^* \times \mathcal{X}_2 \times \reals$ where $\mu \in \mathcal{X}_2$, and $\tau_{\delta} \in \reals$ is given by the following data:
\begin{equation}
C_{\delta} := (C_T \times \{ \textbf{0}_{2n} \} \times \{ -1 \} ) \cup (C_T	\times \reals^{2n} \times [0,T^{\delta}] )
\end{equation}
\begin{equation}
f_{\delta}(x_{\delta}) := f_T(x^*)  \times \{ \textbf{0}_{2n} \} \times \{ - \mbox{min} \{\tau_{\delta} + 1, 1 \} \} \hspace{2mm} \forall x_{T^\delta} \in C_{T^\delta}
\end{equation}
\begin{equation}
D_{\delta} := (D^* \times \{ \textbf{0}_{2n} \} \times \{ -1 \} ) \cup (\mathcal{X}_1 \times
\mathcal{X}_2	\times \reals^{2n} \times \{ 0 \} )
\end{equation}
\begin{equation}
G_{\delta} := \begin{cases}
\widehat{G}(x^*) \times [0,T^{\delta}] \hspace{2mm} \mbox{ if } x_{\delta} \in D^* \times \{ \textbf{0}_{2n} \} \times \{ -1 \} \\
(x_1, \mu, \textbf{0}_{2n}, -1) \hspace{2mm} \mbox{ if } \mathcal{X}_1 \times
\mathcal{X}_2	\times \reals^{2n} \times \{ 0 \}
\end{cases}
\end{equation}
\noindent
The augmented components $\mu$ and $\tau_{\delta}$ define the memory state of the post-jump value of the delayed state $x_2$ and the timer variable that regulates delays, respectively (see \cite{213} for more details). 

For the robustness analysis of the delayed system, we consider stability of $\HS_{\delta}$ to the set 
\begin{equation}
\begin{aligned}
\A_{\delta} := \A^* \times \{ \textbf{0}_{2n} \} \times \{ -1 \} \cup (P(\A \times \mbox{cl}(\A_2^*) ) \times \{0 \} )  \end{aligned}
\end{equation}
\noindent
where $\A_2^*$ is the projection of $\A^*$ onto $\mathcal{X}_2$.

In order to show that the analysis of the augmented system $\HS_{\delta}$ applies to the recast system $\HS^*$ and thus our original system $\HS$, the following results from \cite{213} establish the equivalence between the two systems.

\begin{proposition}Given a solution $\phi_{x^*} = (\phi_{x_1}, \phi_{x_2} )$ of $\HS^*$, let $\{t_j\}^J_{j =1}$ be the sequence of jump times of $\phi_{x^*}$, i.e.,
\begin{equation}
(t_j,j),(t_j,j-1) \in \mbox{dom } \phi_{x^*} \hspace{10mm} \forall j \in \{1,2, \ldots, J \} \cap \N, J := \sup \{ j : \exists t, (t,j) \in \mbox{dom } \phi_{x^*} \}
\end{equation}
Consider the function $\phi_{x_\delta} : \mbox{dom} \phi_{x_\delta} \to \mathcal{X}^* \times \mathcal{X}_2 \times \reals$, where
\begin{equation}
\mbox{dom } \phi_{x_\delta} := \{(t,j) : (t,j/2) \in \mbox{dom } \phi_{x^*} \} \cup \Big ( \bigcup^J_{j=1} \{(t_j, 2j-1) \} \Big )
\end{equation}
and for every $(t, j) \in \mbox{dom } \phi_{x_\delta}$ $\phi_{x_\delta}(t, j) := (\phi_{x^*}(t, j/2), \textbf{0}_{2n} , -1)$
if $j$ is even, and
\begin{equation}
\phi_{x_\delta} (t, j) := (\phi_{x_1}(t, (j + 1)/2), \phi_{x_2}(t, (j - 1)/2), \phi_{x_2}(t, (j + 1)/2), 0)
\end{equation}
if $j$ is odd. Then, $\phi_{x_\delta}$ is a solution to $\HS^*$. Moreover, $\phi_{x_\delta}$ is maximal if $\phi_{x^*}$ is maximal.
\end{proposition}

\begin{lemma}
Let $\phi_{x_{\delta}} = (\phi_{x^*} , \phi_\mu, \phi_{\tau_{\delta}} )$ be a maximal solution of $\HS_{T^{\delta}}$ satisfying $\phi_{\tau}(0,0) = -1$, where $\phi_{x^*} = (\phi_{x_1},\phi_{x_2})$. Then, for every $(t, j) \in \mbox{ dom} \phi_{x_{\delta}}$,
\begin{equation}
(\mu (t,j) , \tau(t,j)) = \begin{cases}
(\textbf{0}_{2n}, -1 ) \hspace{1cm} \mbox{if } j \mbox{ is even} \\
(\phi_{x_2}(t,j+1),0) \hspace{1cm} \mbox{if } j \mbox{ is odd} 
\end{cases}
\end{equation}
Moreover, the function $\phi_{x^*} : dom \phi_{x^*} \to \mathcal{X}^*$ , where $dom \phi_{x^*} := \{ (t,j) : (t,2j) \in \mbox{dom } \phi_{x_{\delta}} \}$
and $\phi_{x^*} (t, j) := \phi_{x_{\delta}}(t, 2j)$ for all $(t, j) \in \mbox{dom } \phi_{x^*}$, is a maximal solution of $\HS$.
\end{lemma}

\begin{lemma}
Given any $T \geq 0$, the augmented hybrid system $\HS_{\delta}$ satisfies the hybrid basic conditions, and the set $F((x + \alpha_F (T ) \rho (x) B) \cap C)$ is convex for all $x \in C_T$.
\end{lemma}

\begin{definition}
Given an open set $\mathcal{U} \subset \reals^n$, a continuous function $w : \mathcal{U} \to \reals_{\geq 0}$ is a proper indicator of a compact set $\A' \subset \mathcal{U}$ on $\mathcal{U}$ if the following hold:
\begin{itemize}
\item $w(x)=0$ if and only if $x \in A$. 
\item Given any sequence $\{x_i \}_{i=0}^{\infty} \in \mathcal{U}$, $\lim_{i \to \infty} w(x_i) = \infty$ if $\lim_{i \to \infty} |x_i| = \infty$ or $\lim_{i \to \infty} |x_i|_{\reals^n \setminus \mathcal{U}} = 0$.
\end{itemize}
\end{definition}

With the formulation of the recast system $\HS^*$ and the augmented delayed system $\HS_{\delta}$ we can apply the following results from \cite{213}. Note that in the following results $\HS_{0}$ refers to the case where $T^{\delta} = 0$ for the system $\HS_{\delta}$ (see \cite{213} for more details).

\begin{proposition} \label{prop:robust_del}
Under Lemmas 2.2 and 4.2, the set $\A^*$ is compact and the basin of pre-attraction $\mathcal{B}_{\A_{\delta}}^p$ of the set $\A_{\delta}$ is open. Moreover, for every proper indicator $w$ of $\A_{\delta}$ on $\mathcal{B}_{\A_{\delta}}^p$, there exists a class-$\mathcal{KL}$ function $\beta$ such that every solution $\phi_0$ of $\HS_{0}$ originating from $\mathcal{B}_{\A_{\delta}}^p$ satisfies
\begin{equation} \label{eqn:beta}
w(\phi_{x_{\delta}}(t,j)) \leq \beta( w(\phi_{x_{\delta}}(0,0)), t+j)  \hspace{1cm} \forall(t,j) \in \mbox{dom } \phi_{x_{\delta}}
\end{equation}
\end{proposition}

\begin{theorem}
Suppose that the conditions of Proposition \ref{prop:robust_del}  hold. Consider the basin of pre-attraction $\mathcal{B}_{\A_{\delta}}^p$  of $\A_{\delta}$, along with any proper indicator $w$ of $\A_{\delta}$ on $\mathcal{B}_{\A_{\delta}}^p$ and any class-$\mathcal{KL}$ function $\beta$ satisfying (\ref{eqn:beta}) for every solution $\phi_0$ of $\HS_{0}$  originating from $\mathcal{B}_{\A_{\delta}}^p$. Then, for every compact set $K \subset \mathcal{B}_{\A_{\delta}}^p$ and every $\epsilon > 0$, there exists $T^{\delta} > 0$ such that every solution $\phi_{\delta}$ of $\HS_{\delta}$ originating from $K$ satisfies
\begin{equation}
w(\phi_{x_{\delta}}(t,j)) \leq \beta( w(\phi_{x_{\delta}}(0,0)), t+j) + \varepsilon \hspace{1cm} \forall(t,j) \in \mbox{dom } \phi_{x_{\delta}}
\end{equation}
\end{theorem}


\begin{theorem}
Given a strongly connected digraph $\DG$, if the parameters $T_2 \geq T_1 > 0$, $\mu > 0$, $h \in \reals$, $\gamma > 0$, and the positive definite symmetric matrices $P_1$, $P_2$, and $P_3$ are such that 
the conditions in Theorem~\ref{thrm1} hold, 
then the set $\A$ is 
semiglobally practically asymptotically stability  with respect to $T^{\delta}$ for $\HS$ in (\ref{eqn:cl_hysys}), namely, for each compact set $K \subset \mathcal{X}$ and each $\epsilon > 0$, 
there exists a $\cal KL$-function $\beta$ and $T^\delta > 0$ such that  each solution $\phi^\delta$ to $\HS$ with $\phi^\delta(0,0) \in K$ and with communication delay no larger than $T^\delta$ satisfies $|\phi^\delta(t,j)|_{\A} \leq \beta(|\phi^\delta(0,0)|_\A,t+j)+\eps$
for all $(t,j) \in \dom \phi^\delta$.
\end{theorem}

\begin{proof}
Pick a solution $\phi \in \mathcal{S}_{\HS}$ with $\phi(0,0) \in C \cup D$. Observe that for any solution $\phi$ to $\HS$ there exists an equivalent solution $\phi^*$ to $\HS^*$ moreover, via notions Proposition 6.8, Lemma 6.9, Lemma 6.10 there exists an equivalent solution $\phi_{x_\delta}$ to $\HS_{\delta}$. Thus, by the result of Theorem 6.13, it follows that the system $\HS$ is semiglobally practically asymptotically stable to $\A$ with respect to $T^{\delta}$.
\end{proof}
\EndNew
}{
\color{black}
{
\begin{theorem}\label{thm:RobustnessToDelays}
Given a strongly connected digraph $\DG$, if the parameters $T_2 \geq T_1 > 0$, $\mu > 0$, $h \in \reals$, $\gamma > 0$, and the positive definite symmetric matrices $P_1$, $P_2$, and $P_3$ are such that 
the conditions in Theorem~\ref{thrm1} hold, 
then the set $\A$ is 
semiglobally practically asymptotically stability  with respect to $T^{\delta}$ for $\HS$ in (\ref{eqn:cl_hysys}), namely, for each compact set $K \subset \mathcal{X}$ and each $\epsilon > 0$, 
there exists a $\cal KL$-function $\beta$ and $T^\delta > 0$ such that  each solution $\phi^\delta$ to $\HS$ with $\phi^\delta(0,0) \in K$ and with communication delay no larger than $T^\delta$ satisfies $|\phi^\delta(t,j)|_{\A} \leq \beta(|\phi^\delta(0,0)|_\A,t+j)+\eps$
for all $(t,j) \in \dom \phi^\delta$.
\end{theorem}

The proof of this result  follows from \cite[Theorem 5.3]{213} using the notions given 
in \cite[Section VII-C]{213}. Complete details and a proof of Theorem~\ref{thm:RobustnessToDelays} can be found in \cite{12}, along with a numerical validation. To emphasize the robustness of $\HS$ to small delays, Figure~\ref{fig:delay_sim} depicts  trajectories to a decentralized delay scenario where the data from each node incurs a maximum delay of $T^{\delta} = 0.1$ seconds, modeled by a delay timer $\tau_{\delta_i}$ for each node.
}

\begin{figure}
\centering
\includegraphics[trim={30mm 0mm 0mm 0mm},clip,width=0.43 \textwidth]{../Figures/Multi_fig_delay_png_2.png}
\vspace{-3mm} \caption{\label{fig:delay_sim} The trajectories of the solution $\phi$ for state component errors $e_i {-} e_k$, $\varepsilon_{a_i}$, $\tau$, and $\tau_{\delta_i}$} \vspace{-5mm}
\end{figure}

}

\ifbool{rep}{
\section{Comparisons and Examples} \label{sec:num}
}{
\section{Comparisons} \label{sec:num}
}

\ifbool{two_col}{ 
In this section, we compare our algorithm to several 
consensus-based clock synchronization algorithms from the literature through a numerical example. 
We consider a four agent setting and simulate each algorithm presented in \cite{carli2008pi} (PI-Consensus), \cite{bolognani2015randomized} (RandSync), and \cite{7} (Average TimeSync) to our hybrid algorithm HyNTP as in (\ref{eqn:cl_hysys}).

Consider $N = 4$ agents with clock dynamics as in (\ref{eqn:clk_dynamics1}) and (\ref{eqn:clk_dynamics2}) over a strongly connected graph with the following adjacency matrix
\ifbool{two_col}{
$\DG_A = [0, 1 , 0 , 1],$
$[1 , 0 , 1 , 0],$
$[0 , 1 , 0 , 1],$
$[1 , 0 , 1 , 0]$
}{
\vspace{-3mm}
\begin{equation} \label{eqn:graph_2} \small
\DG_A = \begin{pmatrix}
0 & 1 & 0 & 1 \\
1 & 0 & 1 & 0 \\
0 & 1 & 0 & 1 \\
1 & 0 & 1 & 0 \\
\end{pmatrix}
\end{equation}}
\noindent
and aperiodic communication events such that successive communications events are lower and upper bounded by $T_1 = 0.1$ and $T_2 = 0.5$, respectively. The initial conditions for the clock rates $\bar{a}_i$ and adjustable clock values $\bar{\tau}_i$ for each $i \in \nodes$ has been randomly chosen within the intervals $(0.5,1.5)$ and $(0, 200)$, respectively. Moreover, consider the case where the system is subjected to a communication noise $m_{\tau_i}(t,j) \in (0,1)$ on the clock measurements. Figure \ref{fig:noise_aper_1a} and \ref{fig:noise_aper_1b} show the trajectories of $\bar{\tau}$ and $\bar{a}_i$, respectively, for agents $i \in \{1,2,3,4\}$ for the \textit{HyNTP} algorithm and each of the comparison algorithms under consideration.

For the \textit{HyNTP} algorithm, setting $\sigma^* = 1$, it can be found that the parameters $h= -2$, $\mu = 9$, $\gamma = 0.06$ and $\epsilon = 4.752$ with suitable matrices $P_1$, $P_2$, and $P_3$  satisfy conditions (\ref{cond:phil_1}) and (\ref{thrm_cond1}) in Theorem \ref{thrm1} with $\bar{\kappa}_1 = 2.02$, $\kappa_1 = 19.22$, $\bar{\kappa}_2 = 1$, and $\alpha_2 = 44.03$.

\begin{figure}
\centering
\subfloat[\label{fig:noise_aper_1a}]{\includegraphics[width=0.25 \textwidth, trim={10mm 11mm 10mm 11mm},clip]{../Figures/adj_clock_aper_measNoise_compare.eps}}
\subfloat[\label{fig:noise_aper_1b}]{\includegraphics[width=0.25 \textwidth, trim={10mm 11mm 10mm 11mm},clip]{../Figures/adj_clock_rate_aper_measNoise_compare.eps}}
\caption{\label{fig:noise_aper_1} The evolution of the trajectories of the adjustable clocks $\bar{\tau}_i$ (\ref{fig:noise_aper_1a}) and adjustable clock rates $\bar{a}_i$ (\ref{fig:noise_aper_1b}) for each clock synchronization algorithm. From top to bottom, \textit{HyNTP}, \textit{Average TimeSync}, \textit{PI-Consensus}, and \textit{RandSync}.}
\vspace{-5mm}
\end{figure}

}{
In this section we compare our algorithm to several contemporary consensus-based clock synchronization algorithms from the literature through a numerical example. In particular, we consider a four agent setting and simulate each algorithm presented in \cite{carli2008pi} (PI-Consensus), \cite{bolognani2015randomized} (RandSync), and \cite{7} (Average TimeSync) to our hybrid algorithm HyNTP as in (\ref{eqn:cl_hysys}). We have restricted our comparison to these algorithms due to their shared assumptions on the underlying communication graph being strongly connected. Our first example considers the nominal case of zero noise and a fixed communication event period. The next example also considers the nominal case but with aperiodic communication events. We then present an example where the systems are subjected to communication noise with aperiodic communication. Our final example considers the case of noise on the clock rate while also being subjected to aperiodic communication events. 

\subsection{Nominal case with fixed communication event period}

Consider $N = 4$ agents with clock dynamics as in (\ref{eqn:clk_dynamics1}) and (\ref{eqn:clk_dynamics2}) over a strongly connected graph with the following adjacency matrix
\begin{equation} \label{eqn:graph_2}
\DG_A = \begin{pmatrix}
0 & 1 & 0 & 1 \\
1 & 0 & 1 & 0 \\
0 & 1 & 0 & 1 \\
1 & 0 & 1 & 0 \\
\end{pmatrix}
\end{equation}
\noindent
and a dwell time between communication events $T = 0.15$. The initial conditions for the clock rates $a_i$ and clock values $\tau_i$ for each $i \in \nodes$ has been randomly chosen within the intervals $(0.5,1.5)$ and $(0, 200)$, respectively. 

For the \textit{HyNTP} algorithm, we let $T_1 = T_2 = T = 0.15$, and $\sigma^* = 1$, then it can be found that the parameters $h= -2$, $\mu = 3$, $\gamma = 0.06$ and $\epsilon = 1.607$ with suitable matrices $P_1$, $P_2$, and $P_3$  satisfy conditions (\ref{cond:phil_1}) and (\ref{thrm_cond1}) in Theorem \ref{thrm1} with $\bar{\kappa}_1 = 6.86$, $\kappa_1 = 22.98$, $\bar{\kappa}_2 = 1$, and $\alpha_2 = 16.93$. 

Figure \ref{fig:nom_fixed_1} shows the trajectories of $e_i - e_k$, $\varepsilon_{a_i}$ for components $i \in \{1,2,3,4,5\}$ of a solution $\phi$ for the case where $\sigma = \sigma^*$

\ifbool{rep}{
\begin{figure}[H]
\centering
\includegraphics[trim={0mm 0mm 0mm 0mm},clip,width=1 \textwidth]{../Figures/adj_clock_compare.eps}
\caption{\label{fig:nom_fixed_1} The evolution of the trajectories of the adjustable clocks $\bar{\tau}_i$ for each clock synchronization algorithm. From top to bottom, \textit{HyNTP}, \textit{Average TimeSync}, \textit{PI-Consensus}, and \textit{RandSync}.}
\end{figure}}{
\begin{figure}
\centering
\includegraphics[trim={0mm 0mm 0mm 0mm},clip,width=0.5 \textwidth]{../Figures/adj_clock_compare.eps}
\caption{\label{fig:nom_fixed_1} The evolution of the trajectories of the adjustable clocks $\bar{\tau}_i$ for each clock synchronization algorithm. From top to bottom, \textit{HyNTP}, \textit{Average TimeSync}, \textit{PI-Consensus}, and \textit{RandSync}.}
\end{figure}}

\ifbool{rep}{
\begin{figure}[H]
\centering
\includegraphics[trim={0mm 0mm 0mm 0mm},clip,width=1 \textwidth]{../Figures/adj_clock_rate_compare.eps}
\caption{\label{fig:nom_fixed_2} The evolution of the trajectories of the adjustable clock rates $\bar{a}_i$ for each clock synchronization algorithm. From top to bottom, \textit{HyNTP}, \textit{Average TimeSync}, \textit{PI-Consensus}, and \textit{RandSync}.}
\end{figure}}{
\begin{figure}
\centering
\includegraphics[trim={0mm 0mm 0mm 0mm},clip,width=0.5 \textwidth]{../Figures/adj_clock_rate_compare.eps}
\caption{\label{fig:nom_fixed_2} The evolution of the trajectories of the adjustable clock rates $\bar{a}_i$ for each clock synchronization algorithm. From top to bottom, \textit{HyNTP}, \textit{Average TimeSync}, \textit{PI-Consensus}, and \textit{RandSync}.}
\end{figure}}

\subsection{Nominal case with aperiodic communication events}

Consider the same $N = 4$ agents with clock dynamics as in (\ref{eqn:clk_dynamics1}) and (\ref{eqn:clk_dynamics2}) over a strongly connected graph with the following adjacency matrix
\ifbool{two_col}{as in (\ref{eqn:graph_2})}{
\begin{equation*}
\DG_A = \begin{pmatrix}
0 & 1 & 0 & 1 \\
1 & 0 & 1 & 0 \\
0 & 1 & 0 & 1 \\
1 & 0 & 1 & 0 \\
\end{pmatrix}
\end{equation*}
\noindent}
and aperiodic communication events such that successive communications events are lower and upper bounded by $T_1 = 0.1$ and $T_2 = 0.5$, respectively. The initial conditions for the clock rates $a_i$ and clock values $\tau_i$ for each $i \in \nodes$ has been randomly chosen within the intervals $(0.5,1.5)$ and $(0, 200)$, respectively. 

For the \textit{HyNTP} algorithm, setting $\sigma^* = 1$, it can be found that the parameters $h= -2$, $\mu = 9$, $\gamma = 0.06$ and $\epsilon = 4.752$ with suitable matrices $P_1$, $P_2$, and $P_3$  satisfy conditions (\ref{cond:phil_1}) and (\ref{thrm_cond1}) in Theorem \ref{thrm1} with $\bar{\kappa}_1 = 2.02$, $\kappa_1 = 19.22$, $\bar{\kappa}_2 = 1$, and $\alpha_2 = 44.03$. 

Figure \ref{fig:1} shows the trajectories of $e_i - e_k$, $\varepsilon_{a_i}$ for components $i \in \{1,2,3,4,5\}$ of a solution $\phi$ for the case where $\sigma = \sigma^*$.

\ifbool{rep}{
\begin{figure}[H]
\centering
\includegraphics[trim={0mm 0mm 0mm 0mm},clip,width=1 \textwidth]{../Figures/adj_clock_aper_compare.eps}
\caption{\label{fig:nom_aper_1} The evolution of the trajectories of the adjustable clocks $\bar{\tau}_i$ for each clock synchronization algorithm. From top to bottom, \textit{HyNTP}, \textit{Average TimeSync}, \textit{PI-Consensus}, and \textit{RandSync}.}
\end{figure}}{
\begin{figure}
\centering
\includegraphics[trim={0mm 0mm 0mm 0mm},clip,width=0.5 \textwidth]{../Figures/adj_clock_aper_compare.eps}
\caption{\label{fig:nom_aper_1} The evolution of the trajectories of the adjustable clocks $\bar{\tau}_i$ for each clock synchronization algorithm. From top to bottom, \textit{HyNTP}, \textit{Average TimeSync}, \textit{PI-Consensus}, and \textit{RandSync}.}
\end{figure}}

\ifbool{rep}{
\begin{figure}[H]
\centering
\includegraphics[trim={0mm 0mm 0mm 0mm},clip,width=1 \textwidth]{../Figures/adj_clock_rate_aper_compare.eps}
\caption{\label{fig:nom_aper_2} The evolution of the trajectories of the adjustable clock rates $\bar{a}_i$ for each clock synchronization algorithm. From top to bottom, \textit{HyNTP}, \textit{Average TimeSync}, \textit{PI-Consensus}, and \textit{RandSync}.}
\end{figure}}{
\begin{figure}
\centering
\includegraphics[trim={0mm 0mm 0mm 0mm},clip,width=0.5 \textwidth]{../Figures/adj_clock_rate_aper_compare.eps}
\caption{\label{fig:nom_aper_2} The evolution of the trajectories of the adjustable clock rates $\bar{a}_i$ for each clock synchronization algorithm. From top to bottom, \textit{HyNTP}, \textit{Average TimeSync}, \textit{PI-Consensus}, and \textit{RandSync}.}
\end{figure}}

\subsection{Communication noise with aperiodic communication events}

Consider the same $N = 4$ agents with clock dynamics as in (\ref{eqn:clk_dynamics1}) and (\ref{eqn:clk_dynamics2}) over a strongly connected graph with the adjacency matrix given in (\ref{eqn:graph_2}) and aperiodic communication events such that successive communications events are lower and upper bounded by $T_1 = 0.1$ and $T_2 = 0.5$, respectively. The initial conditions for the clock rates $a_i$ and clock values $\tau_i$ for each $i \in \nodes$ has been randomly chosen within the intervals $(0.5,1.5)$ and $(0, 200)$, respectively. Moreover, consider the case where the system is subjected to a communication noise $m_{\tau_i}(t,j) \in (0,1)$ on the clock measurements.

For the \textit{HyNTP} algorithm, setting $\sigma^* = 1$, it can be found that the parameters $h= -2$, $\mu = 9$, $\gamma = 0.06$ and $\epsilon = 4.752$ with suitable matrices $P_1$, $P_2$, and $P_3$  satisfy conditions (\ref{cond:phil_1}) and (\ref{thrm_cond1}) in Theorem \ref{thrm1} with $\bar{\kappa}_1 = 2.02$, $\kappa_1 = 19.22$, $\bar{\kappa}_2 = 1$, and $\alpha_2 = 44.03$. 

\ifbool{rep}{
\begin{figure}[H]
\centering
\includegraphics[trim={0mm 0mm 0mm 0mm},clip,width=1 \textwidth]{../Figures/adj_clock_aper_measNoise_compare.eps}
\caption{\label{fig:nom_aper_1} The evolution of the trajectories of the adjustable clocks $\bar{\tau}_i$ for each clock synchronization algorithm. From top to bottom, \textit{HyNTP}, \textit{Average TimeSync}, \textit{PI-Consensus}, and \textit{RandSync}.}
\end{figure}}{
\begin{figure}
\centering
\includegraphics[trim={0mm 0mm 0mm 0mm},clip,width=0.5 \textwidth]{../Figures/adj_clock_aper_measNoise_compare.eps}
\caption{\label{fig:nom_aper_1} The evolution of the trajectories of the adjustable clocks $\bar{\tau}_i$ for each clock synchronization algorithm. From top to bottom, \textit{HyNTP}, \textit{Average TimeSync}, \textit{PI-Consensus}, and \textit{RandSync}.}
\end{figure}}

\ifbool{rep}{
\begin{figure}[H]
\centering
\includegraphics[trim={0mm 0mm 0mm 0mm},clip,width=1 \textwidth]{../Figures/adj_clock_rate_aper_measNoise_compare.eps}
\caption{\label{fig:nom_aper_2} The evolution of the trajectories of the adjustable clock rates $\bar{a}_i$ for each clock synchronization algorithm. From top to bottom, \textit{HyNTP}, \textit{Average TimeSync}, \textit{PI-Consensus}, and \textit{RandSync}.}
\end{figure}}{
\begin{figure}
\centering
\includegraphics[trim={0mm 0mm 0mm 0mm},clip,width=0.5 \textwidth]{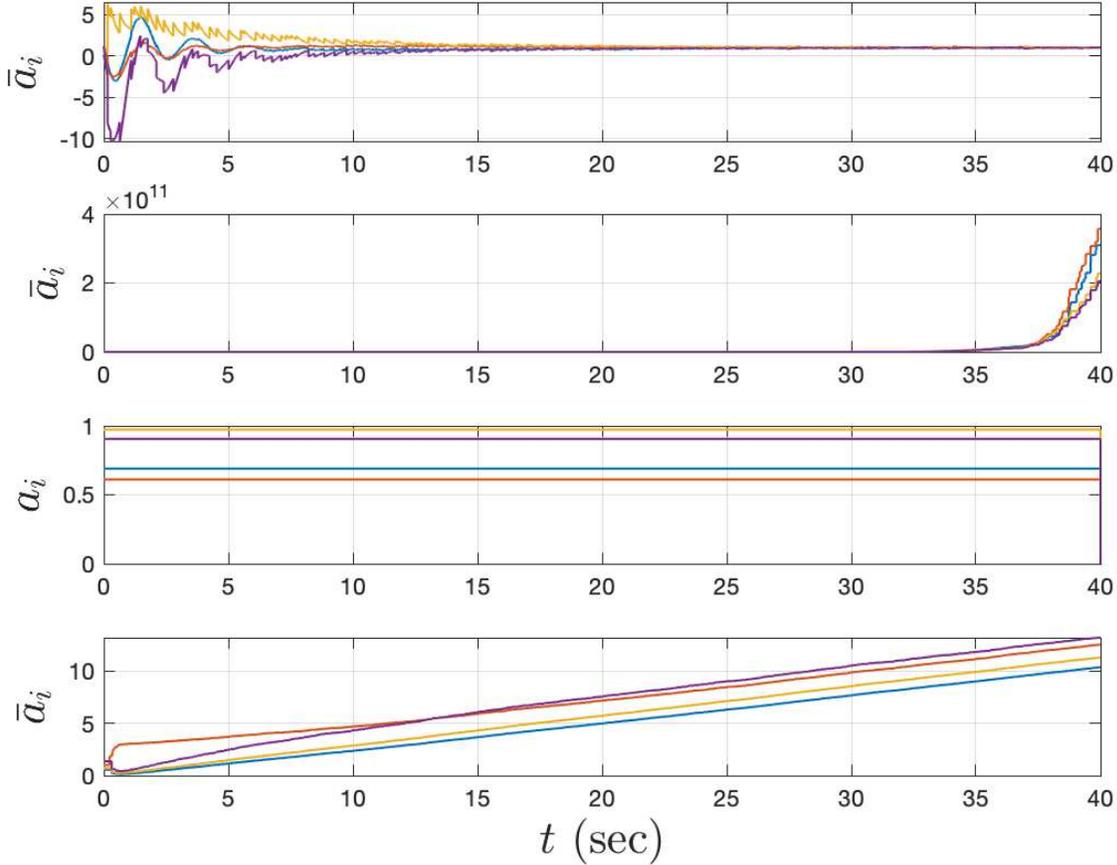}
\caption{\label{fig:nom_aper_2} The evolution of the trajectories of the adjustable clock rates $\bar{a}_i$ for each clock synchronization algorithm. From top to bottom, \textit{HyNTP}, \textit{Average TimeSync}, \textit{PI-Consensus}, and \textit{RandSync}.}
\end{figure}}
}

\ifbool{rep}{
\subsection{Nominal case with N=100 agents}

Consider N=100 agents with clock dynamics as in (\ref{eqn:clk_dynamics1}) and (\ref{eqn:clk_dynamics2}) over a strongly connected graph with adjacency matrix $\DG = \textbf{1}_{N \times N} - I_{N \times N}$ with aperiodic communication events such that successive communications events are lower and upper bounded by $T_1 = 0.001$ and $T_2 = 0.01$, respectively. The initial conditions for the clock rates $a_i$ and clock values $\tau_i$ for each $i \in \nodes$ has been randomly chosen within the intervals $(0.85,1.15)$ and $(-100, 100)$, respectively.

Given $T_1$, $T_2$, and $\sigma^*=1$ then it can be found that the parameters $h_i = -1.3$ for each $i \in \nodes$, $\mu = 3$, and $\gamma = 0.005$, sutiable matrices $P_1$, $P_2$, $P_3$, and $\epsilon = 2.180$ satisfy conditions (\ref{cond:phil_1}) and  (\ref{thrm_cond1}) in Theorem \ref{thrm1} with $\bar{\kappa}_1 = 6.896$, $\kappa_1 = 30.067$, $\bar{\kappa}_2 = 1$, and $\alpha_2 = 55.598$. Figure \ref{fig:nom_100_nodes} shows the evolution of the trajectories of the virtual clocks $\tilde{\tau}_i$ for $i \in \{1,2, \ldots, 100\}$.\footnote{Code at github.com/HybridSystemsLab/HybridClockSync}

\begin{figure}[H]
\centering
\includegraphics[trim={0mm 0mm 0mm 0mm},clip,width=1 \textwidth]{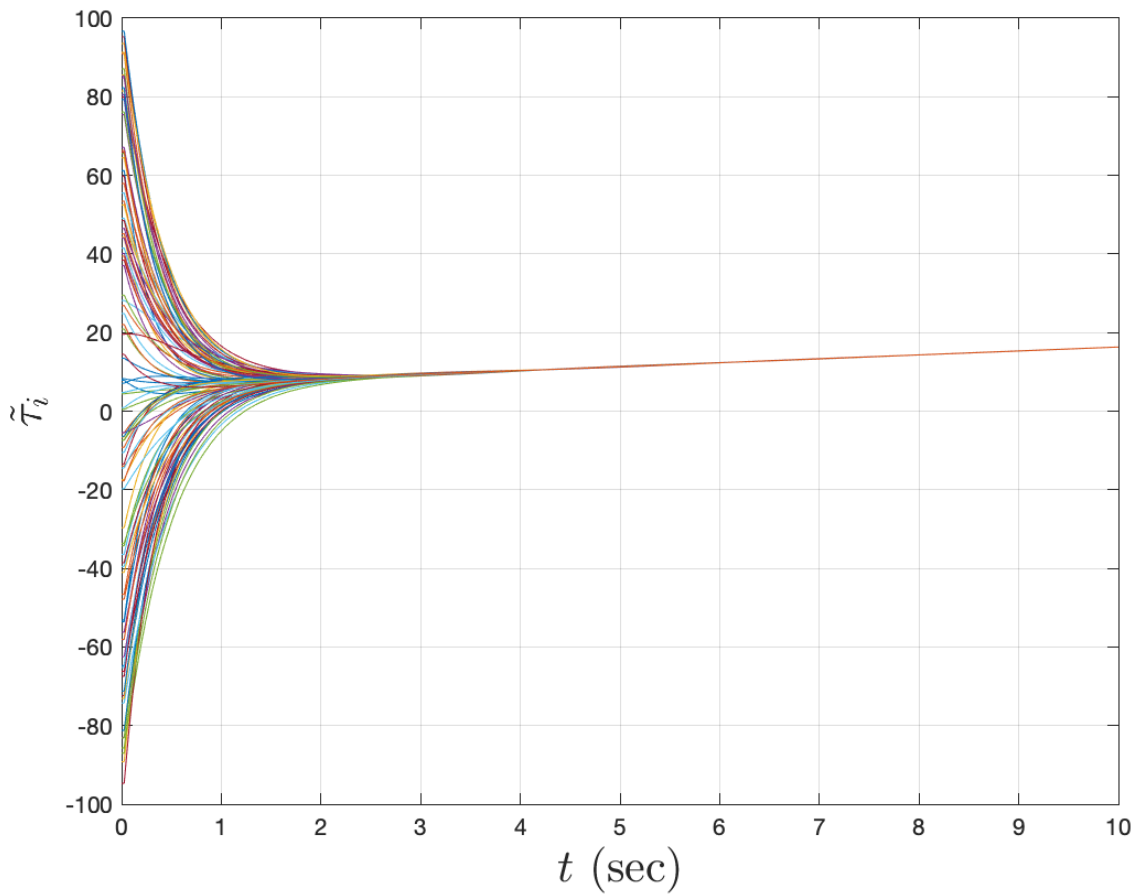}
\caption{\label{fig:nom_100_nodes} The evolution of the trajectories of the virtual clocks $\tilde{\tau}_i$ for $i \in \{1,2, \ldots, 100\}$.}
\end{figure}

\EndNew
}{}

\section{Conclusion}

In this paper, we modeled a network of clocks with aperiodic communication that utilizes a distributed hybrid controller to achieve synchronization, using the hybrid systems framework. Results were given to guarantee and show synchronization of the timers, exponentially fast. Numerical results validating the exponentially fast convergence of the timers were also given. Numerical results were also provided to demonstrate performance against a similar class of clock synchronization algorithms.


\ifbool{two_col}{
\appendix
}{
\appendix
\section{Appendix}
}

\ifbool{two_col}{}{
\subsection{Proof of Lemma \ref{lem:hbc}}
By inspection of the hybrid system data defining $\HS$ given in (\ref{eqn:cl_hysys}) and below it, the following is observed:
\begin{itemize}
\item The set $C$ is a closed subset of $\reals^m$ since, $C = \mathcal{X}$ and $\mathcal{X}$ is the Cartesian product of closed sets. Similar arguments show that $D$ is closed since it can be written as $$D = \reals^n \times \reals^n \times \reals^n \times \reals_{\geq 0}^n \times \reals^n \times \reals_{\geq 0}^n \times \{0\}$$ Thus, (A1) holds.
\item $f: \mathcal{X} \rightarrow \mathcal{X}$ is linear affine in the state and thus continuous on $C$. Moreover, since $\mbox{dom } f = \mathcal{X} = C$, $C \subset \mbox{dom } f$ holds. Thus, (A2) holds.
\item To show that the set-valued map $G$ defined in (\ref{eqn:cl_hysys}) satisfies (A3), note that the graph of $G$ is given by \begin{align*}
\mbox{gph}(G) & = \{ (x,y) : x \in D,  y \in G(x) \} \\
& = D \times \big ( \reals^n \times \reals^n \times \reals^n \times \reals_{\geq 0}^n \times \reals^n \times \reals_{\geq 0}^n \times [T_1,T_2] \big )
\end{align*}
\noindent
is closed. Thus, via \cite[Lemma 5.10]{4}, $G$ is outer semicontinuous and locally bounded at each $x \in D$. Moreover, by definition, we have that $\mbox{dom } G = D$. Hence, (A3) holds.
\end{itemize}
\hfill$\blacksquare$
}

\ifbool{two_col}{}{
\subsection{Proof of Lemma \ref{lem:complete1}}
For each $\xi \in C$, the tangent cone $T_C(\xi)$, as defined in \cite[Definition 5.12]{4}, is given by
\begin{equation*}
\begin{aligned}
T_C(\xi) {=} \begin{cases} \reals^n {\times} \reals^n {\times} \reals^n {\times} \reals^n_{\geq 0} {\times} \reals^n {\times} \reals^n_{\geq 0} {\times} \reals_{\geq 0} \hspace{2mm} \mbox{if } \xi \in \mathcal{X}^1 \\ 
\reals^n {\times} \reals^n {\times} \reals^n {\times} \reals^n_{\geq 0} {\times} \reals^n {\times} \reals^n_{\geq 0} {\times} \reals \hspace{5.5mm} \mbox{if } \xi \in \mathcal{X}^2 \\
\reals^n {\times} \reals^n {\times} \reals^n {\times} \reals^n_{\geq 0} {\times} \reals^n {\times} \reals^n_{\geq 0} {\times} \reals_{\leq 0} \hspace{2mm} \mbox{if } \xi \in \mathcal{X}^3
\end{cases}
\end{aligned}
\end{equation*}
\normalsize
\noindent
where $\mathcal{X}^1 := \{ x \in \mathcal{X} : \tau = 0 \}$, $\mathcal{X}^2 := \{ x \in \mathcal{X} : \tau \in (0, T_2) \}$, and $\mathcal{X}^3 := \{ x \in \mathcal{X} : \tau = T_2 \}$. By inspection, from the definition of $f$ in (\ref{eqn:cl_hysys}, $f(x) \cap T_C(x) \neq \emptyset$ holds  for every $x \in C \setminus D$. Then, since $\HS$ satisfies the hybrid basic conditions, as  shown in Lemma \ref{lem:hbc}, by \cite[Proposition 6.10]{4} there exists a nontrivial solution $\phi$ to $\HS$ with $\phi(0,0) = \xi$. Moreover, every $\phi \in \mathcal{S}_{\HS}$ satisfies one of the following conditions:
\begin{enumerate}
\setlength{\itemindent}{7mm}
\item[a)] $\phi$ is complete;
\item[b)] $\mbox{dom } \phi$ is bounded and the interval $I^J$, where $J = \mbox{sup}_{j} \mbox{dom } \phi$, has nonempty interior and $t \mapsto \phi(t,J)$ is a maximal solution to $\dot{x} \in F(x)$, in fact $\lim_{t \to T} |\phi(t,J)| = \infty $, where $T = \mbox{sup}_{t} \mbox{dom } \phi$;
\item[c)] $\phi(T,J) \notin C \cup D$, where $(T,J) = \mbox{sup dom } \phi$.
\end{enumerate}
\noindent
Now, since $G(D) \subset C \cup D = \mathcal{X}$ due to the definition of $G$, case c) does not occur. Additionally, one can eliminate case b) since $f$ is globally Lipschitz continuous on $C$ due to being linear affine in the state. Hence, only a) holds.
\hfill$\blacksquare$
}

\ifbool{two_col}{
\indent \textit{Proof of Lemma \ref{lem:complete}}:
Pick an initial condition $\xi \in \A$. Let $\phi$ be a maximal solution to $\HS$ with $\phi(0,0) = \xi$.\footnote{Note that for a given solution $\phi(t,j)$ to $\HS$, the solution components are given by $\phi(t,j) =$ $\big (\phi_e(t,j), \phi_u(t,j), \phi_\eta(t,j),$ $\phi_{\tau^*}(t,j), \phi_{\hat{a}}(t,j), \phi_{\hat{\tau}}(t,j),  \phi_{\tau}(t,j) \big )$}
\begin{itemize}
\item Consider the case where $\phi(0,0) \in \A \setminus D$. The initial conditions of the components of $\phi$ satisfy $\phi_{e_i}(0,0) = \phi_{\eta_i}(0,0) = 0$ for the clock errors $e_i$, $\phi_{\hat{\tau}_i}(0,0) = \phi_{\tau_i^*}(0,0)$ for the estimated clocks $\hat{\tau}_i$, $\phi_{\hat{a}_i}(0,0) = \phi_{a_i}(0,0)$ for the clock rates $\hat{a}_i$ and $\phi_{u_i}(0,0) = \phi_{\eta_i}(0,0) -\phi_{\hat{a}_i}(0,0) + \sigma^*$ for the control input for each $i \in \nodes$. With $f$ being linear affine, the constrained differential equation $\dot{x} = f(x)$ $x \in C$ has unique solutions. Let $[0, t_1] \times \{0\}  \subset \mbox{dom } \phi$ with $t_1 > 0$, which exists since $\phi(0,0) \in \A \setminus D$. From the definition of $f$, the solution components of the states $u$, $\eta$, and $e$ during this interval remain constant. From the  definition of $f$ in (\ref{eqn:cl_hysys}) we have that the components of the solution $\phi$ satisfy  $\phi_{e_i}(t,j) = \phi_{e_k}(t,j) $,  $\phi_{\eta}(t,j) = 0$, $\phi_{\hat{a}_i}(t,j) = \phi_{a_i}(t,j)$, $\phi_{\hat{\tau}_i}(t,j) = \phi_{\tau_i^*}(t,j)$, and $\phi_{u_i}(t,j) = \phi_{\eta_i}(t,j) -\phi_{\hat{a}_i}(t,j) + \sigma^*$ for each $(t,j) \in [0, t_1] \times \{0\}$.  Therefore, the solution $\phi$ does not leave the set $\A$ when $\phi(0,0) \in \A \setminus D$.
\item Consider the case where $\phi(0,0) \in \A \cap D$. Since flow is not possible from $\phi(0,0)$ as $\phi_{\tau}(0,0) = 0$, $(\{0\} \times \{0\}) \cup (\{0\} \times \{1\}) \subset \mbox{dom } \phi$ as the solution $\phi$ jumps initially. By inspection, the jump map $G$ in (\ref{eqn:cl_hysys}) only affects the states $\eta$, $u$, and $\tau$, with the other state components unchanged. Since the quantity $- \gamma \mathcal{L} e$ in the jump map is zero at $\phi(0,0)$, we have that $\phi_{\eta}(0,1) = - \gamma \mathcal{L} \phi_{e}(0,0) = 0$. Moreover, since $\hat{a}$ is constant across jumps, $\phi_{\hat{a}}(0,1) = \phi_{\hat{a}}(0,0)$, then, $\phi_{u}(0,1) = - \gamma \mathcal{L} \phi_{e}(0,0) - \phi_{\hat{a}}(0,0) + \sigma^* \textbf{1}_n$ $\phi_{\eta}(0,1) - \phi_{\hat{a}}(0,1) + \sigma^* \textbf{1}_n$
The timer $\tau$ resets to a point in the interval $[T_1,T_2]$, namely, $\phi_{\tau}(0,1) \in [T_1,T_2]$. Then, the full solution $\phi$ at $(0,1)$ satisfies
$\phi(0,1) \in \A$.
\end{itemize}
Since these properties hold for each $\xi \in \A$, $\A$ is forward invariant for $\HS$.
\hfill$\blacksquare$
}{
\subsection{Proof of Lemma \ref{lem:complete}}
Pick an initial condition $\xi \in \A$. Let $\phi$ be a maximal solution to $\HS$ with $\phi(0,0) = \xi$.\footnote{ Note that for a given solution $\phi(t,j)$ to $\HS$, the solution components are given by $\phi(t,j) = \big (\phi_e(t,j), \phi_u(t,j), \phi_\eta(t,j), \phi_{\tau^*}(t,j), \phi_{\hat{a}}(t,j), \phi_{\hat{\tau}}(t,j),  \phi_{\tau}(t,j) \big )$}
\begin{itemize}
\item Consider the case where $\phi(0,0) \in \A \setminus D$. The initial conditions of the components of $\phi$ satisfy $\phi_{e_i}(0,0) = \phi_{\eta_i}(0,0) = 0$ for the clock errors $e_i$, $\phi_{\hat{\tau}_i}(0,0) = \phi_{\tau_i^*}(0,0)$ for the estimated clocks $\hat{\tau}_i$, $\phi_{\hat{a}_i}(0,0) = \phi_{a_i}(0,0)$ for the clock rates $\hat{a}_i$ and $\phi_{u_i}(0,0) = \phi_{\eta_i}(0,0) -\phi_{\hat{a}_i}(0,0) + \sigma^*$ for the control input for each $i \in \nodes$. With $f$ being linear affine and, thus, globally Lipschitz continuous on $C$, the constrained differential equation $\dot{x} = f(x)$ $x \in C$ has unique solutions. Let $[0, t_1] \times \{0\}  \subset \mbox{dom } \phi$ with $t_1 > 0$, which exists since $\phi(0,0) \in \A \setminus D$. Observe that, from the definition of $f$, the solution components of the states $u$, $\eta$, and $e$ during this interval remain constant. This is evident since $\dot{\phi}_u = h \phi_{\eta}(0,0) - \mu \big ( \phi_{\hat{\tau}}(0,0) - \phi_{\tau^*}(0,0) \big ) = 0$ with $\phi_{\eta}(0,0) = 0$, $\dot{\phi}_{\eta} = h \phi_{\eta}(0,0) = 0$, and $\phi_{\hat{\tau}}(0,0) = \phi_{\tau^*}(0,0)$; hence, $\dot{\phi}_e = \phi_a(0,0) + \phi_u(0,0) - \sigma^* \textbf{1}_n = 0$. From the  definition of $f$ in (\ref{eqn:cl_hysys}) we have that the components of the solution $\phi$ satisfy  $\phi_{e_i}(t,j) = \phi_{e_k}(t,j) $,  $\phi_{\eta}(t,j) = 0$, $\phi_{\hat{a}_i}(t,j) = \phi_{a_i}(t,j)$, $\phi_{\hat{\tau}_i}(t,j) = \phi_{\tau_i^*}(t,j)$, and $\phi_{u_i}(t,j) = \phi_{\eta_i}(t,j) -\phi_{\hat{a}_i}(t,j) + \sigma^*$ for each $(t,j) \in [0, t_1] \times \{0\}$.  Therefore, the solution $\phi$ does not leave the set $\A$ during the interval $[0, t_1] \times \{0\}$ when $\phi(0,0) \in \A \setminus D$.
\item Consider the case where $\phi(0,0) \in \A \cap D$. Since flow is not possible from $\phi(0,0)$ as $\phi_{\tau}(0,0) = 0$, $(\{0\} \times \{0\}) \cup (\{0\} \times \{1\}) \subset \mbox{dom } \phi$ as the solution $\phi$ jumps initially. By inspection, the jump map $G$ in (\ref{eqn:cl_hysys}) only affects the states $\eta$, $u$, and $\tau$, whereas the value of the other state components remains unchanged. Since the quantity $- \gamma \mathcal{L} e$ in the jump map is zero at $\phi(0,0)$, we have that $\phi_{\eta}(0,1) = - \gamma \mathcal{L} \phi_{e}(0,0) = 0$. Moreover, since $\hat{a}$ is constant across jumps, $\phi_{\hat{a}}(0,1) = \phi_{\hat{a}}(0,0)$, then,
\begin{align*}
\phi_{u}(0,1) & = - \gamma \mathcal{L} \phi_{e}(0,0) - \phi_{\hat{a}}(0,0) + \sigma^* \textbf{1}_n \\
& = \phi_{\eta}(0,1) - \phi_{\hat{a}}(0,1) + \sigma^* \textbf{1}_n
\end{align*} 
Lastly, we have that the timer $\tau$ resets to a point in the interval $[T_1,T_2]$, namely, $\phi_{\tau}(0,1) \in [T_1,T_2]$. Then, the full solution $\phi$ at $(0,1)$ satisfies
\begin{align*}
\phi(0,1) \in \begin{bmatrix}
\phi_e(0,1) \\
\phi_{\eta}(0,1) - \phi_{\hat{a}}(0,1) + \sigma^* \textbf{1}_n \\
\phi_{\eta}(0,1) \\
\phi_{\tau^*}(0,1) \\
\phi_{\hat{a}}(0,1) \\
\phi_{\hat{\tau}}(0,1) \\
[T_1,T_2]
\end{bmatrix}
\end{align*}
Hence, from the definition of $\A$, $\phi(0,1) \in \A$.
\end{itemize}
Since this property holds for each $\xi \in \A$, we have that solutions from $\A$ cannot flow out of $\A$ and cannot jump out of $\A$ since $G(\A \cap D) \subset \A$. Hence, $\A$ is forward invariant for $\HS$.
\hfill$\blacksquare$
}


\ifbool{two_col}{
\indent \textit{Proof of Lemma \ref{lem:set_dists}}:
For each $x \in \mathcal{X}$, the distance from $x$ to the set $\A$ is given as $|x|_{\A} = \inf_{y \in \A} |x - y|$. Evaluating the distance directly, one has 
$|x|_{\A} =  \inf_{ e^* \in E}  \mbox{\rm sqrt} \big (  (e {-} e^*)^{\top} (e {-} e^*)$ $+ (u {-} \eta + \hat{a} {-} \sigma^*\textbf{1}_n)^{\top}(u {-} \eta + \hat{a} {-} \sigma^*\textbf{1}_n)$ $ + \eta^{\top}\eta + (\hat{a} - a)^{\top}(\hat{a} - a) + (\hat{\tau} - \tau^*)^{\top}(\hat{\tau} - \tau^*) \big )$ where $E := \{ e^* \in \reals^n : e_i^* = e_k^* \hspace{2mm} \forall i,k \in \mathcal{V} \}$. When $u = \eta - \hat{a} + \sigma^* \textbf{1}_n$ we have $|x|_{\A} =  \inf_{ e^* \in E}  \mbox{\rm sqrt} \big ( (e - e^*)^{\top} (e - e^*)  + \eta^{\top} \eta$ $+ (\hat{a} - a)^{\top}(\hat{a} - a) + (\hat{\tau} - \tau^*)^{\top}(\hat{\tau} - \tau^*) \big )$. For each $x_{\varepsilon} \in \mathcal{X}_{\varepsilon}$, the distance from $x_{\varepsilon}$ to the set $\A_{\varepsilon}$ is given as $|x_{\varepsilon}|_{\A_{\varepsilon}} = \inf_{y \in \A_{\varepsilon}} |x_{\varepsilon} - y|$. Evaluating the distance directly, one has $|x_{\varepsilon}|_{\A_{\varepsilon}} =  \inf_{ e^* \in E}  \mbox{\rm sqrt} \big ( (e - e^*)^{\top} (e - e^*)$ $+ \eta^{\top} \eta + \varepsilon_{a}^{\top} \varepsilon_{a} + \varepsilon_{\tau}^{\top}\varepsilon_{\tau} \big )$. Making the appropriate substitutions for $\varepsilon_{\tau}$ and $\varepsilon_a$, we get $|x_{\varepsilon}|_{\A_{\varepsilon}} =  \inf_{ e^* \in E}  \mbox{\rm sqrt} \big ( (e {-} e^*)^{\top} (e {-} e^*) $ $+ \eta^{\top} \eta + (\hat{a} {-} a)^{\top}(\hat{a} {-} a) + (\hat{\tau} - \tau^*)^{\top}(\hat{\tau} - \tau^*) \big )$. Now, for each $(x_{\varepsilon}, \hat{\tau}, \tau^*) \in \mathcal{X}$, the distance from the point $\widetilde{M}(x_{\varepsilon}, \hat{\tau}, \tau^*)$ to the set $\A$ is given by $|\widetilde{M}(x_{\varepsilon}, \hat{\tau}, \tau^* )|_{\A} = \inf_{y \in \A} |\widetilde{M}(x_{\varepsilon}, \hat{\tau}, \tau^*) - y|$. Computing this distance, one has $|\widetilde{M}(x_{\varepsilon}, \hat{\tau}, \tau^* )|_{\A} = \inf_{ e^* \in E , \alpha_{\tau^*} \in \reals^n_{\geq 0},\alpha_{\tau} \in [0,T_2]} | (e, \eta - (a - \varepsilon_a) + \sigma^* \textbf{1}_n, \eta,$ $\hat{\tau} - \varepsilon_{\tau}, a - \varepsilon_a, \varepsilon_{\tau} + \tau^*, \tau)$ $- ( e^* , \eta -\hat{a} + \sigma^*\textbf{1}_n, 0, \alpha_{\tau^*},a,\tau^*,\alpha_{\tau} )|$. Making the appropriate substitutions for $\varepsilon_{\tau}$ and $\varepsilon_a$, we get $|\widetilde{M}(x_{\varepsilon}, \hat{\tau}, \tau^*)|_{\A} =  \inf_{ e^* \in E}  \mbox{\rm sqrt} \big ( (e - e^*)^{\top} (e - e^*)$ $ + \eta^{\top}\eta + (\hat{a} - a)^{\top}(\hat{a} - a)$ $+ (\hat{\tau} - \tau^*)^{\top}(\hat{\tau} - \tau^*) \big )$ Thus, we have that 
$|\widetilde{M}(x_{\varepsilon}, \hat{\tau}, \tau^*)|_{\A} = |x|_{\A} = |x_{\varepsilon}|_{\A_{\varepsilon}}$.
\hfill$\blacksquare$
}{
\subsection{Proof of Lemma \ref{lem:set_dists}}
For each $x \in \mathcal{X}$, the distance from $x$ to the set $\A$ is given as
\begin{equation}
|x|_{\A} = \inf_{y \in \A} |x - y|
\end{equation}
\noindent
Evaluating the distance directly, one has
\begin{align*}
|x|_{\A} & = \inf_{y \in \A} |x - y| \\
& = \inf_{ e^* \in E , \alpha_{\tau^*} \in \reals^n_{\geq 0},\alpha_{\tau} \in [0,T_2]}  | (e,u,\eta, \tau^*, \hat{a}, \hat{\tau}, \tau) \\
& \hspace{25mm} - ( e^* ,\eta -\hat{a} + \sigma^*\textbf{1}_n, 0, \alpha_{\tau^*},a,\tau^*,\alpha_{\tau} )| \\
& = \inf_{ e^* \in E ,\alpha_{\tau^*} \in \reals^n_{\geq 0},\alpha_{\tau} \in [0,T_2]} | ( e - e^*  ,u  - \eta + \hat{a} - \sigma^*\textbf{1}_n, \eta, \\
& \hspace{25mm} \tau^* - \alpha_{\tau^*}, \hat{a} - a, \hat{\tau} - \tau^*, \tau - \alpha_{\tau}) | \\
& =  \inf_{ e^* \in E}  | ( e {-} e^* ,u {-} \eta + \hat{a} {-} \sigma^*\textbf{1}_n, \eta, 0, \hat{a} {-} a, \hat{\tau} {-} \tau^*, 0) | \\
& =  \inf_{ e^* \in E}  \mbox{\rm sqrt} \Big (  (e {-} e^*)^{\top} (e {-} e^*) \\
& \hspace{10mm}  + (u {-} \eta + \hat{a} {-} \sigma^*\textbf{1}_n)^{\top}(u {-} \eta + \hat{a} {-} \sigma^*\textbf{1}_n) \\
& \hspace{10mm} + \eta^{\top}\eta + (\hat{a} - a)^{\top}(\hat{a} - a) + (\hat{\tau} - \tau^*)^{\top}(\hat{\tau} - \tau^*) \Big )
\end{align*}
\noindent
 where $E := \{ e^* \in \reals^n : e_i^* = e_k^* \hspace{2mm} \forall i,k \in \mathcal{V} \}$.  When $u = \eta - \hat{a} + \sigma^* \textbf{1}_n$ we have
\begin{align*}
|x|_{\A} & =  \inf_{ e^* \in E}  \mbox{\rm sqrt} \Big ( (e - e^*)^{\top} (e - e^*)  + \eta^{\top} \eta \\ 
& \hspace{20mm} + (\hat{a} - a)^{\top}(\hat{a} - a) + (\hat{\tau} - \tau^*)^{\top}(\hat{\tau} - \tau^*) \Big )
\end{align*}
\noindent
For each $x_{\varepsilon} \in \mathcal{X}_{\varepsilon}$, the distance from $x_{\varepsilon}$ to the set $\A_{\varepsilon}$ is given as
\begin{equation}
|x_{\varepsilon}|_{\A_{\varepsilon}} = \inf_{y \in \A_{\varepsilon}} |x_{\varepsilon} - y|
\end{equation}
\noindent
Evaluating the distance directly, one has
\begin{align*}
|x_{\varepsilon}|_{\A_{\varepsilon}} & = \inf_{y \in \A_{\varepsilon}} |x_{\varepsilon} - y| \\
& = \inf_{ e^* \in E , \alpha_{\tau^*} \in \reals^n_{\geq 0}, \alpha_{\tau} \in [0,T_2]} | (e, \eta, \varepsilon_{a}, \varepsilon_{\tau}, \tau) \\
& \hspace{50mm} - ( e^* , 0, 0, 0,\alpha_{\tau} )| \\
& = \inf_{ e^* \in E , \alpha_{\tau^*} \in \reals^n_{\geq 0}, \alpha_{\tau} \in [0,T_2]} | ( e - e^* , \eta, \varepsilon_{a}, \varepsilon_{\tau}, \tau - \alpha_{\tau}) | \\
& =  \inf_{ e^* \in E}  | ( e - e^* , \eta, \varepsilon_{a}, \varepsilon_{\tau}, 0) | \\
& =  \inf_{ e^* \in E}  \sqrt{ (e - e^*)^{\top} (e - e^*)  + \eta^{\top} \eta + \varepsilon_{a}^{\top} \varepsilon_{a} + \varepsilon_{\tau}^{\top}\varepsilon_{\tau}}
\end{align*}
\noindent
Making the appropriate substitutions for $\varepsilon_{\tau}$ and $\varepsilon_a$, we get
\begin{align*}
|x_{\varepsilon}|_{\A_{\varepsilon}} & =  \inf_{ e^* \in E}  \mbox{\rm sqrt} \Big ( (e {-} e^*)^{\top} (e {-} e^*)  + \eta^{\top} \eta + (\hat{a} {-} a)^{\top}(\hat{a} {-} a) \\
& \hspace{45mm} + (\hat{\tau} - \tau^*)^{\top}(\hat{\tau} - \tau^*) \Big )
\end{align*}
\noindent
Now, for each $(x_{\varepsilon}, \hat{\tau}, \tau^*) \in \mathcal{X}$, the distance from the point $\widetilde{M}(x_{\varepsilon}, \hat{\tau}, \tau^*)$ to the set $\A$ is given by
\begin{equation}
|\widetilde{M}(x_{\varepsilon}, \hat{\tau}, \tau^* )|_{\A} = \inf_{y \in \A} |\widetilde{M}(x_{\varepsilon}, \hat{\tau}, \tau^*) - y|
\end{equation}
\noindent
Computing this distance, one has
\begin{align*} 
|\widetilde{M}(x_{\varepsilon}, \hat{\tau}, \tau^* )|_{\A} & = \inf_{y \in \A} |\widetilde{M}(x_{\varepsilon}, \hat{\tau}, \tau^*) - y| \\
& \hspace{-20mm} = \inf_{ e^* \in E , \alpha_{\tau^*} \in \reals^n_{\geq 0},\alpha_{\tau} \in [0,T_2]} | (e, \eta - (a - \varepsilon_a) + \sigma^* \textbf{1}_n, \eta, \\
& \hspace{20mm} \hat{\tau} - \varepsilon_{\tau}, a - \varepsilon_a, \varepsilon_{\tau} + \tau^*, \tau) \\ 
& \hspace{5mm} - ( e^* , \eta -\hat{a} + \sigma^*\textbf{1}_n, 0, \alpha_{\tau^*},a,\tau^*,\alpha_{\tau} )|
\end{align*}
Making the appropriate substitutions for $\varepsilon_{\tau}$ and $\varepsilon_a$, we get
\begin{align*} 
|\widetilde{M}(x_{\varepsilon}, \hat{\tau}, \tau^*)|_{\A} & \\ 
& \hspace{-23mm} = \inf_{ e^* \in E , \alpha_{\tau^*} \in \reals^n_{\geq 0},\alpha_{\tau} \in [0,T_2]} | (e, \eta - \hat{a} + \sigma^* \textbf{1}_n, \eta, \tau^*, \hat{a}, \hat{\tau}, \tau) \\
& \hspace{7mm} - ( e^* ,\eta -\hat{a} + \sigma^*\textbf{1}_n, 0, \alpha_{\tau^*},a,\tau^*,\alpha_{\tau} )| \\
& \hspace{-23mm} = \inf_{ e^* \in E , \alpha_{\tau^*} \in \reals^n_{\geq 0},\alpha_{\tau} \in [0,T_2]} | ( e {-} e^* ,\eta {-} \hat{a} + \sigma^* \textbf{1}_n {-} \eta + \hat{a} {-} \sigma^*\textbf{1}_n, \\
& \hspace{7mm} \eta - 0, \tau^* - \alpha_{\tau^*}, \hat{a} - a, \hat{\tau} - \tau^*, \tau - \alpha_{\tau}) | \\
& \hspace{-23mm} =  \inf_{ e^* \in E}  | ( e - e^* , 0, \eta, 0, \hat{a} - a, \hat{\tau} - \tau^*, 0) | \\
& \hspace{-23mm} =  \inf_{ e^* \in E}  \mbox{\rm sqrt} \Big ( (e - e^*)^{\top} (e - e^*)  + \eta^{\top}\eta + (\hat{a} - a)^{\top}(\hat{a} - a) \\
& \hspace{30mm} + (\hat{\tau} - \tau^*)^{\top}(\hat{\tau} - \tau^*) \Big )
\end{align*} 
Thus, we have that 
$$|\widetilde{M}(x_{\varepsilon}, \hat{\tau}, \tau^*)|_{\A} = |x|_{\A} = |x_{\varepsilon}|_{\A_{\varepsilon}}$$
\hfill$\blacksquare$
}


\ifbool{two_col}{
\indent \textit{Proof of Lemma \ref{lem:GES1}}:
Suppose the set $\A_{\varepsilon}$ is GES for $\HS_{\varepsilon}$. By Definition \ref{def:stability} there exist $\kappa, \alpha > 0$ such that each maximal solution $\phi^{\varepsilon}$ to $\HS_{\varepsilon}$ satisfies $|\phi^{\varepsilon}(t,j)|_{\A_{\varepsilon}} \leq \kappa \exp (\minus \alpha (t+ j)) |\phi^{\varepsilon}(0,0)|_{\A_{\varepsilon}}$ for each $(t,j) \in \mbox{dom } \phi^{\varepsilon}$. Now, pick any maximal solution $\phi$ to $\HS$. Through an application of Lemma \ref{lem:equiv}, there exists a corresponding solution $\phi^{\varepsilon}$ to $\HS_{\varepsilon}$ such that $\phi (t,j) = \widetilde{M} \big (\phi^{\varepsilon}(t,j), \phi_{\hat{\tau}}(t,j), \phi_{\tau^*}(t,j) \big )$ for each $(t,j) \in \mbox{\rm dom } \phi$. 
Given that $\phi^{\varepsilon}$ satisfies the given exponential bound, using the relationship between distances in Lemma \ref{lem:set_dists} we have that $\phi$ satisfies $|\phi(t,j)|_{\A} \leq \kappa \exp (\minus \alpha (t+ j)) |\phi(0,0)|_{\A}$. Then, the set $\A$ is GES for $\HS$. 
\hfill$\blacksquare$
}{
\subsection{Proof of Lemma \ref{lem:GES1}}
Suppose the set $\A_{\varepsilon}$ is GES for $\HS_{\varepsilon}$. By Definition \ref{def:stability} there exist $\kappa, \alpha > 0$ such that each maximal solution $\phi^{\varepsilon}$ to $\HS_{\varepsilon}$ satisfies
\begin{equation} \label{eqn:ges_bnd}
|\phi^{\varepsilon}(t,j)|_{\A_{\varepsilon}} \leq \kappa \exp (\minus \alpha (t+ j)) |\phi^{\varepsilon}(0,0)|_{\A_{\varepsilon}}
\end{equation} 
for each $(t,j) \in \mbox{dom } \phi^{\varepsilon}$. Now, pick any maximal solution $\phi$ to $\HS$. Through an application of Lemma \ref{lem:equiv}, there exists a corresponding solution $\phi^{\varepsilon}$ to $\HS_{\varepsilon}$ such that $$\phi (t,j) = \widetilde{M} \big (\phi^{\varepsilon}(t,j), \phi_{\hat{\tau}}(t,j), \phi_{\tau^*}(t,j) \big )$$ for each $(t,j) \in \mbox{\rm dom } \phi$. 
Given that $\phi^{\varepsilon}$ satisfies (\ref{eqn:ges_bnd}), using relationship (\ref{eqn:set_dists_1}) between distances in Lemma \ref{lem:set_dists} we have that $\phi$ satisfies
\begin{equation}
|\phi(t,j)|_{\A} \leq \kappa \exp (\minus \alpha (t+ j)) |\phi(0,0)|_{\A}
\end{equation}
Then, the set $\A$ is GES for $\HS$. 
\hfill$\blacksquare$
}

\ifbool{two_col}{
\indent \textit{Proof of Lemma \ref{lem:equiv_transform}}:
Pick a solution $\tilde{\phi} \in \mathcal{S}_{\tilde{\HS}_{\varepsilon}}$. with $\tilde{\phi} = (\tilde{\phi}_{\bar{z}_1}, \tilde{\phi}_{\bar{z}_2}, \tilde{\phi}_{\bar{w}_1}, \tilde{\phi}_{\bar{w}_2}, \tau)$, however, recall that $\bar{z}_1 := (\bar{e}_1, \bar{\eta}_1)$, $\bar{z}_2 := (\bar{e}_2,\ldots,\bar{e}_N, \bar{\eta}_2, \ldots, \bar{\eta}_N)$, $\bar{w}_1 = (\bar{\varepsilon}_{a_1},\bar{\varepsilon}_{\tau_1})$, and $\bar{w}_2 = (\bar{\varepsilon}_{a_2},\ldots, \bar{\varepsilon}_{a_n},\bar{\varepsilon}_{\tau_2},\ldots, \bar{\varepsilon}_{\tau_n})$ thus, through a reordering of the solution trajectories, one has that with some of the above notation, $\tilde{\phi}$ can be rewritten as $\tilde{\phi} = (\tilde{\phi}_{\bar{e}}, \tilde{\phi}_{\bar{\eta}}, \tilde{\phi}_{\bar{\varepsilon}_a}, \tilde{\phi}_{\bar{\varepsilon}_{\tau}}, \tau)$. Then, recall the change of coordinates $\bar{e} = \mathcal{T}^{\minus 1} e$, $\bar{\eta} = \mathcal{T}^{\minus 1} \eta$, $\bar{\varepsilon}_a = \mathcal{T}^{\minus 1} \varepsilon_a$, and $\bar{\varepsilon}_{\tau} = \mathcal{T}^{\minus 1} \varepsilon_{\tau}$. Since $\Tmat^{\minus 1}$ is an invertible time-invariant linear operator, applying its inverse $\Tmat$ to the components of $\tilde{\phi}$, one has $\big (\Tmat \tilde{\phi}_{\bar{e}}(t,j), \Tmat \tilde{\phi}_{\bar{\eta}}(t,j), \Tmat \tilde{\phi}_{\bar{\varepsilon}_a}(t,j), \Tmat \tilde{\phi}_{\bar{\varepsilon}_{\tau}}(t,j) \big) = 
\big ( \phi_e(t,j),
\phi_{\eta}(t,j),
\phi_{\varepsilon_a}(t,j),
\phi_{\varepsilon_{\tau}}(t,j) \big )$
for each $(t,j) \in \mbox{dom } \tilde{\phi}$. \ifbool{two_col}{Noting that $\tau$ is equivalent for $\HS_{\varepsilon}$ and $\widetilde{\HS}_{\varepsilon}$.}{Note that the dynamics of the variable $\tau$, responsible for governing the flows and the jumps of both $\HS_{\varepsilon}$ and $\widetilde{\HS}_{\varepsilon}$, is identical for the two systems. Thus, the set of solutions for the component $\tau$ is the same between the two system.} Therefore, it follows that$\tilde{\phi}(t,j) = \Gamma^{\minus 1} \phi(t,j)$ for each $(t,j) \in \mbox{dom } \tilde{\phi}$.

Conversely, pick a solution $\phi \in \mathcal{S}_{\HS_{\varepsilon}}$, let $\phi = (\phi_e, \phi_{\eta}, \phi_{\varepsilon_a}, \phi_{\varepsilon_{\tau}}, \tau)$ and recall the change of coordinates $\bar{e} = \mathcal{T}^{\minus 1} e$, $\bar{\eta} = \mathcal{T}^{\minus 1} \eta$, $\bar{\varepsilon}_a = \mathcal{T}^{\minus 1} \varepsilon_a$, and $\bar{\varepsilon}_{\tau} = \mathcal{T}^{\minus 1} \varepsilon_{\tau}$. Since $\Tmat^{\minus 1}$ is a time-invariant linear operator, applying it to the components of $\phi$, one has
$\big ( \Tmat^{\minus 1} \phi_e(t,j), \Tmat^{\minus 1} \phi_{\eta}(t,j), \Tmat^{\minus 1} \phi_{\varepsilon_a}(t,j), \Tmat^{\minus 1} \phi_{\varepsilon_{\tau}}(t,j) \big ) = \big ( \tilde{\phi}_{\bar{e}}(t,j), \tilde{\phi}_{\bar{\eta}}(t,j), \tilde{\phi}_{\bar{\varepsilon}_a}(t,j), \tilde{\phi}_{\bar{\varepsilon}_{\tau}}(t,j) \big )$ for each $(t,j) \in \mbox{dom } \phi$. Thus, it follows that $\phi(t,j) = \Gamma \tilde{\phi}(t,j)$ for each $(t,j) \in \mbox{dom } \phi$.
\hfill$\blacksquare$
}{
\subsection{Proof of Lemma \ref{lem:equiv_transform}} 
Pick a solution $\tilde{\phi} \in \mathcal{S}_{\tilde{\HS}_{\varepsilon}}$ with $\tilde{\phi} = (\tilde{\phi}_{\bar{z}_1}, \tilde{\phi}_{\bar{z}_2}, \tilde{\phi}_{\bar{w}_1}, \tilde{\phi}_{\bar{w}_2}, \tau)$, however, recall that $\bar{z}_1 := (\bar{e}_1, \bar{\eta}_1)$, $\bar{z}_2 := (\bar{e}_2,\ldots,\bar{e}_N, \bar{\eta}_2, \ldots, \bar{\eta}_N)$, $\bar{w}_1 = (\bar{\varepsilon}_{a_1},\bar{\varepsilon}_{\tau_1})$, and $\bar{w}_2 = (\bar{\varepsilon}_{a_2},\ldots, \bar{\varepsilon}_{a_n},\bar{\varepsilon}_{\tau_2},\ldots, \bar{\varepsilon}_{\tau_n})$. Thus, through a reordering of the solution trajectories, one has that with some of the above notation, $\tilde{\phi}$ can be rewritten as $\tilde{\phi} = (\tilde{\phi}_{\bar{e}}, \tilde{\phi}_{\bar{\eta}}, \tilde{\phi}_{\bar{\varepsilon}_a}, \tilde{\phi}_{\bar{\varepsilon}_{\tau}}, \tau)$. Then, recall the change of coordinates $\bar{e} = \mathcal{T}^{\minus 1} e$, $\bar{\eta} = \mathcal{T}^{\minus 1} \eta$, $\bar{\varepsilon}_a = \mathcal{T}^{\minus 1} \varepsilon_a$, and $\bar{\varepsilon}_{\tau} = \mathcal{T}^{\minus 1} \varepsilon_{\tau}$. Since $\Tmat^{\minus 1}$ is an invertible time-invariant linear operator, applying its inverse $\Tmat$ to the components of $\tilde{\phi}$, one has $\big (\Tmat \tilde{\phi}_{\bar{e}}(t,j), \Tmat \tilde{\phi}_{\bar{\eta}}(t,j), \Tmat \tilde{\phi}_{\bar{\varepsilon}_a}(t,j), \Tmat \tilde{\phi}_{\bar{\varepsilon}_{\tau}}(t,j) \big) = 
\big ( \phi_e(t,j),
\phi_{\eta}(t,j),
\phi_{\varepsilon_a}(t,j),
\phi_{\varepsilon_{\tau}}(t,j) \big )$
for each $(t,j) \in \mbox{dom } \tilde{\phi}$. Note that the dynamics of the variable $\tau$, responsible for governing the flows and the jumps of both $\HS_{\varepsilon}$ and $\widetilde{\HS}_{\varepsilon}$, is identical for the two systems. Thus, the set of solutions for the component $\tau$ is the same between the two system. Therefore, it follows that $\tilde{\phi}(t,j) = \Gamma^{\minus 1} \phi(t,j)$ for each $(t,j) \in \mbox{dom } \tilde{\phi}$.

Conversely, we can pick a solution $\phi \in \mathcal{S}_{\HS_{\varepsilon}}$, let $\phi$ $=$ $(\phi_e, \phi_{\eta}, \phi_{\varepsilon_a}, \phi_{\varepsilon_{\tau}}, \tau)$ and recall the change of coordinates $\bar{e} = \mathcal{T}^{\minus 1} e$, $\bar{\eta} = \mathcal{T}^{\minus 1} \eta$, $\bar{\varepsilon}_a = \mathcal{T}^{\minus 1} \varepsilon_a$, and $\bar{\varepsilon}_{\tau} = \mathcal{T}^{\minus 1} \varepsilon_{\tau}$. Since $\Tmat^{\minus 1}$ is a time-invariant linear operator, applying it to the components of $\phi$, one has
$\big ( \Tmat^{\minus 1} \phi_e(t,j), \Tmat^{\minus 1} \phi_{\eta}(t,j), \Tmat^{\minus 1} \phi_{\varepsilon_a}(t,j), \Tmat^{\minus 1} \phi_{\varepsilon_{\tau}}(t,j) \big ) = \big ( \tilde{\phi}_{\bar{e}}(t,j), \tilde{\phi}_{\bar{\eta}}(t,j), \tilde{\phi}_{\bar{\varepsilon}_a}(t,j), \tilde{\phi}_{\bar{\varepsilon}_{\tau}}(t,j) \big )$ for each $(t,j) \in \mbox{dom } \phi$. Thus, it follows that $\phi(t,j) = \Gamma \tilde{\phi}(t,j)$ for each $(t,j) \in \mbox{dom } \phi$.
\hfill$\blacksquare$}


\ifbool{two_col}{\indent \textit{Proof of Lemma \ref{lem:set_equiv}}:  This proof has been omitted as it is simply exploits the property of norms for linear systems. 
}{
\subsection{Proof of Lemma \ref{lem:set_equiv}} 
\sloppy
Pick a point $\tilde{z}' = (\bar{e}_1', \bar{\eta}_1', \bar{e}_2',\ldots,\bar{e}_N',$ $\bar{\eta}_2', \ldots, \bar{\eta}_N', \bar{\varepsilon}_{a_1}', \bar{\varepsilon}_{\tau_1}', \bar{\varepsilon}_{a_2}',\ldots, \bar{\varepsilon}_{a_n}',\bar{\varepsilon}_{\tau_2}',\ldots, \bar{\varepsilon}_{\tau_n}') \in \reals^{4N}$ such that $(\tilde{z}', \tau') \in \tilde{\A}_{\varepsilon}$  for some $\tau' \in [0,T_2]$ , i.e., $\tilde{z}' = ( e^*_1 , 0, \textbf{0}_{N-1}, \textbf{0}_{N-1}, 0, 0, \textbf{0}_{N-1}, \textbf{0}_{N-1})$ with $e^*_1 \in \reals$. Given that the digraph $\DG$ is strongly connected, there exists a nonsingular matrix $\mathcal{T}$ as in (\ref{eqn:T_mat})  that allows for the following coordinate change: $\bar{e} = \mathcal{T}^{\minus 1} e$, $\bar{\eta} = \mathcal{T}^{\minus 1} \eta$, $\bar{\varepsilon}_a = \mathcal{T}^{\minus 1} \varepsilon_a$, and $\bar{\varepsilon}_{\tau} = \mathcal{T}^{\minus 1} \varepsilon_{\tau}$. Now, by left multiplying $(\tilde{z}', \tau')$ by $\Gamma$ one has
\begin{equation} \label{eqn:coord_change} 
\begin{aligned}
e & {=} \mathcal{T} \begin{bmatrix} \bar{e}_1' & \bar{e}_2' & \ldots & \bar{e}_N'  \end{bmatrix}^{\top} {=} \begin{bmatrix} v_1 & \mathcal{T}_1  \end{bmatrix} \begin{bmatrix} e_1^* & \textbf{0}_{N \minus 1}^\top  \end{bmatrix}^{\top} {=} e_1^* \textbf{1}_N \\
\eta & = \mathcal{T} \begin{bmatrix} \bar{\eta}_1' & \bar{\eta}_2' & \ldots & \bar{\eta}_N'  \end{bmatrix}^{\top} = \begin{bmatrix} v_1 & \mathcal{T}_1  \end{bmatrix} \begin{bmatrix} 0 & \textbf{0}_{N \minus 1}^\top  \end{bmatrix}^{\top} = \textbf{0}_{N} \\
\varepsilon_a & = \mathcal{T} \begin{bmatrix} \bar{\varepsilon}_{a_1}' & \bar{\varepsilon}_{a_2}' &\ldots & \bar{\varepsilon}_{a_n}'  \end{bmatrix}^{\top} {=} \begin{bmatrix} v_1 & \mathcal{T}_1  \end{bmatrix} \begin{bmatrix} 0 & \textbf{0}_{N \minus 1}^\top  \end{bmatrix}^{\top} = \textbf{0}_{N} \\
\varepsilon_{\tau} & = \mathcal{T} \begin{bmatrix} \bar{\varepsilon}_{\tau_1}' & \bar{\varepsilon}_{\tau_2}' & \ldots & \bar{\varepsilon}_{\tau_n}' \end{bmatrix}^{\top} {=} \begin{bmatrix} v_1 & \mathcal{T}_1  \end{bmatrix} \begin{bmatrix} 0 & \textbf{0}_{N \minus 1}^\top  \end{bmatrix}^{\top} = \textbf{0}_{N} \\ 
\tau & = 1 \tau' = \tau' 
\end{aligned}
\end{equation}
\noindent
Then, since $e = e_1^* \textbf{1}_N$ we have that $e_i = e_k$ for each $i, k \in \nodes$. Since $\tau'$ was not subject to a coordinate change, then the point $(e, \eta, \varepsilon_a, \varepsilon_{\tau}, \tau) = ( e_1^* \textbf{1}_N, \textbf{0}_N, \textbf{0}_N, \textbf{0}_N,  \tau' )$ is an element of $\A_{\varepsilon}$. 

Now, pick a point $z' = (e, \eta, \varepsilon_a, \varepsilon_{\tau}) \in \reals^{4N}$ such that $z' \in \A_{\varepsilon}$. This requires that $e_i = e_k$, $\eta_i = 0$, $\varepsilon_{a_i} = 0$, and $\varepsilon_{\tau_i} = 0$ for each $i,k \in \nodes$; thus, $z' = (e^*, \textbf{0}_N, \textbf{0}_N, \textbf{0}_N)$, where $e^* \in E := \{ e^* \in \reals^n : e^*_i = e^*_k \hspace{2mm} \forall i,k \in \nodes \}$. Then, by left multiplying $(z', \tau')$ by $\Gamma^{\minus 1}$, one has
\begin{equation}
\begin{aligned}
\bar{e} & = \mathcal{T}^{\minus 1} [e^*] = \begin{bmatrix} e^*_1 & \textbf{0}_{N-1}^{\top} \end{bmatrix}^{\top} \\
\bar{\eta} & = \mathcal{T}^{\minus 1} \textbf{0}_{N} = \begin{bmatrix} 0 & \textbf{0}_{N-1}^{\top} \end{bmatrix}^{\top} \\
\bar{\varepsilon}_{a} & = \mathcal{T}^{\minus 1} \textbf{0}_{N} = \begin{bmatrix} 0 & \textbf{0}_{N-1}^{\top} \end{bmatrix}^{\top} \\
\bar{\varepsilon}_{\tau} & = \mathcal{T}^{\minus 1} \textbf{0}_{N} = \begin{bmatrix} 0 & \textbf{0}_{N-1}^{\top} \end{bmatrix}^{\top} \\
 \tau & = 1 \tau = \tau
\end{aligned}
\end{equation}
\noindent \sloppy
giving the point $(\bar{e}, \bar{\eta}, \bar{\varepsilon}_{a}, \bar{\varepsilon}_{\tau}, \tau) = ( e^*_1, \textbf{0}_{N-1}, \textbf{0}_N, \textbf{0}_N, \textbf{0}_N)$. Rearranging the components into the form $(\bar{z}_1, \bar{z}_2, \bar{w}_1, \bar{w}_2, \tau)$ where $\bar{z}_1 = (\bar{e}_1, \bar{\eta}_1)$, $\bar{z}_2 = ( \bar{e}_2, \ldots, \bar{e}_N, \bar{\eta}_2, \ldots, \bar{\eta}_N)$, $\bar{w}_1 = (\bar{\varepsilon}_{a_1},\bar{\varepsilon}_{\tau_1})$, and $\bar{w}_2 = (\bar{\varepsilon}_{a_2},\ldots, \bar{\varepsilon}_{a_n},\bar{\varepsilon}_{\tau_2},\ldots, \bar{\varepsilon}_{\tau_n})$ one has $( e^*_1, 0, \textbf{0}_{N-1}, \textbf{0}_{N-1}, 0, 0,\textbf{0}_{N-1}, \textbf{0}_{N-1})$ which is an element of $\tilde{\A}_{\varepsilon}$.

To relate the set distances between $|x_{\varepsilon}|_{\A_{\varepsilon}}$ and $|\chi_{\varepsilon}|_{\tilde{\A}_{\varepsilon}}$ for every $x_{\varepsilon} \in \mathcal{X}_{\varepsilon}$ and $\chi_{\varepsilon} \in \mathcal{X}_{\varepsilon}$ , note that by definition, one has
$|x_{\varepsilon}|_{\A_{\varepsilon}} = {\rm inf}_{y \in \A_{\varepsilon}} |x_{\varepsilon} - y|$ and
$|\chi_{\varepsilon}|_{\tilde{\A}_{\varepsilon}} = {\rm inf}_{y \in \tilde{\A}_{\varepsilon}} |\chi_{\varepsilon} - y|$, respectively.  Recall that $\chi_{\varepsilon} = \Gamma^{\minus 1} x_{\varepsilon}$ and $x_{\varepsilon} = \Gamma \chi_{\varepsilon}$. Computing the distance $|\chi_{\varepsilon}|_{\tilde{\A}_{\varepsilon}}$, one has
\begin{align*}
|\chi_{\varepsilon}|_{\tilde{\A}_{\varepsilon}} = |\Gamma^{\minus 1} x_{\varepsilon}|_{\tilde{\A}_{\varepsilon}} & = {\rm inf}_{y \in \tilde{\A}_{\varepsilon}} |\Gamma^{\minus 1} x_{\varepsilon}- y| \\
& = {\rm inf}_{e^* \in \reals} |\Gamma^{\minus 1} x_{\varepsilon} \\
& \hspace{-5mm} - (e^*, 0, \textbf{0}_{N-1}, \textbf{0}_{N-1}, 0, 0,\textbf{0}_{N-1}, \textbf{0}_{N-1})|
\end{align*}
\noindent \sloppy
Then, by using the relation $(e^*, 0, \textbf{0}_{N-1},$ $\textbf{0}_{N-1}, 0, 0,\textbf{0}_{N-1}, \textbf{0}_{N-1})$ $= \Gamma^{\minus 1} (e^*\textbf{1}_N, \textbf{0}_{N}, \textbf{0}_{N}, \textbf{0}_{N})$ one has
\begin{align*}
|\Gamma^{\minus 1} x_{\varepsilon}|_{\tilde{\A}_{\varepsilon}} & = {\rm inf}_{e^* \in \reals} |\Gamma^{\minus 1} x_{\varepsilon}- \Gamma^{\minus 1} (e^*\textbf{1}_N, \textbf{0}_{N}, \textbf{0}_{N}, \textbf{0}_{N})| \\
& = {\rm inf}_{e^* \in \reals} \big |\Gamma^{\minus 1} \big ( x_{\varepsilon} - (e^*\textbf{1}_N, \textbf{0}_{N}, \textbf{0}_{N}, \textbf{0}_{N}) \big ) \big | \\
& \leq |\Gamma^{\minus 1} | \Big ( {\rm inf}_{e^* \in \reals} | x_{\varepsilon} - (e^*\textbf{1}_N, \textbf{0}_{N}, \textbf{0}_{N}, \textbf{0}_{N}) | \Big ) \\
& \leq |\Gamma^{\minus 1} | \Big ( {\rm inf}_{y \in \A_{\varepsilon}} | x_{\varepsilon} - y | \Big ) \\
& \leq |\Gamma^{\minus 1} | | x_{\varepsilon} |_{\A_{\varepsilon}} 
\end{align*}
Conversely, computing the distance $|x_{\varepsilon}|_{\A_{\varepsilon}}$, one has
\begin{align*}
|x_{\varepsilon}|_{\A_{\varepsilon}} = |\Gamma \chi_{\varepsilon}|_{\A_{\varepsilon}} & = {\rm inf}_{y \in \A_{\varepsilon}} |\Gamma \chi_{\varepsilon} - y| \\
& = {\rm inf}_{e^* \in \reals} |\Gamma \chi_{\varepsilon} - (e^*\textbf{1}_N, \textbf{0}_{N}, \textbf{0}_{N}, \textbf{0}_{N})|
\end{align*}
\noindent
Then by using the relation $(e^*\textbf{1}_N, \textbf{0}_{N}, \textbf{0}_{N}, \textbf{0}_{N})$ $= \Gamma (e^*, 0, \textbf{0}_{N-1}, \textbf{0}_{N-1}, 0, 0,\textbf{0}_{N-1}, \textbf{0}_{N-1})$, one has
\small
\begin{align*} 
|\Gamma \chi_{\varepsilon}|_{\A_{\varepsilon}} & = {\rm inf}_{e^* \in \reals} |\Gamma \chi_{\varepsilon} { \minus } \Gamma (e^*, 0, \textbf{0}_{N \minus 1}, \textbf{0}_{N \minus 1}, 0, 0,\textbf{0}_{N \minus 1}, \textbf{0}_{N \minus 1})| \\
& = {\rm inf}_{e^* \in \reals} \big |\Gamma \big ( \chi_{\varepsilon} { \minus } (e^*, 0, \textbf{0}_{N \minus 1}, \textbf{0}_{N \minus 1}, 0, 0,\textbf{0}_{N \minus 1}, \textbf{0}_{N \minus 1}) \big ) \big | \\
& \leq |\Gamma | \Big ( {\rm inf}_{e^* \in \reals} | \chi_{\varepsilon} { \minus } (e^*, 0, \textbf{0}_{N-1}, \textbf{0}_{N \minus 1}, 0, 0,\textbf{0}_{N \minus 1}, \textbf{0}_{N \minus 1}) | \Big ) \\
& \leq |\Gamma| \Big ( {\rm inf}_{y \in \tilde{\A}_{\varepsilon}} | \chi_{\varepsilon} - y | \Big ) \\
& \leq |\Gamma| | \chi_{\varepsilon} |_{\tilde{\A}_{\varepsilon}} 
\end{align*} \normalsize
\hfill$\blacksquare$
}


\ifbool{two_col}{
\indent \textit{Proof of Proposition \ref{lem:GES2}}: 
 The proof of this result has been omitted as it follows similarly to that of Lemma \ref{lem:GES1}. 
 }{
\subsection{Proof of Proposition \ref{lem:GES2}}
First, we prove that GES of $\tilde{\A}_{\varepsilon}$ for $\widetilde{\HS}_{\varepsilon}$ implies GES of $\A_{\varepsilon}$ for $\HS_{\varepsilon}$. Suppose the set $\tilde{\A}_{\varepsilon}$ is GES for $\widetilde{\HS}_{\varepsilon}$. By Definition \ref{def:stability}, there exist $\kappa, \alpha > 0$ such that 
\begin{equation} \label{eqn:lem_ges_bnd_1}
| \tilde{\phi}(t,j) |_{\tilde{\A}_{\varepsilon}} \leq \kappa \exp (\minus \alpha (t+ j)) |\tilde{\phi}(0,0)|_{\tilde{\A}_{\varepsilon}} \hspace{5mm} \forall (t,j) \in \mbox{dom } \tilde{\phi}
\end{equation} 
holds for every solution $\tilde{\phi}$ to $\widetilde{\HS}_{\varepsilon}$. Pick a (maximal) solution $\tilde{\phi} \in \mathcal{S}_{\widetilde{\HS}_{\varepsilon}}$ with initial condition $\tilde{\phi}(0,0) \in \widetilde{C}_{\varepsilon} \cup \widetilde{D}_{\varepsilon}$. According to Lemma \ref{lem:equiv_transform}, there exists a maximal solution $\phi$ to $\HS_{\varepsilon}$ such that 
\begin{equation} \label{eqn:gamma_inv_sol_transform}
\tilde{\phi}(t,j) = \Gamma^{\minus 1} \phi(t,j)
\end{equation} 
\noindent
for each $(t,j) \in \mbox{dom } \tilde{\phi}$, where $\Gamma^{\minus 1} = {\rm diag}(\mathcal{T}^{\minus 1}, \mathcal{T}^{\minus 1}, \mathcal{T}^{\minus 1}, \mathcal{T}^{\minus 1}, 1)$. Given that $\tilde{\phi}$ satisfies (\ref{eqn:lem_ges_bnd_1}), applying (\ref{eqn:gamma_inv_sol_transform}) and the relationship between distances in Lemma \ref{lem:set_equiv} given in (\ref{eqn:set_bound_a}) to the right-hand side of (\ref{eqn:lem_ges_bnd_1}), we have that
\begin{equation} \label{eqn:transform_ges_bnd}
\begin{aligned}
| \tilde{\phi}(t,j) |_{\tilde{\A}_{\varepsilon}} & \leq \kappa \exp (\minus \alpha (t+ j)) |\tilde{\phi}(0,0)|_{\tilde{\A}_{\varepsilon}} =  \kappa \exp (\minus \alpha (t+ j)) |\Gamma^{\minus 1} \phi(0,0)|_{\tilde{\A}_{\varepsilon}} \\
& \leq \kappa \exp (\minus \alpha (t+ j)) |\Gamma^{\minus 1}| |\phi(0,0)|_{\A_{\varepsilon}}
\end{aligned}
\end{equation}
By rearranging the relationship given in (\ref{eqn:set_bound_b}), we obtain
\begin{equation}
\frac{1}{|\Gamma|} |x_{\varepsilon}|_{\A_{\varepsilon}} = \frac{1}{|\Gamma|} |\Gamma \chi_{\varepsilon}|_{\A_{\varepsilon}} \leq |\chi_{\varepsilon}|_{\tilde{\A}_{\varepsilon}}
\end{equation}
\noindent
Applying it to the left-hand side of (\ref{eqn:transform_ges_bnd}), we have
\begin{align*}
\frac{1}{|\Gamma|} |\phi(t,j)|_{\A_{\varepsilon}} \leq | \tilde{\phi}(t,j) |_{\tilde{\A}_{\varepsilon}} & \leq \kappa \exp (\minus \alpha (t+ j)) |\Gamma^{\minus 1}| |\phi(0,0)|_{\A_{\varepsilon}}
\end{align*}
Thus, we have that $\phi$ satisfies
\begin{equation}
| \phi(t,j) |_{\A_{\varepsilon}} \leq \tilde{\kappa} \exp (\minus \alpha (t+ j)) |\phi(0,0)|_{\A_{\varepsilon}} \hspace{5mm} \forall (t,j) \in \mbox{dom } \phi
\end{equation} 
where $\tilde{\kappa} = \kappa |\Gamma| |\Gamma^{\minus 1}|$. Then, the set $\A_{\varepsilon}$ is GES for $\HS_{\varepsilon}$.

Conversely, suppose the set $\A_{\varepsilon}$ is GES for $\HS_{\varepsilon}$. By Definition \ref{def:stability}, there exist $\kappa, \alpha > 0$ such that 
\begin{equation} \label{eqn:lem_ges_bnd_2}
| \phi(t,j) |_{\A_{\varepsilon}} \leq \kappa \exp (\minus \alpha (t+ j)) |\phi(0,0)|_{\A_{\varepsilon}} \hspace{5mm} \forall (t,j) \in \mbox{dom } \phi
\end{equation} 
holds for every maximal solution $\phi$ to $\HS_{\varepsilon}$. Pick a maximal solution $\phi \in \mathcal{S}_{\HS_{\varepsilon}}$ with initial condition $\phi(0,0) \in C_{\varepsilon} \cup D_{\varepsilon}$.  According to Lemma \ref{lem:equiv_transform}, there exists a solution $\tilde{\phi}$ to $\widetilde{\HS}_{\varepsilon}$ such that 
\begin{equation} \label{eqn:gamma_sol_transform}
\phi(t,j) = \Gamma \tilde{\phi}(t,j)
\end{equation}
\noindent
for each $(t,j) \in \mbox{dom } \phi$, where $\Gamma = {\rm diag}(\mathcal{T}, \mathcal{T}, \mathcal{T}, \mathcal{T}, 1)$. Given that $\phi$ satisfies (\ref{eqn:lem_ges_bnd_2}), applying (\ref{eqn:gamma_sol_transform}) and the relationship between distances in Lemma \ref{lem:set_equiv} to the right-hand side of (\ref{eqn:lem_ges_bnd_1}), we have that
\begin{equation} \label{eqn:transform_ges_bnd_2}
\begin{aligned}
| \phi(t,j) |_{\A_{\varepsilon}} & \leq \kappa \exp (\minus \alpha (t+ j)) |\phi(0,0)|_{\A_{\varepsilon}} =  \kappa \exp (\minus \alpha (t+ j)) |\Gamma \tilde{\phi}(0,0)|_{\A_{\varepsilon}} \\
& \leq \kappa \exp (\minus \alpha (t+ j)) |\Gamma| |\tilde{\phi}(0,0)|_{\tilde{\A}_{\varepsilon}} 
\end{aligned}
\end{equation}
By rearranging the relationship given in (\ref{eqn:set_bound_a}), we obtain
\begin{equation}
\frac{1}{|\Gamma^{\minus 1}|} |\chi_{\varepsilon}|_{\tilde{\A}_{\varepsilon}} = \frac{1}{|\Gamma^{\minus 1}|} |\Gamma^{\minus 1} x_{\varepsilon}|_{\tilde{\A}_{\varepsilon}} \leq |x_{\varepsilon}|_{\A_{\varepsilon}}
\end{equation}
\noindent
Applying it to the left-hand side of (\ref{eqn:transform_ges_bnd_2}), we have
\begin{align*}
\frac{1}{|\Gamma^{\minus 1}|} | \tilde{\phi}(t,j)|_{\A_{\varepsilon}} \leq |\phi(t,j) |_{\A_{\varepsilon}} & \leq \kappa \exp (\minus \alpha (t+ j)) |\Gamma| |\tilde{\phi}(0,0)|_{\A_{\varepsilon}}
\end{align*}
Thus, we have that $\tilde{\phi}$ satisfies
\begin{equation}
| \tilde{\phi}(t,j) |_{\tilde{\A}_{\varepsilon}} \leq \kappa' \exp (\minus \alpha (t+ j)) |\tilde{\phi}(0,0)|_{\tilde{\A}_{\varepsilon}} \hspace{5mm} \forall (t,j) \in \mbox{dom } \tilde{\phi}
\end{equation}  
where $\kappa' = \kappa |\Gamma^{\minus 1}| |\Gamma|$. Then, the set $\tilde{\A}_{\varepsilon}$ is GES for $\widetilde{\HS}_{\varepsilon}$.
\hfill$\blacksquare$
}


\ifbool{two_col}{
\indent \textit{Proof of Proposition \ref{prop:est_param}}:
Consider the following Lyapunov function candidate $V_{\varepsilon_r}(\chi_{\varepsilon_r}) = \bar{w}_1^{\top} P_2 \bar{w}_1 + \bar{w}_2^{\top} P_3 \bar{w}_2$ It satisfies
\ifbool{two_col}{$\alpha_{\bar{\omega}_1} |\chi_{\varepsilon_r}|_{\tilde{\A}_{\varepsilon_r}}^2 \leq V_{\varepsilon_r}(\chi_{\varepsilon_r}) \leq \alpha_{\bar{\omega}_2} |\chi_{\varepsilon_r}|_{\tilde{\A}_{\varepsilon_r}}^2$}{
\begin{equation} \label{eqn:Vw_def}
\alpha_{\bar{\omega}_1} |\chi_{\varepsilon_r}|_{\tilde{\A}_{\varepsilon_r}}^2 \leq V(\chi_{\varepsilon_r}) \leq \alpha_{\bar{\omega}_2} |\chi_{\varepsilon_r}|_{\tilde{\A}_{\varepsilon_r}}^2
\end{equation}}
\noindent \sloppy
for each $\chi_{\varepsilon_r} \in$ $\tilde{C}_{\varepsilon_r} \cup \tilde{D}_{\varepsilon_r}$ with $\alpha_1 =$ $\mathrm{min} \big \{ \lambda_{min} (P_2),$ $\lambda_{min} (P_3) \}$ and $\alpha_2 =$ $\mathrm{max} \big \{ \lambda_{max} (P_2),$ $\lambda_{max} (P_3) \}$. For each $\chi_{\varepsilon_r} \in \tilde{C}_{\varepsilon_r}$ $\langle \nabla V_{\varepsilon_r}(\chi_{\varepsilon_r}), \tilde{f}(\chi_{\varepsilon_r}) \rangle$ $\leq \bar{w}_1^{\top} (P_2 A_{f_3}$ $+ A_{f_3}^{\top} P_2) \bar{w}_1$ $+ \bar{w}_2^{\top} (P_3 A_{f_4}$ $+ A_{f_4}^{\top} P_3) \bar{w}_2$
\noindent
The conditions in (\ref{eqn:lyap1}) imply the existence of positive numbers $\beta_1$ and $\beta_2$ such that $P_2 A_{f_3} {+} A_{f_3}^{\top} P_2 \leq - \beta_1 I$ $P_3 A_{f_4} {+} A_{f_4}^{\top} P_3 \leq - \beta_2 I$ Then
\ifbool{two_col}{$\langle \nabla V_{\varepsilon_r}(\chi_{\varepsilon_r}), \tilde{f}_{\varepsilon_r}(\chi_{\varepsilon_r}) \rangle \leq - \frac{\tilde{\beta}}{\alpha_{\bar{\omega}_2}} V_{\varepsilon_r}(\chi_{\varepsilon_r})$}{
\vspace{-1mm}
\begin{equation} \label{eqn:Vw_flows}
\begin{aligned}
\langle \nabla V_{\varepsilon_r}(\chi_{\varepsilon_r}), \tilde{f}_{\varepsilon_r}(\chi_{\varepsilon_r}) \rangle \leq - \frac{\tilde{\beta}}{\alpha_{\bar{\omega}_2}} V_{\varepsilon_r}(\chi_{\varepsilon_r})
\end{aligned}
\end{equation}
\noindent}
where  $\tilde{\beta} = \min \{ \beta_1, \beta_2 \} > 0$. For all $\chi_{\varepsilon_r} \in \tilde{D}_{\varepsilon_r}$ and $g \in \tilde{G}_{\varepsilon_r} (\chi_{\varepsilon_r})$
\ifbool{two_col}{, $V_{\varepsilon_r}(g) - V_{\varepsilon_r}(\chi_{\varepsilon_r}) = 0$.}{
\vspace{-2mm}
\begin{equation} \label{eqn:Vw_jumps}
\begin{aligned}
V_{\varepsilon_r}(g) - V_{\varepsilon_r}(\chi_{\varepsilon_r}) = 0
\end{aligned}
\end{equation}
\noindent}
 Now, pick a solution $\tilde{\phi}$ to $\widetilde{\HS}_{\varepsilon_r}$ with initial condition $\tilde{\phi}(0,0) \in \widetilde{C}_{\varepsilon_r} \cup \widetilde{D}_{\varepsilon_r}$. 
\ifbool{two_col}{Direct}{ 
 As a result of (\ref{eqn:Vw_flows}) and (\ref{eqn:Vw_jumps}), direct} 
integration of $(t,j) \mapsto V_{\varepsilon_r} (\tilde{\phi}(t,j))$ over $\mbox{dom } \tilde{\phi}$ gives 
$V_{\varepsilon_r}( \tilde{\phi} (t,j)) \leq \exp{ \big (\minus \frac{ \tilde{\beta}}{\alpha_{\bar{\omega}_2}} t \big )} V_{\varepsilon_r}( \tilde{\phi} (0,0))$ for each $(t,j) \in \mbox{dom } \tilde{\phi}$. Now, given the relation established in (\ref{eqn:t_upper_bound}), for any solution $\tilde{\phi}$ to $\widetilde{\HS}_{\varepsilon_r}$, we have $j T_2 \leq t \Rightarrow -t \leq -j T_2$. Then, for any $\gamma \in (0,1)$ we have $- \gamma t \leq -\gamma T_2 j$. Moreover, $- t = -(1 - \gamma) t - \gamma t$ $\leq -(1-\gamma)t - \gamma T_2 j$ $\leq - \min \{1-\gamma, \gamma T_2 \}(t+j)$
\noindent
leading to $V_{\varepsilon_r}( \tilde{\phi} (t,j)) \leq \exp{ \Big (\minus \frac{ \bar{\gamma}  \tilde{\beta}}{\alpha_{\bar{\omega}_2}} ( t + j ) \Big )} V_{\varepsilon_r}( \tilde{\phi} (0,0))$
\noindent
for each $(t,j) \in \mbox{dom } \tilde{\phi}$ where $\bar{\gamma} = \min \{1-\gamma, \gamma T_2\}$. 
\ifbool{two_col}{{Then, by combining this inequality with the definition of $V_{\varepsilon_r}$,}}{Then, by combining this inequality with (\ref{eqn:Vw_def}),} one has $\alpha_{\bar{\omega}_1} |\chi_{\varepsilon_r}|_{\tilde{\A}_{\varepsilon_r}}^2 {\leq} V_{\varepsilon_r}( \tilde{\phi}  (t,j)) \leq \exp{ \Big ( \minus \frac{ \bar{\gamma}   \tilde{\beta}}{\alpha_{\bar{\omega}_2}}  ( t + j )  \Big ) } V_{\varepsilon_r}( \tilde{\phi}  (0,0))$
\noindent
then leveraging $V_{\varepsilon_r}( \tilde{\phi}  (0,0)) \leq \alpha_{\bar{\omega}_2} | \tilde{\phi}  (0,0)|_{\tilde{\A}_{\varepsilon_r}}^2$ we have 
\ifbool{two_col}{
$| \tilde{\phi}  (t,j)|_{\tilde{\A}_{\varepsilon_r}} \leq \sqrt{\frac{\alpha_{\bar{\omega}_2}}{\alpha_{\bar{\omega}_1}}} \exp{ \Big ( \minus \frac{ \bar{\gamma}   \tilde{\beta}}{2 \alpha_{\bar{\omega}_2}}  ( t + j )  \Big ) } | \tilde{\phi}  (0,0)|_{\tilde{\A}_{\varepsilon_r}}$}{
\vspace{-3mm}
\begin{equation*}
| \tilde{\phi}  (t,j)|_{\tilde{\A}_{\varepsilon_r}} \leq \sqrt{\frac{\alpha_{\bar{\omega}_2}}{\alpha_{\bar{\omega}_1}}} \exp{ \Big ( \minus \frac{ \bar{\gamma}   \tilde{\beta}}{2 \alpha_{\bar{\omega}_2}}  ( t + j )  \Big ) } | \tilde{\phi}  (0,0)|_{\tilde{\A}_{\varepsilon_r}}
\end{equation*}
\noindent}
Observe that this bound holds for each solution $\tilde{\phi}$ to $\widetilde{\HS}_{\varepsilon_r}$. Maximal solutions to $\widetilde{\HS}_{\varepsilon_r}$ are complete due to the reduction property established in Lemmas \ref{lem:equiv_transform}, \ref{lem:equiv}, and \ref{lem:complete}. \ifbool{two_col}{\hspace{-2mm}}{In particular, Lemma \ref{lem:equiv_transform} establishes the relation between  $\widetilde{\HS}_{\varepsilon}$ and $\HS_{\varepsilon}$, Lemma \ref{lem:equiv} establishes the reduction from $\HS$ to $\HS_{\varepsilon}$, the former for which we have established completeness of solutions in Lemma \ref{lem:complete}.} Therefore, the set $\tilde{\A}_{\varepsilon_r}$ is globally exponentially stable for $\widetilde{\HS}_{\varepsilon_r}$.
\hfill$\blacksquare$
}{
\subsection{Proof of Proposition \ref{prop:est_param}}
\begin{equation}
V_{\varepsilon_r}(\chi_{\varepsilon_r}) = \bar{w}_1^{\top} P_2 \bar{w}_1 + \bar{w}_2^{\top} P_3 \bar{w}_2
\end{equation}
\noindent
It satisfies
\begin{equation} \label{eqn:Vw_def}
\alpha_{\bar{\omega}_1} |\chi_{\varepsilon_r}|_{\tilde{\A}_{\varepsilon_r}}^2 \leq V(\chi_{\varepsilon_r}) \leq \alpha_{\bar{\omega}_2} |\chi_{\varepsilon_r}|_{\tilde{\A}_{\varepsilon_r}}^2 \hspace{5mm} \forall \chi_{\varepsilon_r} \in \tilde{C}_{\varepsilon_r} \cup \tilde{D}_{\varepsilon_r}
\end{equation}
\noindent
with $\alpha_1 = \mathrm{min} \big \{ \lambda_{min} (P_2), \lambda_{min} (P_3) \}$ and $\alpha_2 = \mathrm{max} \big \{ \lambda_{max} (P_2), \lambda_{max} (P_3) \}$. For each $\chi_{\varepsilon_r} \in \tilde{C}_{\varepsilon_r}$ 
\begin{equation}
\begin{aligned}
\langle \nabla V_{\varepsilon_r}(\chi_{\varepsilon_r}), \tilde{f}(\chi_{\varepsilon_r}) \rangle & \leq \bar{w}_1^{\top} (P_2 A_{f_3} + A_{f_3}^{\top} P_2) \bar{w}_1 \\ & \hspace{5mm} + \bar{w}_2^{\top} (P_3 A_{f_4} + A_{f_4}^{\top} P_3) \bar{w}_2 \\ 
\end{aligned}
\end{equation}
\noindent
The conditions in (\ref{eqn:lyap1}) imply the existence of positive numbers $\beta_1$ and $\beta_2$ such that $$P_2 A_{f_3} {+} A_{f_3}^{\top} P_2 \leq - \beta_1 I$$ $$P_3 A_{f_4} {+} A_{f_4}^{\top} P_3 \leq - \beta_2 I$$ Then
\begin{equation} \label{eqn:Vw_flows}
\begin{aligned}
\langle \nabla V_{\varepsilon_r}(\chi_{\varepsilon_r}), \tilde{f}_{\varepsilon_r}(\chi_{\varepsilon_r}) \rangle & \leq -\beta_1 |\bar{w}_1|^2 -\beta_2 |\bar{w}_2|^2 \\
& \leq - \tilde{\beta} \big ( |\bar{w}_1|^2 + |\bar{w}_2|^2 \big ) \\
& \leq - \tilde{\beta}  \big ( |\chi_{\varepsilon_r}|_{\tilde{\A}_{\varepsilon_r}}^2 \big )  \\
& \leq - \frac{\tilde{\beta}}{\alpha_{\bar{\omega}_2}} V_{\varepsilon_r}(\chi_{\varepsilon_r})
\end{aligned}
\end{equation}
\noindent
where  $\tilde{\beta} = \min \{ \beta_1, \beta_2 \} > 0$. For all $\chi_{\varepsilon_r} \in \tilde{D}_{\varepsilon_r}$ and $g \in \tilde{G}_{\varepsilon_r} (\chi_{\varepsilon_r})$
\begin{equation} \label{eqn:Vw_jumps}
\begin{aligned}
V_{\varepsilon_r}(g) - V_{\varepsilon_r}(\chi_{\varepsilon_r}) = 0
\end{aligned}
\end{equation}
\noindent
 Now, pick a solution $\tilde{\phi}$ to $\widetilde{\HS}_{\varepsilon_r}$ with initial condition $\tilde{\phi}(0,0) \in \widetilde{C}_{\varepsilon_r} \cup \widetilde{D}_{\varepsilon_r}$. As a result of (\ref{eqn:Vw_flows}) and (\ref{eqn:Vw_jumps}), direct integration of $(t,j) \mapsto V_{\varepsilon_r} (\tilde{\phi}(t,j))$ over $\mbox{dom } \tilde{\phi}$ gives 
\begin{equation} 
\begin{aligned}
V_{\varepsilon_r}( \tilde{\phi} (t,j)) \leq \exp{ \Big (\minus \frac{ \tilde{\beta}}{\alpha_{\bar{\omega}_2}} t \Big )} V_{\varepsilon_r}( \tilde{\phi} (0,0)) & & \forall (t,j) \in \mbox{dom } \tilde{\phi}
\end{aligned}
\end{equation}
 Now, given the relation established in (\ref{eqn:t_upper_bound}), for any solution $\tilde{\phi}$ to $\widetilde{\HS}_{\varepsilon_r}$, we have $j T_2 \leq t \Rightarrow -t \leq -j T_2$. Then, for any $\gamma \in (0,1)$ we have $- \gamma t \leq -\gamma T_2 j$. Moreover, 
\begin{equation} \label{eqn:t+j_param}
\begin{aligned}
- t = -(1 - \gamma) t - \gamma t & \leq -(1-\gamma)t - \gamma T_2 j \\ & \leq - \min \{1-\gamma, \gamma T_2 \}(t+j) 
\end{aligned}
\end{equation}
\noindent
leading to
\begin{equation} 
\begin{aligned}
V_{\varepsilon_r}( \tilde{\phi} (t,j)) & \leq \exp{ \Big (\minus \frac{ \bar{\gamma}  \tilde{\beta}}{\alpha_{\bar{\omega}_2}} ( t + j ) \Big )} V_{\varepsilon_r}( \tilde{\phi} (0,0))
\end{aligned}
\end{equation}
\noindent
for each $(t,j) \in \mbox{dom } \tilde{\phi}$ where $\bar{\gamma} = \min \{1-\gamma, \gamma T_2\}$. Then, by combining this inequality with (\ref{eqn:Vw_def}), one has
\begin{equation}
\begin{aligned}
\alpha_{\bar{\omega}_1} |\chi_{\varepsilon_r}|_{\tilde{\A}_{\varepsilon_r}}^2 {\leq} V_{\varepsilon_r}( \tilde{\phi}  (t,j)) \leq \exp{ \Big ( \minus \frac{ \bar{\gamma}   \tilde{\beta}}{\alpha_{\bar{\omega}_2}}  ( t + j )  \Big ) } V_{\varepsilon_r}( \tilde{\phi}  (0,0))
\end{aligned}
\end{equation}
\noindent
then leveraging $V_{\varepsilon_r}( \tilde{\phi}  (0,0)) \leq \alpha_{\bar{\omega}_2} | \tilde{\phi}  (0,0)|_{\tilde{\A}_{\varepsilon_r}}^2$ we have 
\begin{equation}
\begin{aligned}
| \tilde{\phi}  (t,j)|_{\tilde{\A}_{\varepsilon_r}}^2 & \leq \frac{\alpha_{\bar{\omega}_2}}{\alpha_{\bar{\omega}_1}}  \exp{\Big ( \minus \frac{ \bar{\gamma}  \tilde{\beta}}{\alpha_{\bar{\omega}_2}}  ( t + j )  \Big ) } | \tilde{\phi} (0,0)|_{\tilde{\A}_{\varepsilon_r}}^2 \\
\end{aligned}
\end{equation}
\noindent
then
\begin{equation}
| \tilde{\phi}  (t,j)|_{\tilde{\A}_{\varepsilon_r}} \leq \sqrt{\frac{\alpha_{\bar{\omega}_2}}{\alpha_{\bar{\omega}_1}}} \exp{ \Big ( \minus \frac{ \bar{\gamma}   \tilde{\beta}}{2 \alpha_{\bar{\omega}_2}}  ( t + j )  \Big ) } | \tilde{\phi}  (0,0)|_{\tilde{\A}_{\varepsilon_r}}
\end{equation}
\noindent
Observe that this bound holds for each solution $\tilde{\phi}$ to $\widetilde{\HS}_{\varepsilon_r}$. Maximal solutions to $\widetilde{\HS}_{\varepsilon_r}$ are complete due to the reduction property established in Lemmas \ref{lem:equiv_transform}, \ref{lem:equiv}, and \ref{lem:complete}. In particular, Lemma \ref{lem:equiv_transform} establishes the relation between  $\widetilde{\HS}_{\varepsilon}$ and $\HS_{\varepsilon}$, Lemma \ref{lem:equiv} establishes the reduction from $\HS$ to $\HS_{\varepsilon}$, the former for which we have established completeness of solutions in Lemma \ref{lem:complete}. Therefore, the set $\tilde{\A}_{\varepsilon_r}$ is globally exponentially stable for $\widetilde{\HS}_{\varepsilon_r}$.
\hfill$\blacksquare$
}

\ifbool{two_col}{}{
\indent \textit{Proof of Theorem \ref{thm:iss2}}:
Consider the same Lyapunov function candidate $V(\chi_m) = V_1(\chi_m) + V_2(\chi_m) + V_{\varepsilon_r} (\chi_m)$ from the proof of Theorem \ref{thrm1}  in Section \ref{sec:thm_proof}.  During flows, there is no contribution from the perturbation thus the derivative of $V$ is unchanged from the proof of Theorem \ref{thrm1}. Thus, one has  that (\ref{eqn:change_in_v}) holds with $\widetilde{f}_{\varepsilon}(\chi_{\varepsilon})$ replaced by $\widetilde{f}_m(\chi_m)$, namely, 
\ifbool{rep}{
\begin{equation*}
\begin{aligned}
\langle \nabla V(\chi_m), \widetilde{f}_m(\chi_m) \rangle & \leq 2 \bar{z}_2^{\top} \big ( \exp (A^{\top}_{f_2} \tau ) P \exp (A_{f_2} \tau) \big ) B_{f_2}  \bar{w}_2 + \bar{w}_1^{\top} (P_1 A_{f_3} + A_{f_3}^{\top} P_1) \bar{w}_1 + \bar{w}_2^{\top} (P_2 A_{f_4} + A_{f_4}^{\top} P_2) \bar{w}_2 \\ 
\end{aligned}
\end{equation*}}{
\begin{equation*}
\begin{aligned}
\langle \nabla V(\chi_m), \tilde{f}(\chi_m) \rangle & \leq 2 \bar{z}_2^{\top} \big ( \exp (A^{\top}_{f_2} \tau ) P \exp (A_{f_2} \tau) \big ) B_{f_2}  \bar{w}_2 \\ 
& \hspace{3mm} + \bar{w}_1^{\top} (P_1 A_{f_3} + A_{f_3}^{\top} P_1) \bar{w}_1 \\
& \hspace{5mm} + \bar{w}_2^{\top} (P_2 A_{f_4} + A_{f_4}^{\top} P_2) \bar{w}_2
\end{aligned}
\end{equation*}}
\noindent
then by following the same notions of the proof in Theorem \ref{thrm1}, one has 
$\langle \nabla V(\chi_m), \widetilde{f}(\chi_m) \rangle \leq \frac{\bar{\kappa}_1}{\alpha_2} V(\chi_m)$  where $\bar{\kappa}_1 = \max \Big \{ \frac{\kappa_1}{2 \epsilon} , \Big ( \frac{\kappa_1 \epsilon}{2} - \beta_2 \Big ) \Big \}$ and $\varepsilon > 0$. At jumps, triggered when $\tau = 0$, one has, for each $\chi_m \in \widetilde{D}_m \setminus \widetilde{\A}_{\varepsilon}$ and $g \in \widetilde{G}_m(\chi_m)$
\ifbool{rep}{
\begin{equation} \label{thrm_iss_jumps}
\begin{aligned}
V(g) {-} V(\chi_m) & \leq {-} \bar{\eta}_1^2 + (A_{g_2} \bar{z}_2 - B_g \bar{m}_{\bar{z}_2})^{\top} Q (A_{g_2} \bar{z}_2 - B_g \bar{m}_{\bar{z}_2}) - \bar{z}_2^{\top} P_1 \bar{z}_2 \\
& \leq {-} \bar{\eta}_1^2 + (A_{g_2} \bar{z}_2)^{\top} \exp{A^{\top}_{f_2} \tau } P_1 \exp{A_{f_2} \tau} (A_{g_2} \bar{z}_2) - 2 (B_g \bar{m}_{\bar{z}_2})^{\top} \exp{A^{\top}_{f_2} \tau } P_1 \exp{A_{f_2} \tau} (A_{g_2} \bar{z}_2) \\ 
& \hspace{60mm} + (B_g \bar{m}_{\bar{z}_2})^{\top} \exp{A^{\top}_{f_2} \tau } P_1 \exp{A_{f_2} \tau} (B_g \bar{m}_{\bar{z}_2}) {-}  \bar{z}_2^{\top} P_1 \bar{z}_2 \\ 
\end{aligned}
\end{equation}}{
\small
\begin{equation} \label{thrm_iss_jumps}
\begin{aligned}
V(g) {\minus} V(\chi_m) 
& \leq \\
& \hspace{-6mm} - \bar{\eta}_1^2 + (A_{g_2} \bar{z}_2)^{\top} \exp (A^{\top}_{f_2} \tau ) P_1 \exp (A_{f_2} \tau) (A_{g_2} \bar{z}_2) \\
& \hspace{-3mm} {-} 2 (B_g \bar{m}_{\bar{z}_2})^{\top} \exp (A^{\top}_{f_2} \tau ) P_1 \exp (A_{f_2} \tau) (A_{g_2} \bar{z}_2) \\ 
& \hspace{0mm} + (B_g \bar{m}_{\bar{z}_2})^{\top} \exp (A^{\top}_{f_2} \tau ) P_1 \exp (A_{f_2} \tau) (B_g \bar{m}_{\bar{z}_2}) \\
& \hspace{3mm} {-}  \bar{z}_2^{\top} P_1 \bar{z}_2 \\ 
\end{aligned}
\end{equation}}
\normalsize
\noindent
From (\ref{cond:phil_1}) and the proof in Theorem \ref{thrm1}, there exists a scalar $\kappa_2$ such that $\bar{z}_2^{\top} (A_{g_2}^{\top} \exp(A^{\top}_{f_2} v ) P_1 \exp(A_{f_2} v) A_{g_2} - P_1) \bar{z}_2 \leq - \kappa_2 \bar{z}_2^{\top} \bar{z}_2$ leading to
\ifbool{rep}{
\begin{equation} \label{thrm_iss_jumps1}
\begin{aligned}
V(g) {-} V(\chi_m) & \leq - \bar{\eta}_1^2 -\kappa_2 \bar{z}_2^{\top} \bar{z}_2 - 2 (B_g \bar{m}_{\bar{z}_2})^{\top} \exp{A^{\top}_{f_2} \tau } P_1 \exp{A_{f_2} \tau} (A_{g_2} \bar{z}_2) + (B_g \bar{m}_{\bar{z}_2})^{\top} \exp{A^{\top}_{f_2} \tau } P_1 \exp{A_{f_2} \tau} (B_g \bar{m}_{\bar{z}_2})
\end{aligned}
\end{equation}}{
\small
\begin{equation} \label{thrm_iss_jumps1}
\begin{aligned}
V(g) {-} V(\chi_m) & \leq - \bar{\eta}_1^2 -\kappa_2 \bar{z}_2^{\top} \bar{z}_2 \\
& \hspace{-4mm} - 2 (B_g \bar{m}_{\bar{z}_2})^{\top} \exp(A^{\top}_{f_2} \tau ) P_1 \exp (A_{f_2} \tau) (A_{g_2} \bar{z}_2) \\
& \hspace{-1mm} + (B_g \bar{m}_{\bar{z}_2})^{\top} \exp (A^{\top}_{f_2} \tau) P_1 \exp(A_{f_2} \tau) (B_g \bar{m}_{\bar{z}_2})
\end{aligned}
\end{equation}}
\normalsize
\noindent
Let $Q = \exp(A^{\top}_{f_2} \tau ) P_1 \exp(A_{f_2} \tau)$, then applying Young's inequality on the third term such that 
\ifbool{rep}{
\begin{equation*}
\begin{aligned}
\bar{m}_{\bar{z}_2}^{\top} B_g^{\top} Q A_{g_2} \bar{z}_2 & \leq \frac{1}{2 \epsilon_2} \Big ( \bar{m}_{\bar{z}_2}^{\top} B_g^{\top} Q A_{g_2} \Big )^{\top} \Big ( \bar{m}_{\bar{z}_2}^{\top} B_g^{\top} Q A_{g_2} \Big ) + \frac{\epsilon_2}{2} \bar{z}_2^{\top} \bar{z}_2 \\
& \leq \frac{1}{2 \epsilon_2} \Big | \big ( B_g^{\top} Q A_{g_2} \big ) \big ( B_g^{\top} Q A_{g_2} \big )^{\top} \Big | \bar{m}_{\bar{z}_2}^{\top} \bar{m}_{\bar{z}_2} + \frac{\epsilon_2}{2} \bar{z}_2^{\top} \bar{z}_2
\end{aligned}
\end{equation*}}{
\small
\begin{equation*}
\begin{aligned}
\bar{m}_{\bar{z}_2}^{\top} B_g^{\top} Q A_{g_2} \bar{z}_2
& \leq \frac{\big | \big ( B_g^{\top} Q A_{g_2} \big ) \big ( B_g^{\top} Q A_{g_2} \big )^{\top} \big |}{2 \epsilon_2}  \bar{m}_{\bar{z}_2}^{\top} \bar{m}_{\bar{z}_2} {+} \frac{\epsilon_2}{2} \bar{z}_2^{\top} \bar{z}_2
\end{aligned}
\end{equation*}}
\normalsize
\noindent
where $\epsilon_2 > 0$. Then, we have
\ifbool{rep}{
\begin{equation} \label{thrm_iss_jumps2}
\begin{aligned}
V(g) \minus V(\chi_m) & \leq - \bar{\eta}_1^2 -\kappa_2 \bar{z}_2^{\top} \bar{z}_2 {-} \Big ( {\small \frac{1}{2 \epsilon_2}} \big | \big ( B_g^{\top} Q A_{g_2} \big ) \big ( B_g^{\top} Q A_{g_2} \big )^{\top} \big | \bar{m}_{\bar{z}_2}^{\top} \bar{m}_{\bar{z}_2} + \frac{\epsilon_2}{2} \bar{z}_2^{\top} \bar{z}_2 \Big ) \\ 
& \hspace{15mm} + \bar{m}_{\bar{z}_2} B_g^{\top} Q B_g \bar{m}_{\bar{z}_2} \\
& \leq - \bar{\eta}_1^2 - \big ( \kappa_2 + \frac{\epsilon_2}{2} \big ) \bar{z}_2^{\top} \bar{z}_2 {+} \big ( | B_g^{\top} Q B_g | \minus \frac{1}{2 \epsilon_2} | ( B_g^{\top} Q A_{g_2} )( B_g^{\top} Q A_{g_2} )^{\top} | \big ) \bar{m}_{\bar{z}_2}^{\top}  \bar{m}_{\bar{z}_2} 
\end{aligned}
\end{equation}}{
\small
\begin{equation} \label{thrm_iss_jumps2}
\begin{aligned}
V(g) \minus V(\chi_m)
& \leq - \bar{\eta}_1^2 - \big ( \kappa_2 + \frac{\epsilon_2}{2} \big ) \bar{z}_2^{\top} \bar{z}_2 {+} \big ( | B_g^{\top} Q B_g | \\
& \hspace{5mm} \minus \frac{1}{2 \epsilon_2} | ( B_g^{\top} Q A_{g_2} )( B_g^{\top} Q A_{g_2} )^{\top} | \big ) \bar{m}_{\bar{z}_2}^{\top}  \bar{m}_{\bar{z}_2} 
\end{aligned}
\end{equation}}
\normalsize
\noindent
By noting $|A_{g_2}|, |B_g|$ $\leq \gamma \lambda_{max}({\bar{\mathcal{L}}})$ let $\kappa_{\bar{m}_2} =$ $\big ( \lambda_{max}({\bar{\mathcal{L}}}) \big )^2 \underset{v \in [0,T_2]}{\mathrm{max}} \Big \{ \lambda_{max} \big (\exp (A^{\top}_{f_2} v ) P_1 \exp (A_{f_2} v ) \big ) \Big \}$, we let $\epsilon_2 = \kappa_2$ and obtain
\ifbool{rep}{
\begin{equation*} 
\begin{aligned}
V(g) - V(\chi_m) & \leq - \bar{\eta}_1^2 - \big ( \kappa_2 + \frac{\kappa_2}{2} \big ) \bar{z}_2^{\top} \bar{z}_2 +  \big ( \gamma^2 \kappa_{\bar{m}_2} - \frac{1}{2 \kappa_2} \gamma^4 \kappa_{\bar{m}_2}^2 \big ) \bar{m}_{\bar{z}_2}^{\top}  \bar{m}_{\bar{z}_2} 
\end{aligned}
\end{equation*}}{
\small
\begin{equation*} 
\begin{aligned}
V(g) \minus V(\chi_m) {\leq} \minus \bar{\eta}_1^2 \minus \frac{3\kappa_2}{2} \bar{z}_2^{\top} \bar{z}_2 {+}  \big ( \gamma^2 \kappa_{\bar{m}_2} {-} \frac{\gamma^4 \kappa_{\bar{m}_2}^2}{2 \kappa_2}   \big ) \bar{m}_{\bar{z}_2}^{\top}  \bar{m}_{\bar{z}_2} 
\end{aligned}
\end{equation*}}
\normalsize
\noindent
Now, let $\tilde{\kappa}_{\bar{m}_2} = \gamma^2 \kappa_{\bar{m}_2} - \frac{1}{2 \kappa_2} \gamma^4 \kappa_{\bar{m}_2}^2$ then at jumps one has
\begin{equation} \label{thrm_iss_jumps3}
\begin{aligned}
V(g) - V(\chi_m) & \leq - \bar{\kappa}_2 (|\bar{\eta}_1|^2 + |\bar{z}_2|^2) + \tilde{\kappa}_{\bar{m}_2} |\bar{m}_{\bar{z}_2}|^2 
\end{aligned}
\end{equation}
\noindent
where $\bar{\kappa}_2 = \max \big \{1, \frac{3 \kappa_2}{2} \big \}$.   Now, recall from (\ref{thrm_zdef}) in the proof of Theorem \ref{thrm1} that
\ifbool{two_col}{$-(|\bar{\eta}_1|^2 + |\bar{z}_2|^2) \leq - \frac{1}{\alpha_2} V(\chi_{\varepsilon}) + |\bar{w}|^2$.}{
\begin{equation}
\begin{aligned}
-(|\bar{\eta}_1|^2 + |\bar{z}_2|^2) \leq - \frac{1}{\alpha_2} V(\chi_{\varepsilon}) + |\bar{w}|^2
\end{aligned}
\end{equation}}
\noindent
Then, plugging (\ref{thrm_zdef}) into (\ref{thrm_iss_jumps3}) one has
\small
\begin{equation*}
\begin{aligned}
V(g) & \leq \Big ( 1 - \frac{3 \kappa_2}{2 \alpha_2} \Big ) V(\chi_m) + \frac{3 \kappa_2}{2} |\bar{w}|^2 + \tilde{\kappa}_{\bar{m}_2} |\bar{m}_{\bar{z}_2}|^2 
\end{aligned}
\end{equation*}
\normalsize
\noindent
Noting $\langle \nabla V(\chi_{\varepsilon}), \widetilde{f}(\chi_{\varepsilon}) \rangle  \leq \frac{\bar{\kappa}_1}{\alpha_2} V(\chi_{\varepsilon})$, one can then pick a solution with initial conditions $\tilde{\phi}(0,0) \in \widetilde{C}_m \cup \widetilde{D}_m$ and find that the trajectory of $V(\tilde{\phi}(t,j))$ is bounded as follows:
\small
\begin{equation*}
\begin{aligned}
V(\tilde{\phi}(t,j)) & \leq \\
& \hspace{0mm} \exp \big ( \frac{\bar{\kappa}_1}{\alpha_2} \hspace{0.5mm} T_2 \big ) \Big ( \exp \big ( \frac{\bar{\kappa}_1}{\alpha_2} \hspace{0.5mm} T_2 \big ) \Big ( 1 - \frac{3 \kappa_2}{2 \alpha_2} \Big ) \Big )^j V(\tilde{\phi}(0,0)) \\
& \hspace{5mm} + \frac{3 \kappa_2}{2} \exp \big ( \bar{\kappa} \hspace{0.5mm} T_2 \big ) \mathrm{sup}_{(t,j) \in \mbox{\footnotesize dom} \phi} |\bar{w}(t,j)|^2 \\
& \hspace{10mm} + \tilde{\kappa}_{\bar{m}_2} \exp \Big (\frac{\kappa}{2 \epsilon_2} T_2 \Big ) \mathrm{sup}_{(t,j) \in \mbox{\footnotesize dom} \phi} |\bar{m}_{\bar{z}_2}|^2
\end{aligned}
\end{equation*}
\normalsize
\noindent
The result follows from an analysis of $V(\tilde{\phi}(t,j))$ over $\mbox{dom } \tilde{\phi}$ utilizing the same approach as in the proof of Theorem \ref{thrm1}.
\hfill$\blacksquare$}

\nocite{*}
\bibliography{Clock_Cite}
\bibliographystyle{ieeetr}


\ifbool{rep}{}
{
\begin{IEEEbiography}[{\includegraphics[width=1in,height=1.25in,clip,keepaspectratio]{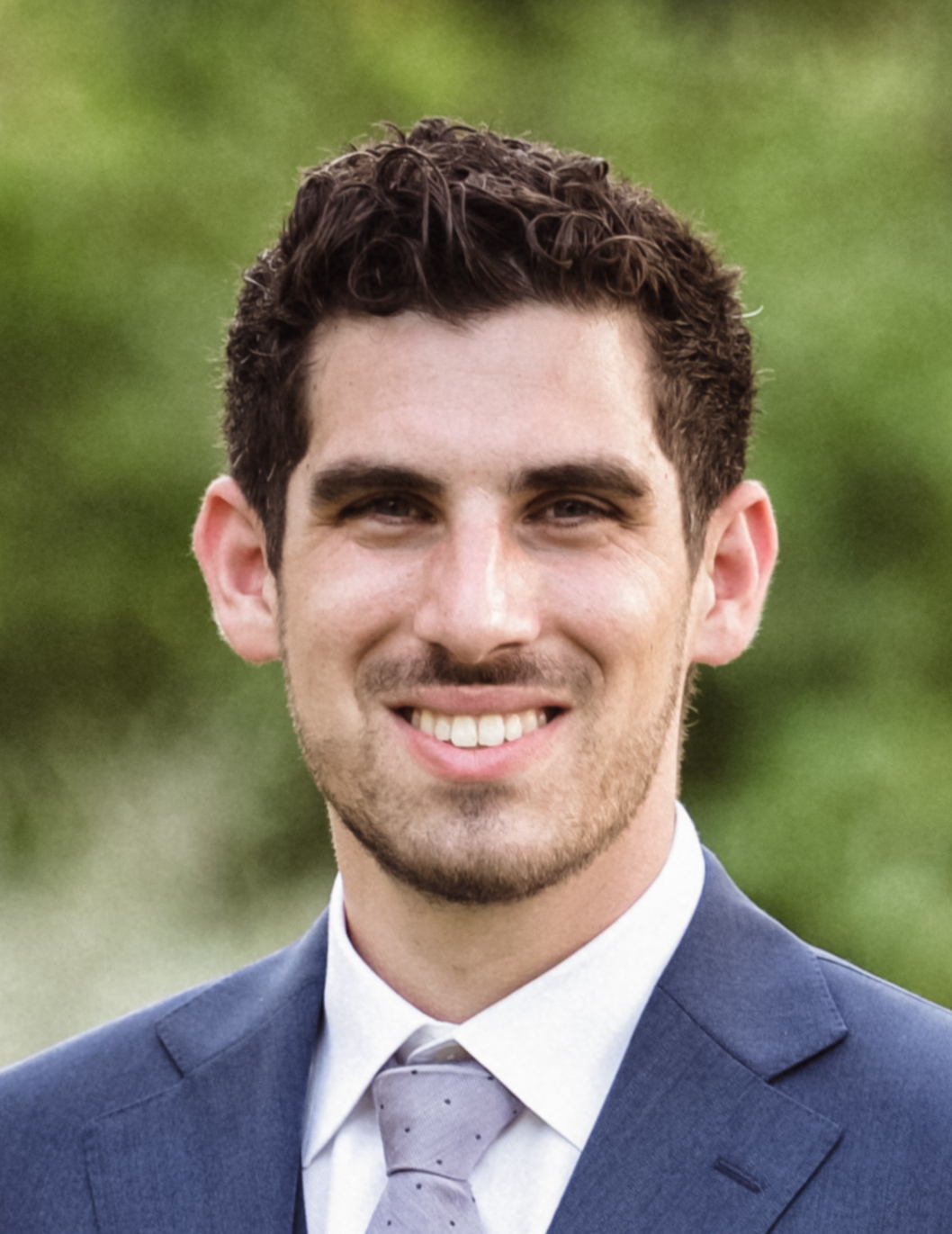}}]{Marcello D. Guarro} received his B.S. degree in Applied Physics in 2010, his M.S. degree in Computer Engineering in 2018, and PhD. in Computer Engineering in 2021 all from the University of California, Santa Cruz. 

His research interests include modeling, observer design, stability, control, and robustness analysis of hybrid systems with applications in cyber-physical systems, clock synchronization, autonomous vehicles, and robotics. 
\end{IEEEbiography}


\begin{IEEEbiography}[{\includegraphics[width=1in,height=1.25in,clip,keepaspectratio]{../Figures/Sanfelice5x7a_2012smallIEEE.eps}}]{Ricardo G. Sanfelice} received the B.S. degree in Electronics Engineering from the Universidad de Mar del Plata, Buenos Aires, Argentina, in 2001, and the M.S. and Ph.D. degrees in Electrical and Computer Engineering from the University of California, Santa Barbara, CA, USA, in 2004 and 2007, respectively. In 2007 and 2008, he held postdoctoral positions at the Laboratory for Information and Decision Systems at the Massachusetts Institute of Technology and at the Centre Automatique et Systèmes at the École de Mines de Paris. In 2009, he joined the faculty of the Department of Aerospace and Mechanical Engineering at the University of Arizona, Tucson, AZ, USA, where he was an Assistant Professor. In 2014, he joined the University of California, Santa Cruz, CA, USA, where he is currently Professor in the Department of Electrical and Computer Engineering. Prof. Sanfelice is the recipient of the 2013 SIAM Control and Systems Theory Prize, the National Science Foundation CAREER award, the Air Force Young Investigator Research Award, the 2010 IEEE Control Systems Magazine Outstanding Paper Award, and the 2020 Test-of-Time Award from the Hybrid Systems: Computation and Control Conference. His research interests are in modeling, stability, robust control, observer design, and simulation of nonlinear and hybrid systems with applications to power systems, aerospace, and biology.
\end{IEEEbiography}
}

\end{document}